\documentclass[11pt,a4paper]{amsart}
\usepackage[utf8]{inputenc}
\usepackage[english]{babel}
\usepackage[T1]{fontenc}
\usepackage{lmodern}
\usepackage[colorlinks=true,citecolor=blue,linkcolor=black]{hyperref}
\usepackage{graphicx}
\usepackage[left=3cm,right=3cm,top=2.5cm,bottom=2.5cm]{geometry}
\usepackage{enumerate}
\usepackage{bm}
\usepackage{array}

\usepackage{amsmath,mathtools,amssymb,extarrows,mathrsfs,amsthm}
\usepackage{tikz-cd}
\usepackage{tikz}
\usetikzlibrary{backgrounds}
\usetikzlibrary{calc}
\usetikzlibrary{hobby}
\usetikzlibrary{decorations.markings}
\usetikzlibrary {arrows.meta}
\usepackage{caption}
\usepackage{subcaption}
\usepackage{xcolor}
\usepackage{float}

\DeclareRobustCommand{\SkipTocEntry}[5]{}
\usepackage{xpatch}
\makeatletter   
\xpatchcmd{\@tocline}
{\hfil\hbox to\@pnumwidth{\@tocpagenum{#7}}\par}
{\ifnum#1<0\hfill\else\dotfill\fi\hbox to\@pnumwidth{\@tocpagenum{#7}}\par}
{}{}
\makeatother
\makeatletter
 \def\l@subsection{\@tocline{2}{0pt}{4pc}{6pc}{}}
\def\l@subsubsection{\@tocline{3}{0pt}{8pc}{8pc}{}}
 \makeatother

\newtheorem{defi}{Definition}[section]
\newtheorem{thm}[defi]{Theorem}
\newtheorem{conj}[defi]{Conjecture}
\newtheorem{prop}[defi]{Proposition}
\newtheorem{lemma}[defi]{Lemma}
\newtheorem{cor}[defi]{Corollary}

\theoremstyle{definition}
\newtheorem{rem}[defi]{Remark}
\newtheorem{assumption}[defi]{Assumption}

\newcommand{\C}{\mathbb{C}}

\newcommand{\PP}{\mathbb{P}}
\newcommand{\R}{\mathbb{R}}
\newcommand{\Z}{\mathbb{Z}}

\newcommand{\cI}{\mathcal{I}}

\newcommand{\cR}{\mathcal{R}}

\newcommand{\sslash}{\mathbin{/\mkern-6mu/}}
\newcommand{\cir}[1]{\langle #1 \rangle}

\newcommand{\defeq}{\vcentcolon=}

\newcommand{\To}{\longrightarrow}

\newcommand{\real}{\mathop{\mathrm{Re}}}

\newcommand{\id}{\mathrm{id}}

\newcommand{\Sto}{\mathrm{Sto}}
\newcommand{\Hom}{\mathrm{Hom}}

\newcommand{\Ical}{\mathcal{I}}
\DeclareMathOperator{\ram}{ram}
\DeclareMathOperator{\Aut}{Aut}
\DeclareMathOperator{\Irr}{Irr}
\DeclareMathOperator{\slope}{slope}
\DeclareMathOperator{\GrAut}{GrAut}

\title[A topological algorithm for the Fourier transform of Stokes data]{A topological algorithm for the Fourier transform of Stokes data at infinity}

\author{Jean Douçot}
\address{(J.D.) Group of Mathematical Physics, Faculty of Sciences, Universidade de Lisboa, Campo Grande, Edifício C6, PT-1749-016 Lisboa, Portugal\newline \emph{Current address:} Institute of Mathematics of the Romanian Academy, Calea Grivi\c{t}ei 21, 010702 Bucharest, Romania}
\email{jeandoucot@gmail.com}
\thanks{The research of J.D. was supported by FCiências.ID}

\author{Andreas Hohl}
\address{(A.H.) Université Paris Cité and Sorbonne Université, CNRS, IMJ-PRG, F-75013 Paris, France; KU Leuven, Departement Wiskunde, Celestijnenlaan 200B, B-3001 Leuven, Belgium\newline \emph{Current address:} Technische Universität Chemnitz, Fakultät für Mathematik, 09107 Chemnitz, Germany}
\email{andreas.hohl@math.tu-chemnitz.de}
\thanks{The research of A.H.\ was supported by the Deutsche Forschungsgemeinschaft (DFG, German Research Foundation), Projektnummer 465657531, and the grant G0B3123N from the Fonds voor Wetenschappelijk Onderzoek -- Vlaanderen (FWO, Research Foundation -- Flanders)}

\date{}


\begin{document}
	
\begin{abstract}
We give a topological description of the behaviour of Stokes matrices under the Fourier transform from infinity to infinity in a large number of cases of one level. This explicit, algorithmic statement is obtained by building on a recent result of T.\ Mochizuki about the Fourier transform of Stokes data of irregular connections on the Riemann sphere and by using the language of Stokes local systems due to P.\ Boalch. In particular, this induces explicit isomorphisms between wild character varieties, in a much larger range of examples than those for which such isomorphisms have previously been written down. We conjecture that these isomorphisms are compatible with the quasi-Hamiltonian structure on the wild character varieties.
\end{abstract}
	
\maketitle

\tableofcontents

\section{Introduction}
Stokes data are generalised monodromy data associated to meromorphic connections with irregular singularities. Considering them is necessary to obtain an irregular analogue of the famous Riemann--Hilbert correspondence. The latter gives, on any smooth algebraic variety, an equivalence between connections with regular singularities (the de Rham side) and monodromy data (the Betti side). Its irregular version, known as the Riemann--Hilbert--Birkhoff correspondence, provides an equivalence between possibly irregular connections and generalised monodromy data. These generalised monodromy data can be encoded in several different ways, one of them being the notion of \emph{Stokes local system} (\cite{Boa01,Boa21}), which is closely related to the wild fundamental group of \cite{MR91}. Via this correspondence, moduli spaces of meromorphic connections with prescribed behaviour at the singularities give rise to the so-called wild character varieties, the irregular analogues of the usual character varieties. These moduli spaces can also be seen, via the nonabelian Hodge correspondence, as moduli spaces of irregular Higgs bundles (the Dolbeault side), and possess very rich geometric (symplectic \cite{Boa01}, hyperkähler \cite{BB04}) structures, leading to the notion of (wild) nonabelian Hodge space \cite{Boa18}.

The Fourier--Laplace transform is an operation on meromorphic connections on the affine line. It plays an important role in many contexts: for example, it underlies many symmetries of isomonodromy systems such as the existence of different Lax pairs for several Painlevé-type equations \cite{Boa05, Boa08, Boa12, Yam16, SW04} and is an essential ingredient in Arinkin's extension \cite{Arinkin} of the Katz algorithm for rigid irregular connections.

In turn, the question of describing directly the action of the Fourier transform on the Betti side, i.e.\ on Stokes data, via the Riemann--Hilbert--Birkhoff correspondence, is natural and of great interest. In particular, a completely intrinsic topological description of the transformation of Stokes data via Fourier transform would be desirable.

This problem already has quite a long history. Explicit expressions for the Stokes data of the Fourier transform have so far been obtained mainly in two cases: the case of a regular singular system (see e.g.\ \cite{BJL81}, \cite[Chap.\ XII]{Mal91}, \cite[§2]{Boa15a} and \cite{DHMS20}) and the case of quadratic monomial exponents at infinity (also called the case of \emph{pure Gaussian type}, see e.g.\ \cite{Sab16}, \cite{Ho22}). Beyond these cases, it has proved difficult to determine the Fourier transform of Stokes data outside of some specific situations (see e.g.\ \cite{HK22} for some computations of Stokes matrices that include cases interchanged by the Fourier transform). For example, in  \cite{Boa15a}, in the so-called simply-laced case featuring one irregular singularity at infinity of order less than three together with regular singularities at finite distance, a symplectic isomorphism has been found between two wild character varieties that are known to be related by the Fourier transform, but it is unclear whether the isomorphism described in loc.\ cit.\ is actually the one induced by the Fourier transform.

On the other hand, there have been several studies aiming at tackling the general case, beginning with Malgrange's book \cite{Mal91}, which provides a complete answer, albeit hardly explicit. The problem has been reconsidered more recently by T.\ Mochizuki, starting with the article \cite{Moc10}, which investigates the general case with a different approach.
In a second recent article \cite{Moc21}, T.\ Mochizuki develops another framework to obtain a more explicit transformation rule for Stokes data under Fourier transform in full generality. The results are formulated using a new description of Stokes data, which the author of loc.~cit.\ calls \emph{Stokes shells}. Although the result obtained there is valid in great generality, understanding the implications of its statement and proof in terms of explicit Stokes matrices is not an obvious task.

In this article, we will -- in a large class of examples -- use the result of \cite{Moc21} to deduce a topological transformation rule for Stokes local systems, by making the connection between the latter and the notion of Stokes shells. Basically, what we do is the following: In his article, T.\ Mochizuki starts from the notions of Stokes filtrations and Stokes local systems (see \cite[Chapter~2]{Moc21} and in particular §2.4.2 for the reference to P.\ Boalch's work). He modifies the presentation of Stokes data and introduces his notion of ``Stokes shells'', which is more adapted to the general formulation of his final result. 
It is our aim in this paper to show how to understand the latter back in the world of Stokes local systems, to draw some concrete conclusions and to give explicit examples.

In particular, our work is not to be seen as an alternative approach to \cite{Moc21}, but deeply relies on the results proved by T.\ Mochizuki and interprets them in a language which we can directly use for explicit computations involving Stokes matrices.\\

\addtocontents{toc}{\SkipTocEntry}
\subsection*{The setting} 
We will consider a large class of examples, going much beyond the cases for which explicit formulas for the Stokes data of the Fourier transform have been known.

Recall that an algebraic connection on $\mathbb{A}^1(\mathbb{C})$ is the same as an analytic meromorphic connection on $\mathbb{P}^1(\mathbb{C})$ with a singular point only at $\infty$ (cf.\ e.g.\ \cite[Ch.~V]{Mal91}). By the formal classification of meromorphic connections (see e.g.\ \cite[Ch.~III]{Mal91} and the references therein), one associates to (the formalisation at $\infty$ of) such a connection a finite set of exponents (sometimes called exponential factors) at $\infty$, which are (possibly multivalued) meromorphic functions with a pole at $\infty$.

The situation we will be considering is the case of (analytic) meromorphic connections on the Riemann sphere where both the initial connection and its Fourier transform only have one singular point, the latter being at infinity. This happens when all the exponents at $\infty$ of the initial connection are of slope $>1$ (and, by the stationary phase formula, this is then the case for its Fourier transform as well).

In order to have a cleaner description of the algorithm, we will make the following further simplifying assumptions (see Assumption~\ref{assumption} in Section~\ref{sec:defData}):
\begin{itemize}
	\item The set of exponents at $\infty$ of the connection is of the form $\{az^k\mid a\in A\}$ for some fixed $k\in\mathbb{Q}_{>1}$ and some finite set $A\subset \C\setminus\{0\}$ (where $z$ is a local coordinate at $0$),
	\item $|a|=|b|$ for any $a,b\in A$, i.e.\ all parameters lie on the same circle centred at the origin.
\end{itemize}
Note that these assumptions are invariant under Fourier transformation: If a system satisfies all our hypotheses, its Fourier transform will do so, too.

Note also that the above assumptions seem to entail no true loss of generality from the weaker assumption that the set of exponents is of pure level $k$ (meaning that each exponent and each difference of two exponents has a leading term of order $k$). Indeed, up to admissible deformations it is always possible to move the coefficients of the polynomial exponents so as to satisfy these conditions. With some caution, we could thus probably consider the case of pure level $k$. We try, however, to keep it simple here and will not make this generalisation, since our aim is to give a description in the clearest possible way, which is already sufficiently involved in this slightly more restricted situation.\\

\addtocontents{toc}{\SkipTocEntry}
\subsection*{Main results} The main result of this work is an algorithm for computing the Fourier transform of Stokes data in the cases covered by our above assumptions. It will be based on the results of \cite{Moc21} and formulated in the language of Stokes local systems: If $(E,\nabla)$ is an algebraic connection on the affine line satisfying our assumptions, and if $(\widehat{E},\widehat{\nabla})$ is its Fourier transform, the algorithm allows us to determine the Stokes local system of $(\widehat{E},\widehat{\nabla})$ from of the one of $(E,\nabla)$ in a completely explicit, algorithmic and topological way. While we refer the reader to the main body of the article for the fully detailed statements, let us sketch the main ideas underlying the algorithm and its structure. \\

It is well-known that the formal data, i.e.\ the Levelt--Turrittin normal form, of $(\widehat{E},\widehat{\nabla})$ are determined by those of $(E,\nabla)$: This is expressed by the stationary phase formula  \cite{Mal91, Sab08, Fa09}, which states in particular that the exponents of $(\widehat{E},\widehat{\nabla})$ are related to those of $(\widehat{E},\widehat{\nabla})$ via a Legendre transform. 

Usually, the Legendre transform is viewed as providing a bijection between the sets of exponents before and after the Fourier transform. The crucial idea behind the algorithm is that the Legendre transform can actually be understood in a stronger sense. This relies on describing the formal data of the connections à la Deligne--Malgrange, in terms of the \emph{exponential local system} $\mathcal I$. In brief, $\mathcal I$ is an infinite collection of circles covering the circle $\partial$ of directions around infinity. These circles represent all the possible exponents of a connection at $\infty$ (see §\ref{sec:explocsys} for a review of the exponential local system).
The stronger version of the Legendre transform is then the following\\

\noindent \textsl{Crucial idea:} The Legendre transform provides a self-homeomorphism $\ell$ of the collection $\Ical_{>1}$ of all circles in $\mathcal{I}$ of slope $>1$ at infinity (see \cite{Mal91,Dou21} and §\ref{sec:Legendre}).\\

Now, the idea of the algorithm is, roughly speaking, that it is possible to use $\ell$ to transport the Stokes data of  $(E,\nabla)$ and obtain those of $(\widehat{E},\widehat{\nabla})$.\\

In order to give a precise formulation of this idea, one needs to describe the formal data and the Stokes data of connections in a suitable way. The description we use is as follows:
\begin{itemize}
\item The formal data of $(E,\nabla)$ are encoded by a \emph{formal local system}, which is a local system (of finite-dimensional vector spaces) $V^0$ on $\partial$ graded by $\mathcal I$ (see §\ref{sec:formalLS}).
\begin{itemize}
    \item Those circles in $\mathcal{I}$ such that the corresponding graded parts of $V^0$ do not vanish are called the \emph{active circles} or \emph{active exponents} of $(E,\nabla)$. The active circles together with the ranks of the corresponding graded pieces of $V^0$ are called the \emph{irregular class} $\Theta$ of $(E,\nabla)$.
    \item The formal local system $V^0$ can equivalently be thought of as a local system on the topological space $\mathcal{I}$ supported on its active circles (see Theorem~\ref{thm:gradedLSvsLSonI}).
    \item The irregular class determines an auxiliary surface $\widetilde{\Sigma}(\Theta)$ (obtained from $\widehat{\Sigma}$ by puncturing at some points), a \textit{fission surface} $\underline{\widetilde\Sigma}(\Theta)$ (see Remark~\ref{rem:fissionLS}) and a finite collection of distinguished paths $\gamma_{i\to j}$ in the fission surface, the so-called \emph{Stokes paths}, whose endpoints $i,j$ are points on some active circles of $\Theta$ (see Remark~\ref{rem:defdataPath}).
\end{itemize}

\item The full Stokes data of $(E,\nabla)$ are encoded in terms of a \textit{Stokes local system} $\mathbb{V}$, which is a local system (of finite-dimensional vector spaces) on $\widetilde{\Sigma}(\Theta)$ satisfying certain conditions (see §\ref{sec:SLSandFission}).
\begin{itemize}
\item The Stokes local system $\mathbb{V}$ of $(E,\nabla)$ determines in particular the associated formal local system $V^0$, which lives on the boundary circle $\partial\subset\widetilde{\Sigma}(\Theta)$.
\item A Stokes local system $\mathbb{V}$ can equivalently be thought of as a local system $\underline{\mathbb{V}}$ on the fission surface $\underline{\widetilde{\Sigma}}(\Theta)$ satisfying certain conditions (see Remark~\ref{rem:fissionLS} for this viewpoint).
\item The Stokes local system $\mathbb V$ can be reconstructed up to isomorphism from its associated formal local system and the data of the parallel transport $\rho_{i\to j}$ of $\underline{\mathbb{V}}$ along the Stokes paths $\gamma_{i\to j}$ in the fission surface (see Remark~\ref{rem:defdataPath} and Proposition~\ref{prop:SLSdefdata} for a more precise statement).

(The linear maps $\rho_{i\to j}$ basically correspond to the nontrivial entries of the Stokes matrices.)
\end{itemize}
\end{itemize}

In this language, our algorithmic description of the Fourier transform of Stokes data, based on \cite{Moc21}, can be formulated as follows: 

\begin{thm}[see Theorem~\ref{thm:trafoRule} for the fully detailed statement]\label{thm:algoPaths}
Let $(E,\nabla)$ be an algebraic connection on the affine line satisfying our assumptions, $\Theta$ its irregular class, and $\underline{\mathbb V}$ its fission local system. Let $(\widehat E,\widehat \nabla)$ be the Fourier transform of $(E,\nabla)$, and $\widehat V^0$, $\widehat{\Theta}$, $\underline{\mathbb V}$ its formal local system, irregular class and fission local system, respectively. Write $k=\frac{s}{r}$ for $s,r\in \mathbb{Z}_{>0}$ coprime, where $k\in\mathbb{Q}_{>1}$ is as in our assumptions. Then we have:
\begin{itemize}
\item The formal local systems on both sides are related by the Legendre transform, up to a possible sign for the formal monodromy, i.e.
\[
\widehat{V}^0\cong\ell_* V^0\otimes W,
\]
where $W$ is a rank one local system on $\partial$ with monodromy $(-1)^{s}$.  
\item The fission local system $\underline{\widehat{\mathbb V}}$ is determined as follows: For any Stokes path $\widehat\gamma_{\hat i\to \hat j}$ in the fission surface $\underline{\widetilde{\Sigma}}(\widehat\Theta)$, the corresponding parallel transport $\widehat\rho_{\hat i\to \hat j}$ is obtained, up to a sign (which is explicitly determined), by pulling back the endpoints $\hat{i}, \hat{j}$ to the initial fission surface $\underline{\widetilde{\Sigma}}(\Theta)$ using the inverse of the Legendre transform and taking the parallel transport in the fission local system $\underline{\mathbb V}$ along the (natural) path from $\ell^{-1}(\hat{i})$ to $\ell^{-1}(\hat{j})$, i.e.
\[
\widehat{\rho}_{\hat{i}\to\hat{j}}=\pm \rho_{\ell^{-1}(\hat{i})\to \ell^{-1}(\hat{j})}.
\] 
\end{itemize}
\end{thm}

Let us note that one has to carefully keep track of the signs in the formal local system here in order to correctly determine $\underline{\widehat{\mathbb V}}$ from the parallel transports along the Stokes paths. In particular, the sign $(-1)^{s}$ comes by gluing the graded pieces of the formal local system from several intervals and introducing multiple signs in this gluing process, which is made precise in Section~\ref{sec:alg}.
 
In particular, this allows us to compute the Stokes matrices of $(\widehat E,\widehat \nabla)$ as a function of the Stokes matrices of $(E,\nabla)$.
In turn, the algorithm induces an explicit isomorphism between the wild character varieties on both sides of the Fourier transform. Given an irregular class $\Theta$, we denote by $\mathcal{E}_{\Theta}$ the associated \emph{reduced representation variety}, parametrising Stokes local systems with irregular class $\Theta$ and a minimal framing (see Definition~\ref{def:reducedRepVar}).

\begin{thm}[see Theorem~\ref{thm:algebraic}]\label{thm:isomWCV}
Let $\Theta$ be an irregular class satisfying the above assumptions and let $\widehat{\Theta}$ be its formal Fourier transform. The algorithm induces an algebraic isomorphism between the reduced representation varieties $\mathcal E_{\Theta}$ and $\mathcal E_{\widehat\Theta}$, provided that they are nonempty.
\end{thm}

Remarkably, the Fourier transform is thus algebraic both on the de Rham side and on the Betti side. This is not obvious since the Riemann--Hilbert--Birkhoff correspondence only gives an analytic isomorphism between both sides in general.\\

\addtocontents{toc}{\SkipTocEntry}
\subsection*{Outlook} Several further questions arise from the results of this article.

For example, we conjecture that the isomorphisms from Theorem~\ref{thm:isomWCV} are compatible with the quasi-Hamiltonian structures of the wild character varieties on both sides of the Fourier transform. This symplectic nature of the Fourier transform is already known in a few particular cases \cite{Boa15a,Sza15}, and this is consistent with the broader perspective of `global Lie theory' taken in \cite{Boa18}, where one views the moduli spaces on both sides of the Fourier transform as two different `representations' of the same abstract nonabelian Hodge space.

Another interesting point is the following: The isomorphisms provided by the algorithm are entirely topological and involve some very nontrivial combinatorics. While they give a very explicit way of computing Stokes data for the Fourier transform, it would be desirable to understand directly (without passing through the Riemann--Hilbert--Birkhoff correspondence and the result of \cite{Moc21}) why they are well-defined, i.e.\ why they produce an element in the correct moduli space of Stokes local systems. It seems possible that some of the combinatorial structures which appear in the study of related spaces, such as braid varieties, and their cluster and symplectic structures \cite{STZ19,CGGS20, GK21} might be used to shed some light on this question, or to understand why the isomorphisms are quasi-Hamiltonian.

Finally, we would like to understand if a similar topological way of computing the Stokes structure for the Fourier transform exists for more general situations. This concerns in particular the case of connections with singularities at finite distance. These singularities will contribute to the Stokes structure of the Fourier transform at infinity. Formally, this is again well-understood in terms of Legendre transforms, and it is natural to ask for a description of the complete Stokes data in terms of Stokes local systems and paths similar to the above.

We hope to address these questions in future work.\\

\addtocontents{toc}{\SkipTocEntry}
\subsection*{Structure of the article} The article is organised as follows. We start by recalling a few facts about the formal data of irregular connections in Section~\ref{sec:formal}, as well as about the description of Stokes data in terms of Stokes local systems, leading to the construction of wild character varieties, in Section~\ref{sec:SLS}.

We then review the description of Stokes data in terms of deformation data used by T.\ Mochizuki in \cite{Moc21} and relate it to Stokes local systems in Section~\ref{sec:defData}, via the observation that to each Stokes arrow in the Stokes diagram corresponds one deformation datum. The most technical work happens in §\ref{sec:defDataAndSplittings}, but it is not necessary to read these details in order to understand the final statement. In Section~\ref{sec:alg}, we formulate an algorithm based on Mochizuki's results which yields the Stokes local system of the Fourier transform starting from the initial Stokes local system, thus inducing an isomorphism between the corresponding wild character varieties.

Finally, in Section~\ref{sec:examples}, we use this algorithm to compute explicitly in several examples the Stokes matrices of the Fourier transform and determine the corresponding isomorphisms, which in all cases are symplectic, strongly suggesting that this will be true in the general case.

\addtocontents{toc}{\SkipTocEntry}
\subsection*{Acknowledgements} We thank Philip Boalch for many discussions that helped us shape the main statements of this work in the language of Stokes local systems. We are indebted to Takuro Mochizuki for answering several questions on his work. We are also grateful to Claude Sabbah, Marco Hien and Andrea D'Agnolo for valuable meetings and seminars about Fourier transforms and, in particular, the article \cite{Moc21}.

\section{Formal data of irregular connections}\label{sec:formal}

The classification of irregular connections up to formal gauge transformations is provided by the well-known Levelt--Turrittin theorem. It will be convenient for us to use a formulation of the formal classification in terms of graded local systems, which we briefly review now, following mainly \cite{Boa21,BY15}.

\subsection{The exponential local system}\label{sec:explocsys}
The idea is to view the exponents as sections of a local system on the circle of directions around the singularity. Let $\Sigma=\mathbb P^1$. Let $\varpi\colon \widehat{\Sigma}\to \Sigma$ the real oriented blow-up at $\infty$ of $\Sigma$. The preimage $\partial\defeq\varpi^{-1}(\infty)$ is a circle whose points correspond to the directions around $\infty$. An open subset of $\partial$ corresponds to a germ of open sector at infinity.

Let $z$ be the standard complex coordinate on $\mathbb C$, so that $1/z$ provides a local coordinate of $\mathbb P^1$ at infinity. The exponential local system $\Ical$ is a local system of sets (i.e.\ a covering space) on $\partial$ whose sections on open intervals in $\partial$ are holomorphic functions of the form
\[
\sum_{i} a_i z^{k_i}
\]
on the corresponding germs of open sectors, where $k_i\in \mathbb{Q}_{>0}$, and $a_i\in\C$.
Each connected component of (the étalé space of) such a local section is a finite order cover of the circle $\partial$. More precisely, let $r$ be the smallest integer such that the expression $q=\sum_{i} a_i z^{k_i}$ is a polynomial in $z^{1/r}$. The corresponding holomorphic function is multivalued, and becomes single-valued when passing to a finite covering $t^r=z$. Therefore, the corresponding connected component of $\cI$, which we denote by $\cir{q}$, is an $r$-sheeted covering of $\partial$. As a topological space, it is homeomorphic to a circle, and $\Ical$ is thus a disjoint union of (an infinite number of) circles. Denote the projection to the base circle by $\pi\colon \mathcal{I}\to \partial$.

More concretely, for $q=\sum_{i} a_i z^{k_i}$ as above, a point on the circle $\langle q\rangle$ is given by a pair $(\theta,[\tilde q])$, where $\theta\in\partial$ is a point and $[\tilde q]$ is the germ on a small sector around the direction $\theta$ of a holomorphic function $\tilde q=\sum_{i} a_i z^{k_i}$, given by a choice of determination of $z^{1/r}$. We denote by $\pi_q\colon \langle q\rangle \to \partial$ the corresponding covering (the restriction of $\pi$ to $\langle q\rangle$). It is given by $(\theta,[\tilde q])\mapsto \theta$.

The integer $r=\vcentcolon\ram(q)$ is the \textit{ramification order} of the circle $\cir{q}$. If $r=1$ we say that the circle $\cir{q}$ is unramified. The degree $s$ of $q$ as a polynomial in $z^{1/r}$ is the \textit{irregularity} $s=\vcentcolon\Irr(\cir{q})$ of $\cir{q}$ and the quotient $s/r=\vcentcolon\slope(I)$ is the \textit{slope} of $\cir{q}$.

If $d\in\partial$, we denote by $\Ical_d$ the fibre of $\Ical$ over the direction $d$. Taking $d$ as a base point of $\partial$, the monodromy of $\Ical$ is the automorphism
\[ \rho\colon \Ical_d \to \Ical_d
\]
of the fibre $\Ical_d$ obtained when going once around $\partial$ in the positive direction.

\subsection{Formal local systems}\label{sec:formalLS}

\begin{defi}
An $\Ical$-graded local system on $\partial$ is a local system $V^0\to \partial$  of finite-dimensional vector spaces, together with a grading of each fibre $V^0_d$, with $d\in\partial$, by the set $\Ical_d$, in such a way that this grading is locally constant. This means that for any $d\in \partial$, there is a direct sum decomposition
\[ V^0_d =\bigoplus_{i\in \Ical_d} V^0_{d}(i)
\]
such that for any path $\gamma$ from $d$ to $d'$ in $\partial$ (which naturally induces actions $\rho_\gamma\colon V^0_d\to V^0_{d'}$ and $\rho_\gamma\colon \cI_d\to\cI_{d'}$), one has
$$\rho_\gamma(V^0_d(i)) = V^0_q(\rho_\gamma(i)).$$
(Note that this definition works completely analogously if we replace the exponential local system $\mathcal{I}$ by any local system of sets on $\partial$.)
\end{defi}

For a fixed $d\in \partial$, let us set $n_i\defeq\dim V^0_{d}(i)$ for $i\in \Ical_d$. Since $V^0$ is finite-dimensional, there is only a finite number of elements $i\in \Ical_d$ for which $n_i$ is nonzero. Furthermore, since the grading is locally constant, $n_i$ only depends on the connected component $\cir{q}$ of $\Ical$ containing $i$, and we set $n_{\cir{q}}\defeq n_i$.

\begin{defi}
An irregular class is a function $\Theta\colon \pi_0(\Ical)\to \mathbb N$ with finite support.
\end{defi}

The dimensions $n_i$ of the graded pieces thus define an irregular class $\Theta$ associated to $V^0$. The connected components $\cir{q}\in \pi_0(\Ical)$ such that $n_{\cir{q}}\neq 0$ are called the \emph{active circles} or \emph{active exponents}. If $i\in \mathcal{I}_d$ and $\cir{q}$ is the circle with $i\in \cir{q}$, we sometimes write $\Theta(i)\defeq \Theta(\cir{q})$ (this is the above number $n_i$).

Let us again fix a base direction $d\in \partial$ and denote by $\rho\in \Aut(V_d)$ the monodromy
\[\rho\colon V^0_d\to V^0_d.
\]
Since the grading is locally constant, the monodromy of $V$ has to be compatible with the monodromy of $\Ical$, that is
\[ \rho(V^0_d(i))=V^0_d(\rho (i)),
\]
for any $i\in\Ical_d$. In particular, the $\cI$-graded local system $V^0$ globally splits as
$$V^0=\bigoplus_{\cir{q}\in \pi_0(\cI)} V^0_{\cir{q}},$$
where $V^0_{\cir{q}}$ is a $\cir{q}$-graded local system and $V^0_{\cir{q}}\neq 0$ if and only if $\cir{q}$ is an active circle.

\begin{figure}[H]
\centering
\begin{tikzpicture}[scale=0.6]
\draw[domain=71.44:(3*360-71.44),scale=1,samples=1000] plot ({-cos(\x)},{-sin(\x)/2+\x/360});
\draw[postaction={decorate}] ({-cos(71.44)},{-sin(71.44)/2+71.44/360}) to [out angle=180, in angle=180,curve through={(-1.5, 3/2)}]({-cos(3*360-71.44)},{-sin(3*360-71.44)/2+(3*360-71.44)/360});
\draw[domain=0:(360),scale=1,samples=1000] plot ({-cos(\x)},{-sin(\x)/2-3});
\draw[->] (0,-1)--(0,-2);
\draw (-2, 3/2) node{$\langle q\rangle$};
\draw (-2,-3) node{$\partial$};
\draw (0.7, -1.5) node{$\pi$};
\node (A1) at (1,0.5) {};
\node (A2) at (1,1.5) {};
\node (A3) at (1,2.5) {};
\draw[fill=black] (A1) circle (0.07cm);
\draw[fill=black] (A2) circle (0.07cm);
\draw[fill=black] (A3) circle (0.07cm);
\draw[fill=black] (1,-3) circle (0.07cm);
\draw (1.4,-3) node {$d$};
\draw ($(A1)+(0.4,0.5)$)--($(A1)+(-0.4,0.3)$)--($(A1)+(-0.4,-0.5)$)--($(A1)+(0.4,-0.3)$)--cycle;
\draw ($(A2)+(0.4,0.5)$)--($(A2)+(-0.4,0.3)$)--($(A2)+(-0.4,-0.5)$)--($(A2)+(0.4,-0.3)$)--cycle;
\draw ($(A3)+(0.4,0.5)$)--($(A3)+(-0.4,0.3)$)--($(A3)+(-0.4,-0.5)$)--($(A3)+(0.4,-0.3)$)--cycle;
\draw ($(A1)+(1.5,0)$) node {$V^0_{d}(k_0)$};
\draw ($(A2)+(1.5,0)$) node {$V^0_{d}(k_1)$};
\draw ($(A3)+(1.5,0)$) node {$V^0_{d}(k_2)$};
\end{tikzpicture}
\caption{A local system $V^0\to\partial$ graded by a single Stokes circle $\cir{q}$, here with ramification order $r=3$. For any direction $d\in \partial$, the fibre $\cir{q}_d=\pi^{-1}(d)=\{k_0, \dots, k_{r-1}\}$ contains $r$ points, and we have a locally constant decomposition $V^0_d=\bigoplus_{i=0}^{r-1}V^0_d(k_i)$. It will often be convenient to view $V^0$ as a local system on the circle $\cir{q}$, in such a way that its fibre over the point $k_i\in \cir{q}$ is the graded piece $V^0_d(k_i)$.}
\label{fig:graded_local_system}
\end{figure}

More explicitly, let $\cir{q}$ be a circle in $\Ical$ with ramification $r$. Fix $d\in\partial$ and choose a numbering on the fibre $\cir{q}_d=\{i_0,\dots,i_{r-1}\}$ such that $\rho(i_k)=i_{k+1}$ for $i=0,\dots,r-1$ (with $i_r\defeq i_0$). Then the monodromy of the piece $V^0_{\cir{q},d}\defeq\bigoplus_{k=0}^{r-1} V^0_d(i_k)$ has the form

\begin{equation}
	\rho_{\cir{q}}=\begin{pmatrix}
		0 & \hdots & \hdots &  0 & \rho_{r-1,0} \\
		\rho_{0,1} & \ddots & & \vdots & 0\\
		0 & \rho_{1,2} & \ddots & \vdots & \vdots \\
		\vdots & \ddots & \ddots & 0 & \vdots\\
		0 & \hdots & 0 & \rho_{r-2,r-1} & 0
	\end{pmatrix}\in \Aut(V^0_{\cir{q}})
	\label{form_formal_monodromy}
\end{equation}
with $\rho_{i_k,i_{k+1}}\colon V^0_d(i_k)\to V^0_d(i_{k+1})$.

The group $\GrAut(V^0_d)$ of graded automorphisms of $V^0_d$ is
\[ 
\GrAut(V^0_d)=\prod_{i\in \Ical_d} \mathrm{GL}(V^0_d(i)).
\]
This is illustrated in Figure~\ref{fig:graded_local_system}.

An $\Ical$-graded local system $V^0$ is entirely determined up to isomorphism by the data of its irregular class $\Theta$ and the equivalence class $\mathcal{C}$ of its monodromy under graded automorphisms (acting by conjugation). We will refer to the pair $(\Theta, \mathcal C)$ as the formal data of $V^0$.

This language enables us to have a more geometric formulation of the Levelt--Turrittin decomposition:

\begin{thm}[{see \cite[Théorème~IV.2.3]{Mal91}}]
The category of connections on the formal punctured disk is equivalent to the category of $\Ical$-graded local systems on $\partial$.
\end{thm}

In turn, to any algebraic connection $(E,\nabla)$ on the affine line is associated (by considering the corresponding connection on the formal punctured disk at $\infty$) a \emph{formal local system} $V^0\to \partial$ (which is $\cI$-graded). The active circles correspond to the exponents of the connection in the Levelt--Turrittin normal form, and the monodromy of $V^0$ is usually called the \emph{formal monodromy}. In particular, a connection is regular if and only if its only active circle is the tame circle $\cir{0}$. 

\medskip

A useful observation is that one can equivalently view an $\Ical$-graded local system $V^0\to \partial$ as a local system on (the étalé space of) $\Ical$, i.e.\ there exists a local system on $\widetilde{V}^0$ on $\mathcal{I}$ such that $V^0\cong \pi_*\widetilde{V}^0$.
If $i\in\Ical$, the fibre of $\widetilde{V}^0$ over $i$ is $\widetilde{V}^0_i\defeq V^0_d(i)$, where $d=\pi(i)\in\partial$. If $\cir{q}$ is a circle of ramification $r$, when going once around $\langle q\rangle$ we go $r$ times around $\partial$. As a consequence, the monodromy of the $\langle q\rangle$-graded local system $V^0_{\cir{q}}$ , as described in \eqref{form_formal_monodromy}, is related to the monodromy $\rho_{\cir{q}}$ of $\widetilde{V}^0_{\cir{q}}$ (now seen as a local system over the circle $\langle q\rangle$) via
\[
\rho_{\cir{q}}=\rho_{r-1,0}\circ \rho_{r-2,r-1}\circ \dots \circ \rho_{0,1} \colon V^0_d(i)\to V^0_d(i),
\]
with $i_0=i$. By abuse of notation, we will also denote $\widetilde{V}^0$ just by $V^0\to\mathcal{I}$, and we will sometimes freely switch between the two perspectives (see again Figure~\ref{fig:graded_local_system} for an illustration).

We thus have (see \cite{BY20} for this viewpoint):

\begin{thm}\label{thm:gradedLSvsLSonI}
The category of connections on the formal punctured disk is equivalent to the category of local systems on $\Ical$ with compact support.
\end{thm}

Finally, notice that the notion of local system on $\partial$ that we will be using here  to describe the formal data of connections is closely related to the one of $2\pi\Z$-equivariant local system on $\R$ which is used in \cite{Moc21}:

\begin{lemma}
	A local system on $S^1$ is the same as a $2\pi\Z$-equivariant local system on $\R$.
\end{lemma}
\begin{proof}
	If $V^0$ is a local system on $S^1$ and $\tau\colon \R\to S^1, x\mapsto e^{ix}$ is the universal covering, then $\tau^{-1}V^0$ is a $2\pi\Z$-equivariant local system on $\mathbb R$. The inverse construction is also easily described. 
\end{proof}
Let $\langle q\rangle\subset \cI$ be a circle in the the exponential local system, and denote by $\pi_q\colon \langle q\rangle \to \partial$ the projection as above. Then, under the correspondence just explained, a $\langle q\rangle$-graded local system $V^0$ on $\partial$ corresponds to a $2\pi\Z$-equivariant local system $L$ on $\R$ such that each stalk is graded as
$$L_x=\bigoplus_{i\in \pi_q^{-1}(\tau(x))} V^0_{x}(i)$$
and the isomorphism $L_x\cong L_{x+2\pi}$ given by the $2\pi\Z$-invariance is a graded isomorphism induced by the permutation that comes from the monodromy of $V^0$ at $\pi_q^{-1}(\tau(x))$.

\section{Stokes local systems and wild character varieties}\label{sec:SLS}

The essence of the Stokes phenomenon is that given an irregular connection $(E,\nabla)$, there exist isomorphisms between the local system $V$ of its analytic solutions and the corresponding formal local system $V^0$, but these isomorphisms are only valid on sectors with limited angular width around the singularities. To have an analogue of the Riemann--Hilbert correspondence, which for connections with regular singularities provides an equivalence between connections and their monodromy data, in the case of irregular singularities one has to add extra data which account for passing from one sector to the next, known as Stokes data. In the literature, one finds many different ways of describing Stokes data. In this work, we will mostly use the formulation of \cite{Boa21} in terms of Stokes local systems, which has the advantage of providing intrinsic explicit presentations of moduli spaces of Stokes data. We briefly review it in this section, and we will relate it to another description used in \cite{Moc21}, namely in terms of `deformation data', in the next section.

We will present everything in our case of interest, where $\Sigma=\PP^1(\C)$ is the Riemann sphere, and the only marked point is at $\infty$. However, the generalisation to an arbitrary Riemann surface with one or several marked points is not difficult.

\subsection{Singular directions, Stokes arrows and Stokes groups}

Let $(E,\nabla)$ be a connection on the affine line and $V^0\to \Ical$ the corresponding formal local system, with irregular class $\Theta$. Let $\mathbb{I}\subset \Ical$ be the finite subcover of $\partial$ consisting of the active circles of $V^0$. 

\begin{defi}
Let $d\in\partial$. We define a partial order $\prec_d$ on $\mathbb{I}_d$ in the following way. Let $i,j\in \mathbb{I}_d$, and let $q_i,q_j$ be the corresponding germs of holomorphic functions.  We say that $i\prec_d j$ if the exponential $e^{q_i-q_j}$ decays fastest in the direction $d$ when $z$ tends to $\infty$. This happens if the leading term of the difference $(q_i-q_j)(z)$ is in $\R_{<0}$ when $z$ has direction (i.e.\ argument) $d$ and $|z|\gg 1$.

In this case, we will say that there is a \emph{Stokes arrow} from $j$ to $i$ over $d$, and we also write $j\to_d i$ instead of $i\prec_d j$.
\end{defi}

\begin{defi}
Let $d\in\partial$. If there exist $i,j\in \mathbb{I}_d$ such that $i\prec_d j$, then $d$ is called a \emph{singular direction} or \emph{anti-Stokes direction} for $(E,\nabla)$. We denote by $\mathbb{A}\subset\partial$ the set of singular directions.
\end{defi}

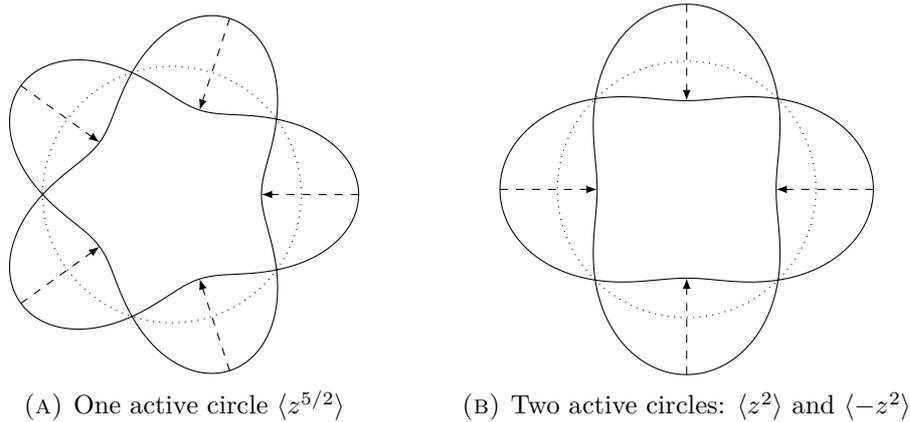
\begin{figure}
\centering

\begin{subfigure}[b]{0.49\textwidth}
\centering
\begin{tikzpicture}[scale=1.7]
\draw[dotted] (0,0) circle (1);
\draw[domain=0:(2*360),scale=1,samples=500] plot (\x:{exp(cos(-5/2*\x)/exp(1))});
\foreach \x in {144,288,...,720} \draw[dashed,->,>=latex] (\x:{exp(cos(-5/2*\x)/exp(1))})--(\x:{exp(-cos(-5/2*\x)/exp(1))});
\end{tikzpicture}
\caption{One active circle $\cir{z^{5/2}}$}
\end{subfigure}
\hspace{-1cm}
\begin{subfigure}[b]{0.49\textwidth}
\centering
\begin{tikzpicture}[scale=1.7]
  \draw[dotted] (0,0) circle (1);
  \draw[domain=0:360,scale=1,samples=500] plot (\x:{exp(cos(-2*\x)/exp(1))});
  \draw[domain=0:360,scale=1,samples=500] plot (\x:{exp(cos(-2*\x+180)/exp(1))});
   \foreach \x in {0,180} \draw[dashed,->,>=latex] (\x:{exp(cos(-2*\x)/exp(1))})--(\x:{exp(cos(-2*\x+180)/exp(1))});
   \foreach \x in {90,270} \draw[dashed,->,>=latex] (\x:{exp(cos(-2*\x+180)/exp(1))})--(\x:{exp(cos(-2*\x)/exp(1))});
\end{tikzpicture}
\caption{Two active circles: $\cir{z^2}$ and $\cir{-z^2}$}
\end{subfigure}
\caption{Stokes diagrams for two examples of irregular classes. The diagram represents the growth rate of the exponents as a function of the direction around the singularity. The dotted circle separates the regions where the exponential is growing or decreasing near $\infty$. The singular directions are the directions where we find a maximal distance between two strands, and the corresponding Stokes arrows are drawn. The \emph{Stokes directions} are the directions where some of the strands of the diagram cross.}
\label{fig:examples_of_Stokes_diagrams}
\end{figure}

When the exponents are monomials of the same degree, their relative dominance order can be nicely visualised in the \textit{Stokes diagram} of the connection (see Figure~\ref{fig:examples_of_Stokes_diagrams}), which is obtained by plotting $|e^{q}|=e^{\real(q)}$ as a function of $d\in\partial$ (for $|z|\gg 1$) for each active circle $\cir{q}$.

\begin{defi}
Let $d\in \mathbb{A}$ be a singular direction. The Stokes group associated to $d$ is the unipotent subgroup $\Sto_d\subset \mathrm{GL}(V^0_d)$ whose Lie algebra is given by
\[ \mathrm{sto}_d\defeq\bigoplus_{i\prec_d j} \Hom(V^0_d(j),V^0_d(i)).
\]
\end{defi}

More concretely, an element of $\mathrm{Sto}_d$ is a block matrix, with the block structure corresponding to the dimensions of the graded pieces $V^0_d(i)$, such that its diagonal blocks are identity matrices, for any Stokes arrow $i\prec_d j$ the corresponding block at position $(i,j)$ is an arbitrary matrix, and all other blocks are zero.

\subsection{Stokes local systems and fission local systems}\label{sec:SLSandFission}

One possible way of encoding Stokes data is given by the notion of \textit{Stokes local systems} \cite{BY15,Boa21}, which we now recall. In our case with $\Sigma=\mathbb{P}^1$, given an irregular class $\Theta$ at $\infty$, the idea is to define a slightly modified surface by introducing a puncture near each singular direction in the real blow-up $\widehat{\Sigma}$ of the Riemann sphere at $\infty$. A Stokes local system is then a local system on this new surface which near $\partial$ is $\cI$-graded (with irregular class $\Theta$) and such that the monodromies around the new, so-called \emph{tangential punctures} are elements in the corresponding Stokes groups.

More precisely, we define the new surface $\widetilde{\Sigma}(\Theta)$ as follows. We define the halo $\mathbb{H}\subset \widehat{\Sigma}$ as a tubular neighbourhood of $\partial$. One of the boundaries of $\mathbb{H}$ is the circle $\partial$. Let $\partial'$ be the other boundary of $\mathbb{H}$. Let $e\colon\partial\to \partial'$ be a homeomorphism preserving the orientation.

\begin{figure}
\centering
\begin{tikzpicture}
  \draw[loosely dotted] (0,0) circle (2);
    \draw (0,0) circle (1);
    \node at (1.065,0.5) {$\partial$};
    \node at (2.025,-1) {$\partial'$};
  \foreach \x in {72,144,...,360} \draw[fill=white] (\x:2) circle (0.1);
\end{tikzpicture}
\caption{Local picture: The halo and the tangential punctures at a singularity.}
\label{fig:halo}
\end{figure}
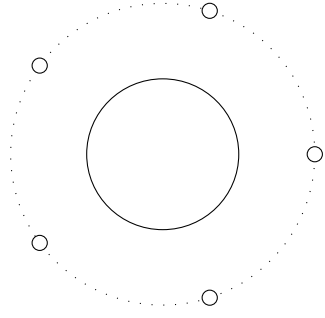

The irregular class $\Theta$ defines a set $\mathbb{A}\subset \partial$ of singular directions. We define a new surface as
\[
\widetilde{\Sigma}\defeq\widetilde{\Sigma}(\Theta)\defeq\widehat{\Sigma}\setminus e(\mathbb{A}),
\]
that is, for each singular direction $d\in \mathbb{A}$, we remove from $\widehat{\Sigma}$ the corresponding \textit{tangential puncture} $e(d)$ (see Figure~\ref{fig:halo}). For each $d\in \mathbb{A}$, let us then denote by $\gamma_d$ a small positive loop in $\widetilde{\Sigma}$ starting from $d$, going around the tangential puncture $e(d)$ (and no other puncture) in a positive sense and going back to $d$.

\begin{defi}
A Stokes local system is a pair $(\mathbb{V},\Theta)$ where $\Theta$ is an irregular class at infinity and
$\mathbb{V}$ is a local system of vector spaces on $\widetilde{\Sigma}(\Theta)$ equipped with an $\mathbb{I}$-grading (of
dimension $\Theta$) over the halo $\mathbb{H}$, where $\mathbb{I}\subset \mathcal I$ denotes the finite subcover corresponding to the active circles of $\Theta$, such that the monodromy $\mathbb{S}_d \defeq \rho(\gamma_
d)$ is in $\Sto_d \subset \mathrm{GL}(\mathbb{V}_d)$ for each $d\in \mathbb{A}$.
\end{defi}

Stokes local systems yield a topological description of the category of algebraic connections on the affine line $\mathbb C=\mathbb P^1\setminus \{\infty\}$ (cf.\ e.g.\ \cite{Boa21}):

\begin{thm}
The category of algebraic connections on $\mathbb{A}^1(\mathbb  C)$ is equivalent to the category of Stokes local systems $(\mathbb{V},\Theta)$.
\end{thm}

(We refer the reader to \cite[§8]{Boa21} for the definition of a morphism between two Stokes local systems, which is not completely obvious.)

\begin{rem}\label{rem:fissionLS}
	There is a different way of thinking about a Stokes local system, namely as a ``local system on a fission surface''. Let us briefly reflect on this notion here (as already alluded to in \cite[§4]{Boa09}, see also \cite[§3.1]{Boa14} and the picture on the title page of \cite{Boa21}). The principal idea is to split the above halo into several pieces given by the active exponents. 
    
    More precisely, we can construct a \emph{fission surface} $\underline{\widetilde{\Sigma}}(\Theta)$ as follows:
		
	Let $\Sigma=\mathbb{P}^1$ be the Riemann sphere, and fix an irregular class $\Theta$ at $\infty$.
	Let $\widehat{\Sigma}$ be the real blow-up of $\Sigma$ at $\infty$, with boundary circle $\partial$.
	
	Each active exponent in $\Theta$ (i.e.\ each connected component of the covering $\mathbb{I}\to\partial$ determined by $\Theta$) gives a covering circle $\langle q\rangle$ of $\partial$. We can consider a small strip $\langle q\rangle\times [0,\varepsilon]$ (for some $\varepsilon\ll 1$), and glue it to $\widehat{\Sigma}$ via the projection $\pi_q\colon \langle q \rangle\times \{0\} \simeq \langle q\rangle \to \partial$. In other words, we construct the pushout
	$$\begin{tikzcd}
		\langle q\rangle \arrow{r}{\pi_q} \arrow[equal]{d} & \partial\arrow[hook]{r} & \widehat{\Sigma}\arrow{dd}\\
		\langle q \rangle \times \{0\}\arrow[hook]{d}\\
		\langle q\rangle \times [0,\varepsilon]\arrow{rr} & & \widehat{\Sigma}\cup_{\pi_q} \big(\langle q\rangle \times [0,\varepsilon]\big)
	\end{tikzcd}$$
	Doing this for all active exponents and removing all the points $d\in\mathbb{A}\subset \partial$ yields the fission surface that we denote by $\underline{\widetilde{\Sigma}}=\underline{\widetilde{\Sigma}}(\Theta)$, i.e.\ concretely we have
	$$\underline{\widetilde{\Sigma}} \defeq \left(\widehat{\Sigma}\cup_{\pi_{q_1}} \big(\langle q_1\rangle\times [0,\varepsilon]\big) \cup_{\pi_{q_2}} \ldots \cup_{\pi_{q_n}} \big(\langle q_n\rangle\times [0,\varepsilon]\big)\right)\setminus \mathbb{A},$$
	where $\langle q_1\rangle,\ldots,\langle q_n\rangle$ are the connected components of $\mathbb{I}$.
	
	A \emph{fission local system} on this fission surface consists of a local system $\underline{\mathbb V}$ on $\widehat{\Sigma}$, a local system $\underline{\mathbb V}_i$ on any $\langle q_i\rangle \times [0,\varepsilon]$, and for any connected component $\delta$ of $\partial\setminus\mathbb{A}$ an isomorphism $\bigoplus_{i=1}^n\underline{\mathbb V}_i|_{\delta}\cong \underline{\mathbb V}|_\delta$.
	
	The $\underline{\mathbb V}_i$ define a local system on $\mathbb{I}\times [0,\varepsilon]$, and since a local system on $\mathbb{I}$ is the same as an $\mathbb{I}$-graded local system on $\partial\times[0,\varepsilon]$, we see that a fission local system is nothing but a local system on $\widetilde{\Sigma}(\Theta)$ graded by $\mathbb{I}$ on the halo $\mathbb{H}$. In particular, any Stokes local system can naturally be regarded as a fission local system. (Note that the notion of a fission local system associated to a Stokes-filtered local system is indeed exactly what has been defined in \cite[Definition~8.1]{Boa21}.)
	
	The advantage of this viewpoint is that the maps induced by the so-called deformation data later are naturally captured by the notion of a path in the fission surface, which we will explain in Remark~\ref{rem:defdataPath}.
\end{rem}

\subsection{Wild character varieties}
\label{subsection:wild_character_varieties}

In the same way as moduli spaces of local systems are character varieties, moduli spaces of Stokes local systems give rise to wild character varieties.

Keeping previous notations, let $(\mathbb{V},\Theta)$ be a Stokes local system. Let us denote by $\langle q_1\rangle,\dots,\langle q_k\rangle$ the active circles of $\Theta$, by $n_l$ the multiplicity of $\cir{q_l}$, and by $r_l$ its ramification order, for $l=1,\dots, k$. Choose a base point $b \in\partial$, and let us also fix a framing of $\mathbb{V}$ at $b$, that is, an isomorphism of vector spaces
\[
\mathbb{V}_b \cong \mathbb{F}
\]
respecting the $\Ical_b$-grading on both sides, where
\[
\mathbb{F}=\C^{\Theta}\defeq\bigoplus_{j\in \Ical_{b}} \C^{\Theta(j)}.
\]

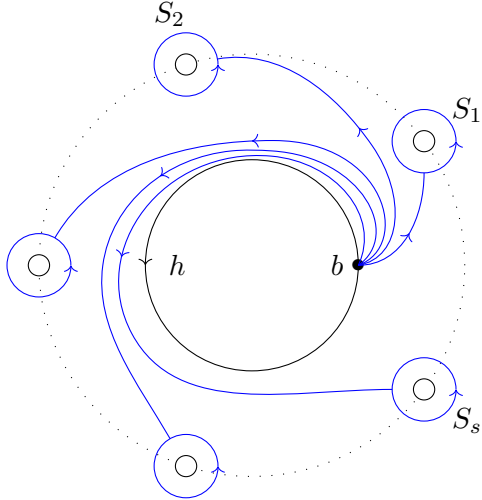
\begin{figure}[H]
\centering
	\begin{tikzpicture}[scale=1.4]
		
		\draw[loosely dotted] (0,0) circle (2);
		\foreach \x in {72,144,...,360} 
		{\draw[fill=white] ({\x+36}:2) circle (0.1);
			\draw[blue, ->] ({\x+36}:2)++(0.3,0) arc (0:360:0.3);
		} 
		\draw (1,0) node {$\bullet$};
		\draw (0.85,0) node {$b$};
		\begin{scope}[decoration={
				markings,
				mark=at position 0.5 with {\arrow{>}}}
			]  
			\draw [postaction={decorate}] (0,0) circle (1);
		\end{scope}   
		\draw (180:0.75) node {$h$};
		\draw[blue] (36:2.5) node {$S_1$};
		\draw[blue] (108:2.5) node {$S_2$};
		\draw[blue] (-36:2.5) node {$S_s$};
		\begin{scope}[blue,decoration={
				markings,
				mark=at position 0.5 with {\arrow{>}}}
			]  
			\draw[postaction={decorate}] (1,0) to [out=0, in=-90]($(0,-0.3)+(36:2)$);
			\draw[postaction={decorate}] (1,0) to [out angle=10, in angle=10,curve through={(45:1.65)}]($(10:0.3)+(3*36:2)$);
			\draw[postaction={decorate}] (1,0) to [out angle=25, in angle=60,curve through={(80:1.45)(100:1.5)}]($(60:0.3)+(5*36:2)$);
			\draw[postaction={decorate}] (1,0) to [out angle=45, in angle=120,curve through={(80:1.3)(100:1.3)(180:1.5)}]($(120:0.3)+(7*36:2)$);
			\draw[postaction={decorate}] (1,0) to [out angle=65, in angle=180,curve through={(80:1.15)(210:1.3)}]($(180:0.3)+(9*36:2)$);
		\end{scope} 
	\end{tikzpicture}
\caption{The paths for expressing the monodromies of a Stokes local system.}
\label{fig:SLS}
\end{figure}

Let $\Pi\defeq\pi_1(\widetilde{\Sigma},b)$, the fundamental group of $\widetilde{\Sigma}$ with base point $b$. Let  $\Hom(\Pi,G)$ be the set of representations of this group on $\mathbb{F}$, i.e.\ $G\defeq\mathrm{GL}(\mathbb F)\cong\mathrm{GL}_{n}(\mathbb{C})$, with $n=\mathop{\mathrm{rk}}\mathbb{V}=\sum_{l=1}^k r_ln_l$. A framed local system on $\widetilde\Sigma$ determines via its monodromy an element in $\Hom(\Pi,G)$, see Figure~\ref{fig:SLS}. Actually, because of the conditions in the definition of a Stokes local system, the representation associated to a framed Stokes local system lives in a subset $\Hom_{\mathbb{S}}(\Pi,G)$ of \textit{Stokes representations}, which we now describe.

Let us view the boundary circle $\partial$ as a loop based at $b$. (The orientation is given by that of $\PP^1$, i.e.\ counterclockwise in the local chart $t=z^{-1}$.)
Let $V^0=\mathbb{V}|_{\partial}\to \partial$ be the $\Ical$-graded local system on $\partial$ with irregular class $\Theta$ associated to $\mathbb{V}$. In the direct sum decomposition $V^0_{b}=\bigoplus_{l=1}^k (V^0_{\cir{q_l}})_{b}$, with  ${(V^0_{\langle q_l\rangle})}_{b}=\bigoplus_{j\in \langle q_l\rangle_{b}}V^0_b(j)$, its monodromy (taking $b$ as the base point) is of the form

\[\label{eq:formalMonblocks} \rho=\rho(\partial)=
\begin{pmatrix}
\rho_{\langle q_1\rangle} & \hdots & 0\\
\vdots & \ddots &  \vdots\\
0 & \hdots & \rho_{\langle q_k\rangle}\\
\end{pmatrix} \in\mathrm{GL}(V^0_{b})\cong \mathrm{GL}(\mathbb{F}),
\]
with each $\rho_{\langle q_l\rangle}$ of the form \eqref{form_formal_monodromy}.
The set of such matrices (seen as elements in $\mathrm{GL}(\mathbb{F}$)) is a twist of the group
\[ 
H\defeq\GrAut(\mathbb{F})\cong\prod_{j\in \Ical_{b}} \mathrm{GL}_{\Theta(j)}(\C),
\]
that we denote by $H(\partial)$ (see \cite[§4.2]{BY15} for details).

If $d\in \mathbb{A}$ is a singular direction, let $\lambda_d\subset \partial$ be an arc from $b$ to $d$ (in a positive sense). By parallel translation along $\lambda_d$, we may identify the Stokes group $\Sto_d\subset \mathrm{GL}(\mathbb{V}_d)$ with a subgroup of $G=\mathrm{GL}(\mathbb{F})\cong\mathrm{GL}_n(\C)$: Let us define
\[ \widehat{\gamma}_d\defeq\lambda_d^{-1} \circ \gamma_d \circ \lambda_d \in \Pi
\]
to be a simple loop around the tangential puncture $e(d)$ based at $b$, where $\gamma_d$ denotes, as before, the simple loop around $e(d)$ based on $d$. From the definition of a Stokes local system we immediately get that the monodromy $\rho(\widehat{\gamma}_d)$ of $\mathbb{V}$ around $\widehat{\gamma}_d$ belongs to $\rho(\lambda_d^{-1}).\Sto_d.\rho(\lambda_d)$, which via the framing $\mathbb F\cong \mathbb{V}_b$ corresponds to a unipotent subgroup of $G$ that we also call $\Sto_d$ by a slight abuse of notation. In summary:

\begin{lemma}\label{lemma:Stokesrep}
Let $\rho\in \Hom(\Pi,G)$ the representation associated to a Stokes local system $(\mathbb{V},\Theta)$. Then $\rho$ satisfies the following conditions:
\begin{enumerate}
\item $\rho(\partial)\in H(\partial)$,
\item $\rho(\widehat{\gamma}_d)\in \Sto_d$ for any singular direction $d\in \mathbb{A}$.
\end{enumerate}
\end{lemma}

\begin{defi}
A representation $\rho\in \Hom(\Pi,G)$ is a \textit{Stokes representation} if it satisfies the conditions (1) and (2) from Lemma~\ref{lemma:Stokesrep}. The subspace of $\Hom(\Pi,G)$ consisting of Stokes representations will be denoted by $\Hom_{\mathbb{S}}(\Pi,G)$.
\end{defi}

We will call the space $\cR_\Theta\defeq\mathrm{Hom}_{\mathbb S}(\Pi,G)$ of Stokes representations with irregular class $\Theta$ the \emph{representation variety} of $\Theta$. It has the following explicit description:
\[
\cR_\Theta=\{ (h, S_1,\dots, S_s)\in H(\partial)\times \Sto_{d_1}\times \dots \times \Sto_{d_s}\; |\; hS_s\dots S_1=\id\},
\]
where the set of singular directions $\mathbb{A}=\{d_1,\dots,d_s\}$ is ordered in a positive sense. (Here, the condition $hS_s\dots S_1=1$ comes from the fact that the composition of paths $\partial\circ \widehat{\gamma}_{d_s}\circ\ldots\circ\widehat{\gamma}_{d_1}$ is contractible in $\widetilde{\Sigma}$, hence the monodromy along this path must be the identity.)

When represented by explicit matrices, the elements $S_1,\ldots,S_s$ are usually called \emph{Stokes factors} (or Stokes matrices) and the element $h$ is the \emph{formal monodromy}.

There is a natural action of the group $H$ on $\cR_\Theta$, which amounts to changing the framing. 

The representation variety actually has a richer structure:

\begin{thm}[\cite{Boa14,BY15}]
The representation variety $\cR_\Theta$ is a smooth affine complex algebraic variety, and is a twisted quasi-Hamiltonian $H$-space, with moment map given by 
\begin{align*}
\mu\colon \cR_\Theta &\to H(\bar{\partial}),\\
\rho &\mapsto \{\rho(\bar{\partial})\},
\end{align*}
with $\bar{\partial}$ referring to the loop around $\partial$ in the negative sense.
\end{thm}

Finally, the (symplectic) wild character variety $\mathcal{M}_B(\Theta)=\mathcal M_B(\Sigma,\Theta)$ is defined as the affine GIT quotient of $\cR_\Theta$ by $H$. Concretely, the points of the wild character variety are the closed orbits, and they correspond to irreducible Stokes representations \cite{BY23}).

\begin{thm}[\cite{Boa14,BY15}]
The wild character variety $\mathcal M_B(\Sigma,\Theta)$ has an algebraic Poisson structure. Its symplectic leaves, the symplectic wild character varieties are the multiplicative symplectic quotients
\[
\mathcal M_B(\Sigma,\Theta,\mathcal C)=\cR_\Theta\sslash_{\mathcal C} H=\mu^{-1}(\mathcal C)/H,
\]
where $\mathcal C\subset H(\bar{\partial})$ is a twisted conjugacy class.
\end{thm}

If $(E,\nabla)$ is an algebraic connection on $\Sigma^\circ\defeq\Sigma\setminus\{\infty\}\subset \Sigma$, it determines an irregular class $\Theta\colon\pi_0(\Ical)\to \mathbb N$ and a formal local system $V^0\to \partial$ with irregular class $\Theta$. The isomorphism class of the inverse of the monodromy of $V^0\to \partial$ is a twisted conjugacy class $\mathcal C\subset H(\bar{\partial})$. The connection $(E,\nabla)$ therefore determines a symplectic wild character variety $\mathcal M_B(E,\nabla)\defeq\mathcal M_B(\Sigma,\Theta,\mathcal C)$.\\

It will also be useful to consider a partial reduction of the representation variety $\cR_\Theta$, obtained in the following way (see the last paragraph of \cite{BY20}): Instead of considering a framing of the full fibre $\mathbb{V}_b$, let us just choose a point $i_l\in \cir{q_l}_b$ and a framing of $V^0_{b}(i_l)$ for each active circle $\cir{q_l}$, which we will refer to as a \textit{minimal framing}. Forgetting the framings of all other graded pieces of the fibre $\mathbb V_b$ amounts to acting on $\cR_\Theta$ with a subgroup $H^{\perp} = \mathrm{GrAut}(\bigoplus_{j\in\mathcal
{I}_b\setminus\{i_1,\ldots,i_k\}} V^0_b(j))$ of $H$. Note that $H^{\perp} \cong \prod_{l=1}^k \mathrm{GL}_{n_l}(\mathbb C)^{r_l-1}$. Concretely, this action can be used to bring the pieces $\rho_{\cir{q_l}}$ of the formal monodromy $\rho(\partial)$ in \eqref{eq:formalMonblocks} to the form
\begin{equation}
	\rho_{\cir{q_l}}=\begin{pmatrix}
		0 & \hdots & \hdots &  0 & \rho_{r_l-1,0} \\
		1 & \ddots & & \vdots & 0\\
		0 & 1 & \ddots & \vdots & \vdots \\
		\vdots & \ddots & \ddots & 0 & \vdots\\
		0 & \hdots & 0 & 1 & 0
	\end{pmatrix}\in \Aut(V ^0_{\cir{q_l},b}),
\label{eq:monodromy_minimal_framing}
\end{equation}
where the first index corresponds to the point $i_l$ (cf.\ equation \eqref{form_formal_monodromy}).

\begin{defi}\label{def:reducedRepVar}
The reduced representation variety of $\Theta$ is the quotient $\mathcal{E}_\Theta\defeq\mathcal R_\Theta/H^{\perp}$.
\end{defi}

Furthermore, changing the minimal framing at $i_l$ amounts to acting with the residual group $\breve{H}\defeq \mathrm{GrAut}(\bigoplus_{j\in\{i_1,\ldots,i_k\}}V^0_{b}(j))\cong \prod_{l=1}^k \mathrm{GL}_{n_l}(\mathbb C)$. In turn, $\mathcal{E}_{\Theta}$ is a quasi-Hamiltonian $\breve H$-space, with moment map given by the monodromies around the $\cir{q_l}$ in the negative sense, and the symplectic wild character variety $\mathcal M_B(\Sigma,\Theta,\mathcal C)$ is its  (usual) quasi-Hamiltonian reduction $\mathcal M_B(\Sigma,\Theta,\mathcal C)=\mathcal{E}_{\Theta}\sslash_{\breve{\mathcal C}}\breve{H}$, where $\breve{\mathcal{C}}\subset \breve H$ is the (usual) conjugacy class in $\breve{H}$ induced by $\mathcal C$.

We investigate the concept of minimal framings more systematically in §\ref{sec:isomWCV}, where we reflect in particular on the relation between minimal framings on both sides of the Fourier transform.

\section{Deformation data}\label{sec:defData}

We now discuss a different description of Stokes data, used in \cite{Moc21} and formulated in terms of `deformation data' in loc.~cit. We will slightly  reformulate this description in order to relate it to Stokes local systems.

As before, we always consider the Riemann sphere $\Sigma=\mathbb{P}^1(\mathbb C)$ and an irregular class $\Theta$ at $\infty$.

Throughout this section and the rest of this article, we will make the following assumption:
\begin{assumption}
\label{assumption}
There is a fixed $k\in\mathbb{Q}_{>1}$ such that any active circle $\langle f\rangle$ of the irregular class $\Theta$ is of the form $\langle f\rangle = \langle a z^k\rangle$ for some $a\in\C\setminus\{0\}$, and for any two active circles $\langle a z^k\rangle, \langle b z^k\rangle\subset \cI$, one has $|a|=|b|$.
\end{assumption}
(We remark that Assumption \ref{assumption} is not necessary in §\ref{sec:distIntervals} and §\ref{sec:Legendre}--\ref{sec:LegendreDistIntervals}, unless explicitly stated.)

\subsection{Distinguished intervals}\label{sec:distIntervals}

One of the main ingredients in the deformation data description is the observation that on each active circle one has distinguished intervals along which the associated exponential is either growing or decreasing. We define them in this paragraph. Recall the notation from §\ref{sec:explocsys}.

Any circle $\langle f\rangle\subset\cI$ with $\cir{f}\neq \langle 0\rangle$ has a finite set $S_0(f)$ of distinguished points defined by
$$S_0(f) \defeq \{(\theta,[\tilde{f}])\mid \lim_{r\to\infty} \real \tilde{f}(r e^{i\theta})=0\}.$$
It is not difficult to see that $\# S_0(f)= 2m$ if $f=ax^{m/n}$ for $m,n$ coprime. The points of $S_0(f)$ correspond to the directions $d\in\partial$ such that $|e^{f(z)}|$ (for $z$ near $\infty$) passes from increasing to decreasing (or vice versa) when crossing the direction $d$. In a Stokes diagram as in Figure~\ref{fig:examples_of_Stokes_diagrams}, these are exactly the points where the circle $\cir{f}$ intersects the dashed one; indeed, this dashed circle corresponds to $f=0$, for which $e^{f(z)}$ is constant (notice, however, that in this picture each of the three intersection points of $\cir{f}$ with the dashed circle defines two elements of $S_0(x^{3/2})\subset \langle x^{3/2} \rangle$.) We also observe that the set of points $\pi_f(S_0(f))\subset \partial$ is precisely the set of Stokes directions of the pair $(\langle 0\rangle ,\langle f\rangle)$.

The distinguished intervals of a circle $\langle f\rangle$ are then defined as the intervals between the points of $S_0(f)$:

\begin{defi}
    Let $\langle f\rangle\subset \cI$ be a circle, $\cir{f}\neq\cir{0}$. We define $A(f)$ to be the set of connected components of $\langle f\rangle\setminus S_0(f)$. We denote by $A_+(f)\subset A(f)$ (resp.\ $A_-(f)\subset A(f)$) the subset of those intervals $I$ such that for some (and hence any) $(\theta,[\tilde{f}])\in I$, we have $$\real \tilde{f}(r e^{i\theta})>0 \qquad \text{(resp.\ $<0$)}$$
    for $r \gg 0$.
\end{defi}

Clearly, the intervals contained in $A(f)$ all have length $\frac{\pi\cdot n}{m}$ if $f=ax^{m/n}$.

In the following, given an interval $I=(a,b)$ in a circle (i.e.\ $I$ is to be understood as the image of this real interval in $S^1\simeq \mathbb{R}/2\pi\mathbb{Z}$) such that the interior of a positively oriented path from $a$ to $b$ lies entirely in $I$, we will call $a$ the \emph{startpoint} of $I$, $b$ the \emph{endpoint} of $I$ and $\frac{a+b}{2}$ the \emph{midpoint} of $I$.

An observation that will be essential in relating deformation data to Stokes local system is that pairwise intersections of distinguished intervals are in one-to-one correspondence with Stokes arrows: 

\begin{lemma}\label{lemma:intersSingDir}
	Let $k\in\mathbb{Q}_{>1}$, let $\cir{f}=\cir{ax^k}$ and $\cir{g}=\cir{bx^k}$ be exponents with $|a|=|b|$, and let $I\in A_+(f)$, $J\in A_-(g)$ be two sectors such that their projections to $\partial$ intersect nontrivially, i.e.\ $\pi_f(I)\cap \pi_g(J)\neq \emptyset$. Then the midpoint $d\in \partial$ of the intersection $\pi_f(I)\cap \pi_g(J)$ is a singular direction for the pair $(\cir{f},\cir{g})$, and there is a Stokes arrow $d_I \to_d d_J$ over $d$, where $d_I,d_J$ are the unique points of $I$ and $J$ lying over $d$.
	
	Conversely, for every Stokes arrow $d_f\to_d d_g$ there is a pair of sectors $I\in A_+(f)$, $J\in A_-(g)$ with $d_f\in I$, $d_g\in J$ such that $d$ is the midpoint of $\pi_f(I)\cap \pi_g(J)$.
\end{lemma}
\begin{proof}
    This is clear from the definition of a Stokes arrow: having a Stokes arrow means that the real parts of  the corresponding branches of the exponents in question have a maximal difference (their difference is a point of maximal decay). This certainly happens at a point where the two real parts differ in sign, and since the absolute values of the coefficients are equal, the maximal difference must appear in the middle of the intersection for symmetry reasons.
    
    (Note that the length of $I$ and $J$ is smaller than $\pi$, which follows from our assumption $k>1$, hence $\pi_f(I)\cap \pi_g(J)$ is connected so that it is actually reasonable to speak about ``the midpoint''.)
\end{proof}

\subsection{Preferred splittings on sectors}

We now briefly recall a few more facts that come into play in the definition of deformation data. The latter are certain morphisms constructed in \cite[§3.4.2]{Moc21} starting from the data of a Stokes-filtered local system, which is yet another way of describing Stokes data (the local system of analytic solutions of an algebraic connection on the affine line comes naturally equipped with the structure of a Stokes-filtered local system). While we refer the reader, for example, to \cite{Mal91,Sab13,Boa21} for the full definition and further details about Stokes-filtered local systems, we will only use here the main fact allowing to pass from Stokes-filtered local systems to deformation data or Stokes local systems, which is the following property.

Let $(V,\Theta,F)$ be a Stokes-filtered local system on $(\PP^1,\infty)$ with irregular class $\Theta$, and let $V^0\to \partial$ be the corresponding formal local system (obtained by taking the associated graded on $\partial$). 
Then any direction $\theta\in\partial$ that is not a singular direction determines a preferred isomorphism
\begin{equation*}
V|_S\cong V^0|_S
\end{equation*}
on a sufficiently small neighbourhood $S\subset \partial$ of $\theta$, which consequently gives a splitting of the Stokes filtration
\begin{equation}\label{eq:cansplitting}
	V|_S\cong \bigoplus_{I\in\pi_0(\pi^{-1}(S))} V^0|_I.
\end{equation}

\begin{rem}
	The right-hand side of \eqref{eq:cansplitting} needs some explanation: In order to restrict $V^0$ to $I\subset \mathcal{I}$, we consider the $\mathcal{I}$-graded local system $V^0$ as a local system on the covering space $\mathcal{I}$ (as explained in §\ref{sec:formalLS}). The restriction $V^0|_I$ then has to be understood as a local system on $S$ again via the homeomorphism $\pi|_I\colon I\overset{\sim}{\to} S$ for any $I\in\pi_0(\pi^{-1}(S))$.	
	Hence, without any abuse of notation, this formula would be written as
	$$V|_S\cong \bigoplus_{I\in\pi_0(\pi^{-1}(S))} (\pi|_I)_*(\widetilde{V}^0|_I).$$
	We will, however, use the simplified notation \eqref{eq:cansplitting}, with implicit identifications as just explained, in the rest of this section, in order to keep notation reasonable.
\end{rem}

The meaning of ``sufficiently small neighbourhood $S$'' for \eqref{eq:cansplitting} to hold can be made more precise: Let us introduce some terminology. A \emph{singular sector} is an open interval between two neighbouring singular directions. In particular, any $\theta\in \partial$ which is not a singular direction lies in a unique singular sector. A \emph{supersector} is an open interval obtained from a singular sector by enlarging it by $\frac{\pi}{2k}$ on both sides. The preferred isomorphism \eqref{eq:cansplitting} locally around $\theta$ can be extended to the singular sector containing $\theta$ and even the associated supersector.

In general, a small interval $S$ can be part of different supersectors, yielding different isomorphisms of the form \eqref{eq:cansplitting} on $S$, all coming from preferred splittings for some $\theta$. However, if an interval $S\subset \partial$ is large enough, one has a unique choice for such a splitting valid on the whole of $S$, as the following lemma shows.
\begin{lemma}
    Let $S=(\theta_0,\theta_0+\frac{\pi}{k}+\varepsilon)\subset \partial$ for some $0<\varepsilon \ll 1$. Then there is a unique supersector containing $S$ and the midpoint of $S$ lies in the associated singular sector.
\end{lemma}
\begin{proof}
    The overlap of two consecutive supersectors is an open interval of angle $\frac{\pi}{k}$. This shows that $S$ cannot be contained in two supersectors. On the other hand, it is not difficult to see that $S$ is certainly contained in one supersector, and that the midpoint of $S$ will lie in the central part of the supersector, which is the associated singular sector.
\end{proof}

The lemma implies in particular the following: Let $S=\pi_f(I)$ for $\langle f\rangle\in\cI$ and $I\in A(f)$, and let $S_+$ be the interval $S$ closed on the right (at its endpoint). (Alternatively, $S_+$ can be thought of as the interval $S$ extended slightly across its endpoint to an open interval). Then one has a unique splitting of the form \eqref{eq:cansplitting} on $S_+$, and it is the splitting induced by that on the singular sector containing the midpoint of $S$ or, if the midpoint is a singular direction, the singular sector starting there. This splitting is called the \emph{canonical splitting} on $S_+$. A similar statement obviously holds for $S_-$.

\subsection{Deformation data and comparison of canonical splittings}\label{sec:defDataAndSplittings}
    We are now in the position to reinterpret the notion of deformation data from \cite{Moc21} in the language of canonical splittings introduced above. In what follows, we tacitly make the following identification: If $I$ is a contractible set, $p\in I$ a point and $V_1,V_2$ local systems on $I$, then it is equivalent to give one of the following: a morphism $V_1\to V_2$; a morphism $(V_1)_p\to (V_2)_p$; a morphism $V_1(I)\to V_2(I)$; a morphism $V_1|_J\to V_2|_J$ for a contractible subset $J\subset I$. Moreover, if $I\in A(f)$, we will write $I_\partial\defeq \pi_f(I)\subset \partial$. For each active circle $\langle f\rangle$, we denote by $V^0_{\langle f\rangle}\subset V^0$ the corresponding graded piece of a formal local system $V^0$.\\

Let $V$ be a Stokes-filtered local system on $(\PP^1,\infty)$ with irregular class $\Theta$, and let $V^0$ be the associated formal local system on $\partial$.
In \cite[Ch.~3]{Moc21}, a concept equivalent to that of a Stokes-filtered local system on $S^1$ is introduced: the notion of a Stokes shell. We note that, given our assumptions, the local systems $\mathcal{K}_\lambda$ (for an exponent $\lambda$) from loc.~cit.\ correspond to the formal pieces of the formal local system $V^0$ in our terminology. 

In our case, such a Stokes shell consists of the data of the corresponding formal local system together with so-called \emph{deformation data}. The latter are linear maps of the following form: Let $\cir{f},\cir{g}\subset \cI$ be circles and $I\in A_+(f)$, $J\in A_-(g)$ such that $I_\partial\cap J_\partial\neq\emptyset$, and denote by $m_I,m_J$ their midpoints, then a deformation datum is a linear map
$$\cR^I_J\colon (V^0_{\langle f\rangle})_{m_I} \to (V^0_{\langle g\rangle})_{m_J}.$$

By Lemma~\ref{lemma:intersSingDir}, this means that there will be exactly one deformation datum for every Stokes arrow. Here, we will relate the deformation data morphism $\cR^I_J$ to the comparison of the canonical splittings on $I$ and $J$.

\begin{prop}\label{prop:defDataTransitionIJ}
Let $\langle f\rangle$ and $\langle g\rangle$ be active exponents with respect to the irregular class $\Theta$. Let $I\in A_+(f)$ and $J\in A_-(g)$ be two intervals such that $I_\partial\cap J_\partial\neq \emptyset$.
Consider the splittings 
$$V|_{I_\partial} \underset{\sim}{\overset{\varphi}{\To}} \bigoplus_{D\in\pi_0(\pi^{-1}(I_\partial))} V^0|_D \qquad \text{and} \qquad V|_{J_\partial} \underset{\sim}{\overset{\psi}{\To}} \bigoplus_{D\in\pi_0(\pi^{-1}(J_\partial))} V^0|_D$$
coming from the canonical splittings on $(I_\partial)_-$ and $(J_\partial)_+$, respectively. Then the associated morphism
\begin{align}
V^0_{\langle f\rangle}|_{I\cap \pi^{-1}(J_\partial)} \hookrightarrow &\bigoplus_{U\in\pi_0(\pi^{-1}(I_\partial\cap J_\partial))} V^0|_U\notag\\
&\underset{\sim}{\overset{\psi\circ \varphi^{-1}}{\To}} \bigoplus_{U\in\pi_0(\pi^{-1}(I_\partial\cap J_\partial))} V^0|_U \twoheadrightarrow V^0_{\langle g\rangle} |_{J\cap \pi^{-1}(I_\partial)}.\label{eq:defdatanatural}
\end{align}
is equivalent to Mochizuki's deformation datum $\cR^I_J$ associated to this situation.
\end{prop}

\begin{proof}
	If $I_\partial = J_\partial$, the statement is clear from the construction in \cite[§3.4.2]{Moc21}.
	
	In the case where $I_\partial\neq J_\partial$, let us first remark that the morphism \eqref{eq:defdatanatural} does not change if we replace the splitting coming from $(I_\partial)_-$ by the one from $(I_\partial)_+$ and the splitting coming from $(J_\partial)_+$ by that from $(J_\partial)_-$.
    
    In this case, the author of loc.~cit.\ does not directly define deformation data in terms of $\psi\circ \varphi^{-1}$. In fact, he gives two different constructions, depending on if $J_\partial$ ``comes after'' $I_\partial$ (we write $J_\partial > I_\partial$, by which we mean that $J_\partial$ contains the endpoint of $I_\partial$, thinking in terms of the given orientation of $\partial$) or not. We will therefore need to distinguish two cases, and we will start with the one where $J_\partial>I_\partial$. In this case, the construction in \cite[§3.4.2]{Moc21} is the following (this is a slight reformulation of p.\ 22 in loc.~cit.\ in our situation and language):
	
	Let $A\defeq \bigsqcup_{\Theta(\cir{q})\neq 0} A(q)$ be the set of all the distinguished sectors on all active circles. One considers the situation at the point $\theta_0$, the endpoint of $I_\partial$ (which is contained in the interior of $J_\partial$). For any sector $U\in A$ such that $\theta_0\in (U_\partial)_-$ there is a splitting
	\begin{equation}\label{eq:decompTheta0}
	V|_{U_\partial} \underset{\sim}{\overset{\varphi_U}{\To}} \bigoplus_{S\in\pi_0(\pi^{-1}(U_\partial))} V^0|_S
	\end{equation}
	coming from that on $(U_\partial)_-$.
	Then, in a small neighbourhood $W$ of $\theta_0$ in $\partial$, we have an isomorphism (see \cite[Lemmas~2.3.1 and 2.3.2]{Moc21} for such a decomposition)
	$$\nu\colon V|_W \overset{\sim}{\longrightarrow} \bigoplus_{U\in A, \theta_0\in (U_\partial)_-} V^0|_U|_W$$
	defined by combining all those from \eqref{eq:decompTheta0} in the natural way. That is, for any $U\in A$ with $\theta_0\in(U_\partial)_-$, we consider the morphism $V^0|_U \to V|_{U_\partial}$ induced by the corresponding splitting \eqref{eq:decompTheta0}. It induces a morphism $V^0|_U|_W\to V|_W$ (where we enlarge $U$ slightly if necessary in order to ensure $W\subset U$), and all these morphisms together give a morphism
    $$\bigoplus_{U\in A, \theta_0\in (U_\partial)_-} V^0|_U|_W\to V|_W,$$
    which can be shown to be an isomorphism. The above isomorphism $\nu$ is then defined to be its inverse. In order to make this clear, let us describe the construction of the morphism $\nu$ differently:
	Take the natural total order $\leq$ on the sectors $U\in A$ with $\theta_0\in (U_\partial)_-$ given by comparing the real parts of the corresponding exponents at $\theta_0$. Start with a maximal $U_1$ and the associated splitting
	$$V|_{(U_1)_\partial} \underset{\sim}{\overset{\varphi_{U_1}}{\To}} \bigoplus_{S\in\pi_0(\pi^{-1}((U_1)_\partial))} V^0|_S.$$
	Now take a maximal element $U_2$ of the remaining sectors. Then $\varphi_{U_2}\circ \varphi_{U_1}^{-1}$ induces on $(U_1)_\partial\cap (U_2)_\partial$ an isomorphism (since it has to be compatible with the filtration induced by the ordering)
	$$ \bigoplus_{D\in\pi_0(\pi^{-1}((U_1)_\partial))\setminus \{U_1\}} V^0|_D \cong \bigoplus_{D\in\pi_0(\pi^{-1}((U_2)_\partial))\setminus \{U_1\}} V^0|_D.$$
	Composing this isomorphism (taking it to be the identity on the summand corresponding to $S=U_1$) with $\varphi_{U_1}$, we find
	$$V|_{(U_1)_\partial\cap (U_2)_\partial} \overset{\sim}{\To} V^0|_{U_1}|_{(U_1)_\partial\cap (U_2)_\partial} \oplus \bigoplus_{D\in\pi_0(\pi^{-1}((U_2)_\partial))} V^0|_D|_{(U_1)_\partial\cap (U_2)_\partial}.$$
	Continuing this procedure yields the desired isomorphism $\nu$. In what follows, we will shrink $W$ if necessary.
	
	Now, the deformation datum is defined in \cite[§3.4.2]{Moc21} (see also \cite[§2.3.4]{Moc21}) as
	$$V^0_{\langle f\rangle} |_I|_W \hookrightarrow \bigoplus_{D\in\pi_0(\pi^{-1}(I_\partial))} V^0|_D|_W \underset{\sim}{\overset{\nu\circ \varphi^{-1}}{\To}} \bigoplus_{U\in A, \theta_0\in (U_\partial)_-} V^0|_U|_W\twoheadrightarrow V^0_{\langle g\rangle} |_J|_W.$$
    (Here, we consider $I$ slightly enlarged on the right, so that $I_\partial\cap W$ actually includes $\theta_0$.)
    
	It is not too difficult to see that, by the construction of $\nu$ given above, this corresponds to the desired morphism \eqref{eq:defdatanatural}, since for any $U\in A$, the other components of $\nu\circ\varphi^{-1}(V^0_{\langle f\rangle}|_I|_U)$ give no other contribution to the component $V^0_{\langle g\rangle} |_J|_W$ due to the fact that the overlap of any such $U$ with $J$ contains a Stokes direction of the pair $(f_U,g)$. This completes the proof for the case $J_\partial> I_\partial$.
	
	Let us now indicate the proof for the other case, when $I_\partial > J_\partial$.
    The construction of Mochizuki goes differently here, namely through the interval $I-\frac{\pi}{k}$ as follows: Denote by $\widetilde\varphi$ the splitting on this shifted sector induced by the one on $((I-\frac{\pi}{k})_\partial)_+$. Let $\theta_0$ be the startpoint of $I_\partial$. Then we can construct a splitting $\nu$ as in the proof above. By the construction in \cite[§3.4.2]{Moc21}, the deformation datum is $(-1)$ times the morphism
    $$V^0_{\langle f\rangle}|_{I}|_W \cong V^0_{\langle f\rangle}|_{I-\frac{\pi}{k}}|_W \hookrightarrow \bigoplus_{D\in\pi_0(\pi^{-1}((I-\frac{\pi}{k})_\partial))} V^0|_D|_W \overset{\nu\circ \widetilde{\varphi}^{-1}}{\To} \bigoplus_{U\in A, \theta_0\in (U_\partial)_-} V^0|_U|_W\twoheadrightarrow V^0_{\langle g\rangle}|_J|_W.$$
    (Here again, intervals are considered slightly enlarged on both sides if necessary.)

    It now remains to see that the image of $V^0_{\langle f\rangle}|_{I-\frac{\pi}{k}}|_W$ under $\nu\circ \widetilde{\varphi}^{-1}$ contributes a summand of $-\lambda(V^0_{\langle f\rangle}|_{I}|_W)$ to the $J$-th graded piece, where $\lambda\colon V^0_{\langle f\rangle} |_{I}|_W \to V^0_{\langle g\rangle}|_{J}|_W$ corresponds to the map \eqref{eq:defdatanatural} that we want to get. Finally, to conclude, it suffices to see that during this procedure no other, ``unwanted'' contributions are made to the $J$-factor of the target.
\end{proof}

\begin{rem}
		Since the idea and proof of the previous proposition might be hard to understand, let us illustrate it in a concrete example:
		
		Consider a Stokes-filtered local system of rank $3$ with an irregular class consisting of $3$ exponents $f$, $g$ and $h$, such that the picture looks as follows:
		$$
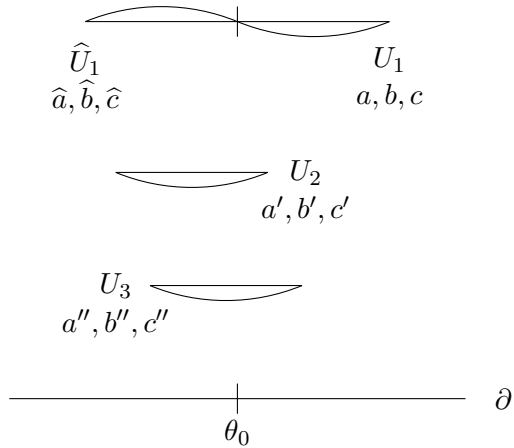
\begin{tikzpicture}
			\draw (0,0) -- (2,0);
			\draw (0,-0.5) node {$\widehat{U}_1$};
			\draw (0,-1) node {$\widehat{a},\widehat{b},\widehat{c}$};
			\draw (2,0) -- (4,0);
			\draw (2,-0.2) -- (2,0.2);
            \draw (2,-5.2) -- (2,-4.8);
			\draw (2,-5.45) node {$\theta_0$};
			\draw (2,0) arc (247.5:292.5:2.6131);
			\draw (2,0) arc (67.5:112.5:2.6131);
			\draw (4,-0.5) node {$U_1$}	;
			\draw (4,-1) node {$a,b,c$};		
			\draw (0.85,-3.5) -- (2.85,-3.5);
			\draw (0.85,-3.5) arc (247.5:292.5:2.6131);
			\draw (0.4,-3.5) node {$U_3$};
			\draw (0.4,-4) node {$a'',b'',c''$};
			
			\draw (0.4,-2) -- (2.4,-2);
			\draw (0.4,-2) arc (247.5:292.5:2.6131);
			\draw (2.9,-2) node {$U_2$};
			\draw (2.9,-2.5) node {$a',b',c'$};
			
			\draw (-1,-5) -- (5,-5);
			\draw (5.5,-5) node {$\partial$};
		\end{tikzpicture}$$
		
		For each of these intervals $U_i$, one has a canonical splitting for the local system $V$ on $((U_i)_\partial)_-$, and one also has a canonical splitting on $(\widehat{U}_\partial)_+$. Each splitting corresponds to choosing a basis of the local system on this interval downstairs. We denote these bases by $a,b,c$ etc., as shown in the picture. On the other hand, Mochizuki's splitting, denoted $\nu$ in the above proof, corresponds to considering the basis $a,b', c''$. The statement of the proposition in the case $J_\partial>I_\partial$ is then the following: If we express $\widehat{a}$ in the basis $a'',b'',c''$ and in the basis $a,b',c''$, then the coefficient of $c''$ is the same in both expressions.
		
		Let us also explain the case $I_\partial>J_\partial$: 		
		Let us assume that, if we express $a$ in terms of $\widehat{a},\widehat{b},\widehat{c}$, the coefficient of $\widehat{a}$ is $1$ (this can always be achieved simply by rescaling one of these bases).
		The statement of the second proposition then says the following in this case: Express $\widehat{a}$ in terms of the basis $a,b',c''$, and express $a$ in terms of the basis $a'',b'',c''$, then the coefficients of $c''$ in both expressions only differ by a sign.
		
		Both statements can be checked by an easy computation, and it is tedious but not difficult to generalise these arguments to more complicated cases, as we indicated in the above proof above.
	\end{rem}

Let us note that, if $I_\partial\neq J_\partial$, then the deformation datum in the above proposition does not change if we replace the splitting coming from $(I_\partial)_-$ by the one from $(I_\partial)_+$, and similarly for the splitting on $J_\partial$.
However, in the case $I_\partial=J_\partial$, the deformation datum really depends on the particular choice of splitting.

\subsection{Stokes local systems and deformation data}

Having understood the deformation data as in the previous subsection, we can also associate similar data to Stokes local systems. (For the equivalence between Stokes-filtered local systems and Stokes local systems, we refer to \cite{Boa21}.) First of all, in view of the above observations, let us give the following definition:

\begin{defi}
    Let $\mathbb V$ be a Stokes local system with irregular class $\Theta$, and let $V^0$ be the corresponding formal local system.
	Let $\langle f\rangle$ be an active circle and let $I\in A(f)$ be a distinguished interval.
	
	If $I\in A(f)$, we define the map
	$$\iota_I\colon (V^0_{\langle f\rangle})_{m_I} \hookrightarrow \mathbb V_{\pi_f(m_I)}$$
	induced by the inclusion of a direct summand and the splitting on the singular sector just to the left of $\pi_f(m_I)$ (or at $\pi_f(m_I)$ if it is not a singular direction). (Here, ``to the left'' means ``going backwards with respect to the given orientation on $\partial$''.)
	
	If $I\in A(f)$, we define the map
	$$\pi_I\colon \mathbb V_{\pi_f(m_I)} \twoheadrightarrow (V^0_{\langle f\rangle})_{m_I}$$
	induced by the splitting on the singular sector just to the right of $\pi_f(m_I)$ (or at $\pi_f(m_I)$ if it is not a singular direction) and the projection to a direct summand.
\end{defi}

In view of Proposition~\ref{prop:defDataTransitionIJ}, we define as follows the deformation datum associated to a Stokes local system:

\begin{defi}\label{def:defDataSLS}
	Let $d$ be a singular direction and let $d_f\to_d d_g$ be a Stokes arrow over it. Then there exist sectors $I$ and $J$ as in Lemma~\ref{lemma:intersSingDir}. Denote their midpoints by $p\defeq m_I$ and $q\defeq m_J$, respectively. Then the associated deformation datum is the map
	$$R^p_q\colon (V^0_{\langle f\rangle})_{m_I}\longrightarrow (V^0_{\langle g\rangle})_{m_J}$$
	given by the composition
	$$(V^0_{\langle f\rangle})_{m_I}\overset{\iota_I}{\hookrightarrow} \mathbb V_{\pi_f(m_I)}\xlongrightarrow{\sim} \mathbb V_{\pi_g(m_J)} \overset{\pi_J}{\twoheadrightarrow} (V^0_{\langle g\rangle})_{m_J},$$
	where the middle arrow is parallel transport in the local system in the interior (i.e.\ not in the halo) of $\widetilde{\Sigma}(\Theta)$.
\end{defi}

Proposition~\ref{prop:defDataTransitionIJ} implies that the deformation data in the sense of \cite{Moc21} associated to a Stokes-filtered local system $V$ coincides with the deformation data associated by the definition above to the Stokes local system $\mathbb V$ corresponding to $V$.

\begin{rem}\label{rem:defdataPath}
    We can relate the deformation data defined above to the parallel transport in the Stokes local system along distinguished paths on $\widetilde{\Sigma}(\Theta)$ as follows:
    Given a Stokes arrow, it determines two sectors $\pi_f(I)$ and $\pi_g(J)$ on $\partial$, where the corresponding real parts of the exponents are positive and negative, respectively. We can connect their midpoints by a path in $\widetilde\Sigma(\Theta)$. Parallel transport along this path in the local system gives an isomorphism of stalks
    $$\mathbb{V}_{\pi_f(m_I)} \overset{\sim}{\longrightarrow} \mathbb{V}_{\pi_g(m_J)}.$$
    Pre- and post-composing this with the inclusion $\iota_I$ and projection $\pi_J$ associated to the $\langle f\rangle$- and $\langle g\rangle$-graded pieces, respectively, yields the deformation datum from Definition~\ref{def:defDataSLS}.

    In what follows, we will illustrate this by a path as drawn in Figure~\ref{fig:defdatum}.
    
    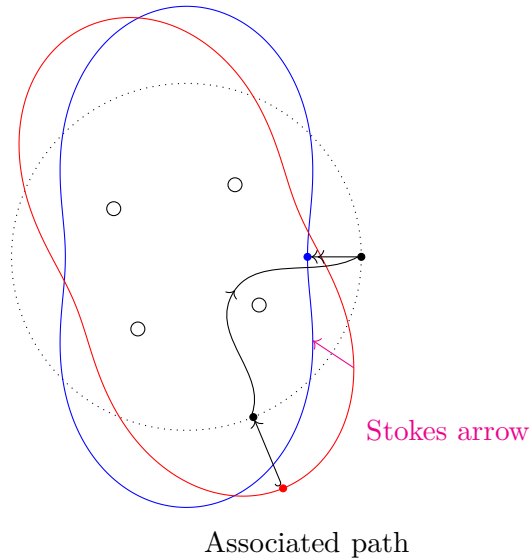
\begin{figure}
    	\centering
    \begin{tikzpicture}[scale=2.3]
    \draw[dotted] (0,0) circle (1);
      \draw[blue,domain=0:(360),scale=1,samples=1000] plot (\x:{exp(cos(2*\x+180)/exp(1))});
       \draw[red,domain=0:(360),scale=1,samples=1000] plot (\x:{exp(cos(2*\x+180-45)/exp(1))});  
    \foreach \n in {1,2,3,4} 
    \draw ({45-\n*90+180/16}: 0.5) circle (0.04cm);
        
       \draw[red, fill=red] (-67.5:{exp(1/exp(1))}) circle (0.02cm);
       \draw[fill=black] (-67.5:1) circle (0.02cm);
        \draw[fill=black] (0:1) circle (0.02cm);
        \draw[blue, fill=blue] (0:{exp(-1/exp(1))}) circle (0.02cm);
       \draw[{Hooks[round,right]}->] (-67.5:{exp(1/exp(1))-0.02}) to (-67.5:1.02); 
       \draw[decoration={markings,mark=at position 0.5 with {\arrow{>}}},postaction={decorate}] (-67.5:0.98) to [out angle=70, in angle=-150,curve through={(-45:0.35)}](0:0.98);
       \draw[->>] (0:0.98) to (0:{exp(-1/exp(1))+0.02}); 
       \draw[->, magenta] ({-45+180/16}: {exp(cos(2*(-45+180/16)+180-45)/exp(1))}) to ({-45+180/16}: {exp(cos(2*(-45+180/16)+180)/exp(1))});
       \draw[magenta] ({-40+180/16}: 1.65) node {Stokes arrow};
       \draw  (-67.5:1.8) node  {Associated path};
    \end{tikzpicture}
    \caption{A deformation datum associated to a Stokes arrow over a singular direction can be pictured as a path from a leaf to another. Note that at the boundary the surface is equipped with a covering by the exponents, so it is best to picture the blue and red circles above the black one (in the third dimension).}
    \label{fig:defdatum}
    \end{figure}

    All this can be formulated in an even cleaner way if we view a Stokes local system as a fission local system as explained in Remark~\ref{rem:fissionLS}: One can draw a path on the fission surface as follows: Recall that the fission surface $\underline{\widetilde{\Sigma}}(\Theta)$ looks like the real blow-up $\widehat\Sigma$ glued to a small strip for each exponent at the boundary $\partial$, and punctured at the singular directions. Start at the point $(m_I,\varepsilon) \in \langle f\rangle \times [0,\varepsilon]$. Then connect it by a straight line to $(m_I,0)\in \langle f\rangle \times [0,\varepsilon]$, which is identified with a point in $\partial$. We can draw a path in $\widehat{\Sigma}$ from this point to $(m_J,0)\in \langle g\rangle \times [0,\varepsilon]$ and connect it again by a straight line to $(m_J,\varepsilon)$. This is a path in the surface $\underline{\widetilde{\Sigma}}(\Theta)$, and if we interpret every transition of the circle $\partial$ (where the surface ``splits'') as the obvious inclusion or projection, the parallel transport in the fission local system along this path naturally represents the deformation data. A path that is obtained in this way from a Stokes arrow is called a \emph{Stokes path}. The algorithm from Theorem~\ref{thm:trafoRule} is explained in this spirit in the introduction (see Theorem~\ref{thm:algoPaths}).
\end{rem}

A Stokes local system can be reconstructed from its deformation data together with its formal local system  (i.e.\ the data referred to as a Stokes shell by in \cite{Moc21}).

\begin{prop}\label{prop:SLSdefdata}
	A Stokes local system $\mathbb V$ with irregular class $\Theta$ is determined by the following data:
	\begin{itemize}
		\item For each active exponent $\langle f\rangle$ the ($\langle f\rangle$-graded) local system $V^0_{\langle f\rangle}$ on $\partial$. It can be viewed as a local system on $\langle f\rangle$, and is hence given by the following data:
		
		\noindent for any $I\in A(f)$, the vector space $K^f_I\vcentcolon = (V^0_{\langle f\rangle})_{m_I} \cong \Gamma(I,V^0_{\langle f\rangle} )$, and for any two adjacent intervals $I,J\in A(f)$, $J$ following $I$ with respect to the given orientation, an isomorphism $K^f_I\overset{\sim}{\To} K^f_J$,
		\item for each Stokes arrow, the corresponding deformation datum $R^p_q$, as in Definition~\ref{def:defDataSLS}.
	\end{itemize}
\end{prop}
\begin{proof}
    This observation follows, for example, from the fact that the notion of Stokes shell is equivalent to that of a Stokes-filtered local system and hence to that of a Stokes local system (combining \cite[Proposition~3.4.3]{Moc21} and \cite[§9 and Theorem~11.3]{Boa21}).
\end{proof}

\section{The algorithm}\label{sec:alg}
In this section, we describe the algorithm for computing the Fourier transform of Stokes data in the class of cases we are interested in here. It is based on the results of \cite{Moc21}, but we (re)formulate it in the language of Stokes local systems, so that it yields in particular isomorphisms between wild character varieties.

\subsection{The Legendre transform as a homeomorphism of circles}\label{sec:Legendre}

A well-known fact about the Fourier transform on irregular connections on the affine line is that there exists an explicit way of determining from the irregular class $\Theta$ of a connection the irregular class $\widehat{\Theta}$ of its Fourier transform. This relation is called the stationary phase formula \cite{Mal91, Sab08}, and it states that the active circles of $\widehat\Theta$ are obtained from those of $\Theta$ by a Legendre transform.  

For exponents at infinity of slope $>1$ (i.e.\ with leading term $az^{k}$, $k\in\mathbb{Q}_{>1}$, which includes our case of interest), the Legendre transform is given as follows: If $f(z)$ is such an exponent, its Legendre transform $g(w)=\widehat{f}(w)$ is determined by the equations $$w=f'(z)\quad \text{and} \quad g(w)=f(z)-zw.$$

The Legendre transform provides a permutation of the set of all circles of slope $>1$ in $\mathcal{I}$, and this is the way it is often presented in the literature. The crucial idea underlying the description of \cite{Moc21} of the Fourier transform of Stokes data is that the Legendre transform can be seen in a stronger way, not only as a bijection between the sets of active circles of $\Theta$ and $\widehat{\Theta}$, but for any active circle $\cir{f}$ of $\Theta$, as a homeomorphism between $\cir{f}$ and the corresponding circle $\cir{\widehat f}$ of $\widehat{\Theta}$. 

This can be observed as follows (what follows is similar to the presentation in \cite[§3.1]{Dou21}, to which we refer for details): Let us denote the circle of directions at $\infty$ before Fourier transform by $\partial$ and the one after Fourier transform by $\partial'$. As explained in §\ref{sec:explocsys}, a point of $\langle f\rangle$ is a pair $(\theta,[\tilde{f}])$ with $\theta\in\partial$ and $\tilde{f}$ a branch of $f$ near $\theta$ (so $\tilde{f}(z)$ is a holomorphic function on a small sector $S(U)$ defined by a small neighbourhood $U\subset \partial$ of $\theta$; we shrink $U$ and $S(U)$ if necessary).

Now we consider the function $\phi(z)\defeq\tilde{f}'(z)=\frac{\mathrm{d}\tilde{f}}{\mathrm{d}z}(z)$. It defines a biholomorphic map between $S(U)$ and some small sector $S(U')$ for $U'\subset \partial'$. In particular, it sends $\theta\in \partial$ to a unique point $\theta'\in\partial'$. Concretely, we can set $\theta' = \lim_{r\to\infty} \arg \tilde{f}'(r\cdot e^{i\theta})$.

Next, define a function on $S(U')$ by $\tilde{g}(w)\defeq \tilde{f}(\phi^{-1}(w))-\phi^{-1}(w)w$. Continuing this function holomorphically around $\infty$ defines a circle $\langle g\rangle$ (a covering of $\partial'$) on which $(\theta',[\tilde{g}])$ is a point. It can be checked that this defines a homeomorphism of circles $\ell_f\colon \langle f\rangle \overset{\sim}{\to} \langle g\rangle$ preserving the orientation.

\subsection{Correspondence between distinguished intervals}\label{sec:LegendreDistIntervals}

We now show that the Legendre transform behaves well with respect to the distinguished intervals. This appears to be a main reason why T.~Mochizuki prefers them to the actual Stokes or singular directions and builds his theory around them (see \cite[§1.4.5]{Moc21}).

\begin{lemma}\label{lemma:LegendreS0}
	If $\langle g\rangle$ is the Legendre transform of an exponent $\langle f\rangle$ of slope $>1$, then the homeomorphism $\ell_f\colon \langle f \rangle \overset{\sim}{\to} \langle g\rangle$ induces a bijection between $S_0(f)$ and $S_0(g)$.
\end{lemma}

\begin{proof}
	First of all, note that the statement makes sense since $\# S_0(f) = \# S_0(g)$: If the highest power appearing in $f$ is $z^{m/n}$, then the highest power appearing in $g$ is $w^{m/(m-n)}$. Since the Legendre transform $\ell_f$ is a homeomorphism between two circles, it suffices to show that any point of $S_0(f)$ is mapped to a point of $S_0(g)$.
	
	Let $(\theta,[\tilde{f}])\in S_0(f)$, i.e.\ $\lim_{r\to\infty}\real \tilde{f}(re^{i\theta})=0$, and denote by $(\theta',[\tilde{g}])$ the corresponding point in the Legendre transform. Note that for $f=z^{m/n}$, we have $z\tilde{f}'(z) = \frac{m}{n}\tilde{f}(z)$, and hence in general one has
	$$ \lim_{|z|\to\infty} z\tilde{f}'(z) = \lim_{|z|\to\infty} \frac{m}{n} \tilde{f}(z),$$
    where the limit means approaching $\infty$ along a fixed path (the same on both sides of the equality). 
	Hence, we can conclude
	\begin{align*}
    \lim_{r'\to\infty} \real \tilde{g}(r'e^{i\theta'}) &= \lim_{r\to\infty}\real \tilde{g}(\tilde{f}'(re^{i\theta}))\\
    &= \lim_{r\to\infty} \left(\real \tilde{f}(re^{i\theta}) - \frac{m}{n}\real \tilde{f}(re^{i\theta}) \right) = 0.
    \end{align*}
	The first equality holds since $\tilde{g}$ is a continuous function (with values in $\C\sqcup\{\infty\}$) defined up to the boundary of the real blow-up space, and the sequences $r' e^{i\theta'}$ and $\tilde{f}'(re^{i\theta})$ approach the same point on the boundary for $r'\to\infty$ and $r\to\infty$, respectively.
\end{proof}

The following is a direct consequence of Lemma~\ref{lemma:LegendreS0} and its proof.
\begin{cor}
	If $\langle g\rangle$ is the Legendre transform of an exponent $\langle f\rangle$ of slope $>1$, then the Legendre transform induces a bijection between $A(f)$ and $A(g)$. Moreover, this bijection sends $A_+(f)$ to $A_-(g)$ and $A_-(f)$ to $A_+(g)$.
\end{cor}
It is not difficult to see that this corresponds to the maps $\nu^\pm_0$ from \cite[§5.4.4]{Moc21}. To be more precise, if $J\in A_-(g)$ and $I\in A_+(f)$ is the corresponding sector via Legendre transform, we have $\pi_f(I)=\nu^+_0(\pi_g(J))$ as sectors on $\partial$.

\subsection{The algorithm for computing the Fourier transform}

We are now ready to translate Mochizuki's description of the Stokes data of the Fourier transform into the language of Stokes local systems. 

Let $(E,\nabla)$ be an algebraic connection on the affine line and $(\widehat E,\widehat\nabla)$ its Fourier transform\footnote{To be consistent with \cite{Moc21} (see §1.2.1 therein), we consider the Fourier transform given by the automorphism of the Weyl algebra $\C[z]\langle\partial_z\rangle\cong\C[w]\langle\partial_w\rangle$ where by $z\mapsto \partial_w,\partial_z\mapsto -w$.}. Let $\mathbb V$ be the Stokes local system associated to $(E,\nabla)$ via the Riemann--Hilbert--Birkhoff correspondence (describing flat sections of $(E,\nabla)$), and by $\widehat{\mathbb V}$ the one of $(\widehat E,\widehat\nabla)$. We assume that the irregular class $\Theta$ of $\mathbb{V}$ satisfies Assumption \ref{assumption}.

If $\mathbb V$ has irregular class $\Theta$ and associated formal local system $V^0$ with graded pieces $V^0_{\langle f\rangle} $ for the active circles $\langle f\rangle$, then the $V^0_{\langle f\rangle}$ can be viewed as local systems on the circles $\langle f\rangle$. Therefore, they are determined by vector spaces $K^{f}_I\defeq \Gamma(I,V^0_{\langle f\rangle})$ for any $I\in A(f)$ and identifications $\kappa_I\colon K^{f}_I \overset{\sim}{\to} K^f_J$ for any neighbouring sectors $I,J\in A(f)$, where $J$ follows $I$ with respect to the given orientation. Note that $K^f_I$ is naturally identified with $(V^0_{\langle f\rangle})_{m_I}$ for the midpoint $m_I$ of $I$.

Analogously, the formal local system associated with $\widehat{\mathbb V}$ has graded pieces $\widehat V^0_{\langle \widehat{f}\rangle}$ for the active circles $\langle \widehat{f}\rangle$ of $\widehat{\Theta}$, and  they are determined by vector spaces $\widehat{K}^{\widehat{f}}_{\widehat{I}}\defeq \Gamma(\widehat{I},\widehat V^0_{\langle \widehat{f}\rangle})$ for any $\widehat{I}\in A(\widehat{f})$ and identifications $\widehat{\kappa}_{\widehat{I}}\colon \widehat{K}^{\widehat{f}}_{\widehat{I}} \overset{\sim}{\to} \widehat{K}^{\widehat{f}}_{\widehat{J}}$ for neighbouring sectors $\widehat{I},\widehat{J}\in A(\widehat{f})$ with $\widehat{J}$ following $\widehat{I}$.

For $\langle f\rangle\in\cI$, we will denote by $\langle \widehat{f}\rangle$ its Legendre transform and denote by $\ell_f\colon \langle f\rangle \overset{\sim}{\to} \langle \widehat{f}\rangle$ the homeomorphism of circles defined by the Legendre transform.
Therefore, the active circles of $\widehat{\mathbb V}$ are the $\langle \widehat{f}\rangle$ for $\langle f\rangle\subset\cI$ an active circle of $\mathbb{V}$.

The following theorem is a reformulation, in our situation and language set up above, of the main result of \cite{Moc21} for the Fourier transform from $\infty$ to $\infty$ (see §8.9 and §10.5 in loc.~cit.). Given $\mathbb V$, we can express the data as in Proposition~\ref{prop:SLSdefdata} for $\widehat{\mathbb V}$, which in turn determines $\widehat{\mathbb V}$.

\begin{thm}\label{thm:trafoRule}
In the above situation, the Stokes local system $\widehat{\mathbb V}$ is the Stokes local system determined by the following data:
	\begin{itemize}
		\item The data encoding the local system $\widehat V^0_{\langle \widehat{f}\rangle}$ on $\langle \widehat{f}\rangle$ are obtained as follows: \\
		For any $\widehat{I}\in A(\widehat{f})$, we set
		$$\widehat{K}^{\widehat{f}}_{\widehat{I}} \defeq K^f_{\ell_f^{-1}(\widehat{I})}$$
		and for neighbouring sectors $\widehat{I},\widehat{J}\in A(\widehat{f})$, $\widehat{J}$ following $\widehat{I}$, the isomorphism
		$$\widehat{\kappa}_{\widehat{I}}\colon \widehat{K}^{\widehat{f}}_{\widehat{I}} \overset{\sim}{\to} \widehat{K}^{\widehat{f}}_{\widehat{J}}$$
		is given by
		\begin{itemize}
			\item[(1)] $\kappa_{\ell^{-1}(\widehat{I})}\colon K^f_{\ell_f^{-1}(\widehat{I})} \overset{\sim}{\to} K^f_{\ell_f^{-1}(\widehat{J})}$ if $\widehat{I}\in A_-(\widehat{f})$ and $\widehat{J}\in A_+(\widehat{f})$,
			\item[(2)] $(-1)$ times the identification $\kappa_{\ell^{-1}(\widehat{I})}\colon K^f_{\ell_f^{-1}(\widehat{I})} \overset{\sim}{\to} K^f_{\ell_f^{-1}(\widehat{J})}$ otherwise.
		\end{itemize}
		
		\item For any Stokes arrow $\widehat{d}_{\widehat{f}}\to_{\widehat{d}} \widehat{d}_{\widehat{g}}$ over a singular direction $\widehat{d}$ of $\widehat{\Theta}$, denote by $\widehat{I}\in A(\widehat{f})$ and $\widehat{J}\in A(\widehat{g})$ the sectors whose intersection has $\widehat{d}$ as a midpoint. Then the deformation datum
		$$\widehat{R}^{m_{\widehat{I}}}_{m_{\widehat{J}}}\colon (\widehat V^0_{\langle f\rangle})_{m_{\widehat{I}}} \to (\widehat V^0_{\langle g\rangle})_{m_{\widehat{J}}}$$
		is obtained by considering $I=\ell_f^{-1}(\widehat{I})$ and $J=\ell_g^{-1}(\widehat{J})$ and taking the linear map given by the composition
		$$(\widehat V^0_{\langle f\rangle})_{m_I}\overset{\iota_{I}}{\hookrightarrow} \mathbb V_{m_I} \overset{\sim}{\to} \mathbb V_{m_J}  \overset{\pi_J}{\twoheadrightarrow} (\widehat V^0_{\langle g\rangle})_{m_J}$$
		if $\widehat{J}$ comes before $\widehat{I}$, and $(-1)$ times this morphism otherwise (in particular if $\pi_f(\widehat{I})=\pi_g(\widehat{J})$).
		\end{itemize}
	\end{thm}
		
		More intuitively, the second part of the above theorem means the following: To obtain a deformation datum  of $\widehat{\mathbb V}$, map the corresponding path in the fission surface of $\widehat{\Theta}$ (cf.\ Remark~\ref{rem:defdataPath}) to the fission surface of $\Theta$ via the Legendre transform and read off the corresponding map given by parallel transport in $\mathbb V$ (up to a sign). Note that the path from $m_I$ to $m_J$ in the picture before Fourier transform that we obtain is not a path associated to a Stokes arrow of $\mathbb V$. We will make this explicit in many examples in Section~\ref{sec:examples}.
		
	\begin{proof}
		With the relations developed in the previous subsections, this is the statement of the results in \cite[§10.5]{Moc21} (see also \cite[§8.9]{Moc21} for more explicit descriptions of the Stokes shells) in our case.
		
		As remarked above, one observes that the objects $\mathcal{K}_\lambda$ from loc.~cit.\ correspond to the graded pieces of the formal local system and that the notion of deformation data from loc.~cit.\ correspond to that for Stokes local systems introduced above, as shown in Proposition~\ref{prop:defDataTransitionIJ}.
		
		A major difference in our formulation is that we consider everything over the circle $\partial$, instead of $2\pi\Z$-equivariant objects on the universal covering space $\R$. This does not pose any problem in our case since the only singularity of our system is at the point $\infty$, hence the topological monodromy (the monodromy of the Stokes local system around the origin) is equal to the identity. Therefore, one can identify an interval on $\partial$ with its copy shifted by $2\pi$ without getting a contribution in the formulas for the deformation data in loc.~cit.
		
		Moreover, a priori, the formulae in loc.~cit.\ describe the Fourier transform of a system localised (or co-localised) at $0$, while there is no singularity at $0$ in our case. However, the Fourier transform of the system without a singular point at $0$ can be constructed from these systems. This is a special case of a procedure called ``extension of local systems with Stokes structure'' in \cite{Moc21}, and here it amounts to ``forgetting'' (i.e.\ setting to zero) the graded piece for the exponent $\langle 0\rangle$ and all the associated deformation data.
	\end{proof}

    \subsection{Isomorphisms of wild character varieties}\label{sec:isomWCV}
	
	Note that Theorem~\ref{thm:trafoRule} gives us a completely topological way of passing from the data of the Stokes local system $\mathbb{V}$ (given, for example, by monodromies as in \S \ref{subsection:wild_character_varieties}) to the formal data and deformation data of $\widehat{\mathbb{V}}$. On the other hand, we know how these deformation data are extracted from $\widehat{\mathbb{V}}$ (see Definition~\ref{def:defDataSLS}) and that the formal data and deformation data together determine $\widehat{\mathbb{V}}$ (see Proposition~\ref{prop:SLSdefdata}), which leads to an explicit way of reconstructing a description of $\widehat{\mathbb{V}}$ in terms of wild monodromies as in \S \ref{subsection:wild_character_varieties} from these data. Altogether, this gives a completely explicit algorithmic way of computing the Stokes data of the Fourier transform in terms of the Stokes data of the original system, with Stokes data represented by wild monodromies as in Figure~\ref{fig:SLS}. 
    
    To see this in a precise way, we have to consider (minimal) framings for the Stokes local systems, since a choice of framing is necessary to write down explicit Stokes matrices. The main observation is that given a minimal framing for $\mathbb V$, the algorithm then yields a preferred minimal framing for $\widehat{\mathbb V}$. To show this, it will be helpful to introduce several slightly different notions of framings, building on our discussion in §\ref{subsection:wild_character_varieties}.

\begin{defi} Let $\Theta$ be an irregular class with active circles $\cir{q_1},\ldots,\cir{q_k}$, and for $l\in\{1,\ldots,k\}$, set $n_l\defeq \Theta(\cir{q_l})$ and $r_l\defeq\ram(q_l)$.
\begin{itemize}
\item A basepoint datum for $\Theta$ is a $k$-tuple $\mathbf b=(b_1, \dots, b_k)$, where $b_l\in \cir{q_l}$ for $l\in\{1,\ldots,k\}$.
\item Let $V^0\to\partial$ be a formal local system with irregular class $\Theta$.
\begin{itemize}
   \item  If $\mathbf b$ is a basepoint datum for $\Theta$, and $d_l\defeq\pi(b_l)\in \partial$ for any $l\in\{1,\ldots,k\}$, a \emph{minimal framing} of $V^0$ at $\mathbf b$ is a $k$-tuple $\bm{\phi}=(\phi_{b_1}, \dots, \phi_{b_k})$, where $\phi_{b_l}\colon V^0_{d_l}(b_l)\overset{\sim}\rightarrow \mathbb C^{n_l}$ is an isomorphism.
   \item If $b\in \partial$ is a direction, a \emph{framing}  of $V^0$ at $b$ is an isomorphism $\phi_b\colon V^0_b\overset{\sim} {\rightarrow}\mathbb C^n$, with $n=\sum_{l=1}^k n_l r_l$. The framing $\phi_b$ is said to be \emph {minimal} if in the matrix representing the monodromy $\rho\in \Aut(V^0_b)$ of $V^0_b$ via $\phi_b$, each $\cir{q_l}$-block is of the form \eqref{eq:monodromy_minimal_framing}, for all $l\in\{1,\ldots,k\}$.
\end{itemize}
\end{itemize}
\end{defi}

It is possible to pass between minimal framings at a basepoint datum $\mathbf b$ and minimal framings at a direction $b\in\partial$. Indeed, if $b\in \partial$, from any pair $(\mathbf b, \bm \phi)$, where $\mathbf b$ is a basepoint datum for $V^0$ and $\bm\phi$ is a minimal framing of $V^0$ at $\mathbf b$, we can construct a framing $\tau_b(\mathbf b, \bm \phi)$ of $V^0$ at $b$ in the following way: For any $l\in\{1,\ldots,k\}$, it is clear that there exists a unique isomorphism $\phi_{l}\colon V^0_{\cir{q_l},d_l}\overset{\sim}\rightarrow \mathbb C^{n_l r_l}$ extending $\phi_{b_l}$ to $V^0_{\cir{q_l},d_l}$, such that the matrix representing the monodromy $\rho_{\cir{q_l}}\in \Aut(V^0_{\cir{q_l},d_l}$) via $\phi_{l}$ and the decomposition $V^0_{\cir{q_l}, d_l}\cong \bigoplus_{j=0}^{r_l-1}V^0_{d_l}(\rho^{j}(b_l))$ is of the form \eqref{eq:monodromy_minimal_framing}.
In turn, we define the framing $\tau_b(\mathbf b, \bm \phi) $ of $V^0$ at $b$ by setting
\[
\tau_b(\mathbf b, \bm \phi)|_{V^0_{\cir{q_l}, b}}\defeq\phi_{l}\circ\rho_{b\to d_l}^{\cir{q_l}}, \qquad \text{ for } l\in\{1,\ldots,k\},
\]
where $\rho_{b\to d_l}^{\cir{q_l}}\colon V^0_b\to V^0_{d_l}$ denotes the parallel transport in $V^0_{\cir{q_l}}$ from $b$ to $d_l$ in the positive direction. We say that $\tau_b(\mathbf b, \bm \phi)$ is the framing of $V^0$  at $b$ associated to the minimal framing $\bm \phi$ at $\mathbf b$. 

It is straightforward to see that, via this process, the notions of minimal framing at a basepoint datum $\mathbf b$ and at a direction $b$ are equivalent:

\begin{lemma}
Let $\Theta$ be an irregular class, and $V^0\to\partial$ a formal local system with irregular class $\Theta$. Let $b\in\partial$ be a point. Then, if $\mathbf b$ is a basepoint datum for $\Theta$ and $\bm\phi$ is a minimal framing of $V^0$ at $\mathbf b$, then the framing $\tau_b(\mathbf b, \bm \phi)$ of $V^0$ at $b$ is minimal. Conversely, if $\phi_b$ is a minimal framing of $V^0$ at $b$, then for any basepoint datum $\mathbf b$ for $\Theta$ there exists a unique minimal framing $\bm \phi$  of $V^0$ at $\mathbf b$ such that $\phi_b=\tau_b(\mathbf b, \bm \phi)$.
\end{lemma}

Notice that, given a choice of basepoint datum $\mathbf b$ for $\Theta$, the set of minimal framings of $V^0$ at $\mathbf b$ is a torsor for the group $\breve H=\prod_{i=1}^k \mathrm{GL}_{n_l}(\mathbb C)$, and, similarly, for any $b\in\partial$, the set of minimal framings of $V^0$ at $b$ is a torsor for $\breve H$. Furthermore the map $\bm{\phi}\mapsto \tau_b(\mathbf b,\bm\phi)$ is equivariant with respect to the $\breve H$-action on both sides.\\

Let us now discuss the interplay of the Fourier transform of Stokes local systems with choices of minimal framings. Let $\Theta$ be an irregular class, and $\widehat{\Theta}$ the corresponding irregular class via the Fourier transform (i.e.\ the Legendre transform of $\Theta$). Let us fix choices of reference directions $b\in \partial$, $\widehat{b}\in \widehat{\partial}$, and a basepoint datum $\mathbf b$ for $\Theta$, which we assume to be \emph{generic}, i.e.\ such that, keeping previous notations, for any $l\in\{1,\ldots,k\}$, we have $b_l\notin S_0(q_l)$ (that is, the point $b_l$ is not a boundary point of a distinguished interval for $\cir{q_l}$, cf.\ the notation in §\ref{sec:distIntervals}). Applying the Legendre transform (which gives an isomorphism of circles, as discussed in §\ref{sec:Legendre}), we set $\widehat{\mathbf b}\defeq(\widehat{b}_1, \dots, \widehat{b}_k)$ with $\widehat{b}_l\defeq\ell(b_l)$. This defines a basepoint datum for $\widehat{\Theta}$. (In practice it will often be convenient to choose $\mathbf b$ such that $\pi(b_1)=\ldots=\pi(b_l)=\vcentcolon b\in \partial$.)

Now, it follows from §\ref{subsection:wild_character_varieties} that a pair $(\mathbb V, \phi_b)$, where $\mathbb V$ is a Stokes local system with irregular class $\Theta$, and $\phi_b$ is a minimal framing at $b$ of the corresponding formal local system $V^0$, is identified with a point $\rho\in\mathcal E_\Theta$ of the reduced representation variety $\mathcal E_\Theta$. Similarly, a pair $(\widehat{\mathbb V}, \widehat{\phi}_{\widehat{b}})$, where $\mathbb V$ is a Stokes local system with irregular class $\widehat{\Theta}$, and $\widehat{\phi}_{\widehat{b}}$ is a minimal framing at $\widehat{b}$ of the corresponding formal local system $\widehat{V}^0$ is identified with a point $\widehat \rho$ of the reduced representation variety $\mathcal E_{\widehat\Theta}$.

Using this, from the algorithm we obtain a map $\Phi\colon\mathcal E_\Theta\to \mathcal E_{\widehat\Theta}$ as follows: 

Let $\rho\in \mathcal E_\Theta$ be a reduced Stokes representation. It corresponds to a pair $(\mathbb V, \phi_b)$ as above, and there exists a unique minimal framing $\bm \phi$  of $V^0$ at $\mathbf b$ such that $\phi_b=\tau_b(\mathbf b, \bm \phi)$. Now, let $\widehat{\mathbb V}$ be the Fourier transform of $\mathbb V$, and $\widehat{V}^0$ its formal local system.

Consider $l\in\{1, \dots, k\}$. Since we chose $\mathbf b$ generic, there exists a unique distinguished interval $I\in A(q_l)$ such that $b_l\in I$. In turn, $\widehat{b}_l\in \ell(I)=\vcentcolon\widehat{I}$ which is a distinguished interval for $\cir{\widehat{q}_l}\defeq\ell(\cir{q_l})$. The algorithm says that, keeping the notations of Theorem~\ref{thm:trafoRule}, we have an identification $K^{q_l}_I=\widehat K^{\widehat{q}_l}_{\widehat{I}}$, which induces an isomorphism 
\begin{equation}
\psi_l\colon V^0_{d_l}(b_l)\overset{\sim}{\rightarrow} \widehat{V}^0_{\widehat{d}_l}(\widehat{b}_l)
\label{eq:new_minimal_framing}
\end{equation}
because $V^0_{\cir{q_l}}|_I$ and $\widehat{V}^0_{\cir{\widehat{q}_l}}|_{\widehat{I}}$ are constant sheaves.
This allows us to define an isomorphism 
\[
\widehat{\phi}_{\widehat{b}_l}\defeq\phi_{b_l}\circ\psi^{-1}_l\colon  \widehat{V}^0_{\widehat{d}_l}(\widehat{b}_l) \overset{\sim}{\rightarrow} \mathbb C^{n_l}.
\]
The $k$-tuple $\bm{\widehat\phi}\defeq(\widehat{\phi}_{\widehat{b}_1}, \dots, \widehat{\phi}_{\widehat{b}_k})$ is a minimal framing of $\widehat V^0$ at the basepoint datum $\mathbf{\widehat{b}}$, and we define $\widehat\phi_{\widehat b}\defeq\tau_{\widehat b}(\widehat b, \bm{\widehat \phi})$ as the corresponding minimal framing of $\widehat V^0$ at $\widehat b$. The pair $(\widehat{\mathbb V}, \widehat \phi_{\widehat b})$ is identified with a reduced Stokes representation $\widehat\rho\in\mathcal E_{\widehat \Theta}$, and we set
\[
\Phi(\rho)\defeq \widehat{\rho} \in \mathcal E_{\widehat\Theta}.
\]

\begin{thm}\label{thm:algebraic}
Let $\Theta$ be an irregular class at infinity satisfying Assumption \ref{assumption}, and $\widehat{\Theta}$ its formal Fourier (i.e.\ Legendre) transform. For any choices of directions $b\in\partial$, $\widehat b\in \widehat{\partial}$, and basepoint datum $\mathbf b$ for $\Theta$, the map $\Phi$ between the reduced representation varieties $\mathcal E_\Theta$ and $\mathcal E_{\widehat\Theta}$ is an algebraic isomorphism. Furthermore, the isomorphism is $\breve{H}$-equivariant with respect to the action of $\breve{H}$ on $\mathcal E_\Theta$ and $\mathcal E_{\widehat\Theta}$, hence induces an algebraic isomorphism between the wild character varieties $\mathcal M_B(\Theta)$ and $\mathcal M_B(\widehat\Theta)$.
\end{thm}

Note that for a given $\Theta$, while the notion of Stokes local system is completely intrinsic, the exact form of $\mathcal E_\Theta$ and $\mathcal E_{\widehat\Theta}$ and that of the isomorphism $\Phi$ depend on some choices (of $b$, $\widehat{b}$ and $\mathbf b$). This is not reflected in our above notation. However, it is clear that these objects, when constructed using different such choices, are canonically related by similar constructions as above.

\begin{proof}
The previous discussion shows that the map $\Phi$ is well-defined. That $\Phi$ is $\breve H$-equivariant follows directly from \eqref{eq:new_minimal_framing}, and from the fact that passing between minimal framings at a baspoint datum and at a direction is $\breve H$-equivariant. That $\Phi$ is a bijection follows immediately from the fact that the Fourier transform is invertible. Finally, that $\Phi$ is algebraic follows from the fact that, following the algorithm in Theorem~\ref{thm:trafoRule}, the entries of the Stokes matrices of the Fourier transform, corresponding via the choice of minimal framing to the monodromies of the Stokes local system $\mathbb V$ along the paths shown in Figure~\ref{fig:SLS} are -- up to a sign -- entries of products of the initial Stokes matrices and formal monodromy matrix, and in turn polynomial expressions in the entries of the initial Stokes data.
\end{proof}

\begin{rem}
    Since the algorithm is entirely topological, a natural question which arises at this point is whether there is a simpler, direct topological proof of the fact that $\Phi(\rho)$ is an element of $\mathcal E_{\widehat\Theta}$, which does not involve the Riemann--Hilbert--Birkhoff correspondence.
\end{rem}

Recall that the representation varieties $\mathcal E_\Theta$ and $\mathcal E_{\widehat\Theta}$ possess a quasi-Hamiltonian structure. We conjecture that these structures are preserved by the isomorphism induced by the Fourier transform.

\begin{conj}
The isomorphism $\Phi$ is compatible with the quasi-Hamiltonian structures on $\mathcal E_{\Theta}$ and $\mathcal E_{\widehat\Theta}$.
\end{conj}

This conjecture is quite natural: Indeed, it is known that the Fourier transform induces a symplectic isomorphism of wild character varieties in the case of connections on $\mathbb P^1$ with an unramified irregular singularity at infinity of order 2 together with regular singularities at finite distance \cite{Boa15a}. Furthermore, in \cite{Sza15} it is shown that in this case the Fourier transform preserves the full hyperkähler structure of the moduli space.

\begin{rem}
We also expect that, in the case where there is only one active exponent of the form $\cir{z^{s/r}}$ with $s>r$ and $r,s$ coprime, by considering the recessive or subdominant solutions as in \cite{Boa15b}, the wild character variety is isomorphic to a moduli space of configuration of points in a projective space, and that the Fourier transform of Stokes data is related to the Gale transform of point configurations. This  will be the object of future work. 
\end{rem}

\section{Examples}\label{sec:examples}

Is this section,  we consider a few examples to illustrate some fully explicit computations with Stokes matrices, using our algorithm. In particular, we will see that in all these examples the conjecture stated above is satisfied.

\subsection{The case of pure Gaussian type}

Consider first a system of pure Gaussian type, with irregular class 
\[
\Theta=\left\langle \frac{1}{2} z^2\right\rangle+\left\langle \frac{1+i}{2\sqrt{2}}z^2\right\rangle.
\]
(This short notation shall mean that the two summands are the only active circles, both of multiplicity one, i.e.\ mapped to $1$ by the irregular class $\Theta$.)

There are four singular directions, and hence any Stokes local system with irregular class $\Theta$ is determined by four Stokes factors and the formal monodromy. Explicitly, if we choose a base point $b\in\partial$ and a framing at $b$, the Stokes data of any connection with irregular class $\Theta$ are of the form:
\begin{align*}
	S_1 = \begin{pmatrix}1&0\\s_1&1\end{pmatrix}&& S_2 = \begin{pmatrix}1&s_2\\0&1\end{pmatrix}&& S_3 = \begin{pmatrix}1&0\\s_3&1\end{pmatrix}&& S_4 = \begin{pmatrix}1&s_4\\0&1\end{pmatrix} && h = \begin{pmatrix}\tau&0\\0&\tau'\end{pmatrix}
\end{align*}
and must satisfy $hS_4S_3S_2S_1=\id$. Let us denote by $\Sto_i$, $i=1,\dots, 4$ the Stokes groups corresponding to these forms of the Stokes matrices. Furthermore, let $H\subset \mathrm{GL}_2(\mathbb C)$ be the subgroup of diagonal matrices.

In turn, the representation variety associated to $\Theta$ is
\[
\cR_{\Theta}=\{(h, S_1,S_2, S_3,S_4)\in H\times \Sto_1\times \Sto_2\times \Sto_3\times \Sto_4 \mid hS_4S_3S_2S_1=\id\}.
\]
Any conjugacy class $\mathcal C\subset H$ is just a singleton $\mathcal C=\{h\}$ for some $h\in H$, so the wild character variety associated to the formal data $(\Theta, \mathcal C)$ corresponds to the quotient
\[
\mathcal M_B(\Theta,\mathcal C)=\{(S_1,S_2, S_3,S_4) \in \Sto_1\times \Sto_2\times \Sto_3\times \Sto_4\mid hS_4S_3S_2S_1=\id\}/H. 
\]

Since all Stokes matrices have determinant 1, it follows that $\mathcal M_B(\Theta,\mathcal C)$ can only be nonempty for $\mathcal{C}=\{h\}$ with $\det h=1$, i.e.\
\[h=\begin{pmatrix} \tau & 0\\ 0 & \tau^{-1}
\end{pmatrix},\]
for some $\tau\in \mathbb C^*$, which we assume to be the case in the rest of this example.

\begin{figure}
\centering
\begin{tikzpicture}[scale=3.4]
  \draw[dotted] (0,0) circle (1);  
  \begin{scope}[decoration={
  		markings,
  		mark=at position 0 with {\arrow{<}}}
  	]  
  	\draw [postaction={decorate},dotted] (0,0) circle (1);
  \end{scope}  
	\draw (0.9,0) node {$h$};
  \draw[blue,domain=0:(360),scale=1,samples=1000] plot (\x:{exp(cos(2*\x)/exp(1))});
   \draw[red,domain=0:(360),scale=1,samples=1000] plot (\x:{exp(cos(2*\x+45)/exp(1))});
   \foreach \n in {1,3} \draw[blue] ({-90*\n+90}:{exp(1.2*cos(-2*(90*\n+90))/exp(1))}) node {$J_{\n}$};
   \draw[blue] ({-90*2+90}:{exp(1.4*cos(-2*(90*2+90))/exp(1))}) node {$J_{2}$};
   \draw[blue] ({-90*4+90+5}:{exp(1.45*cos(-2*(90*4+90))/exp(1))}) node {$J_{4}$};
    \foreach \n in {1,3} \draw[red] ({-90*\n+90-45/2}:{exp(1.2*cos(2*(90*\n+90-45/2)+45)/exp(1))}) node {$J'_{\n}$};
    \draw[red] ({-90*2+90-45/2}:{exp(1.4*cos(2*(90*2+90-45/2)+45)/exp(1))}) node {$J'_{2}$};
    \draw[red] ({-90*4+90-45/2+5}:{exp(1.5*cos(2*(90*4+90-45/2)+45)/exp(1))}) node {$J'_{4}$};
    \draw[->,blue] ({45+5}:{1.06*exp(cos(2*(45+5))/exp(1))}) to ({45-5}:{1.06*exp(cos(2*(45-5))/exp(1))}) ; \draw[blue] (45+2:1.11) node {$\tau$};
    \draw[->,blue] ({-45+5}:{1.06*exp(cos(2*(-45+5))/exp(1))}) to ({-45-5}:{1.06*exp(cos(2*(-45-5))/exp(1))}) ; \draw[blue] (-45-2:1.14) node {$1$};
    \draw[->,blue] ({-135+5}:{1.06*exp(cos(2*(-135+5))/exp(1))}) to ({-135-5}:{1.06*exp(cos(2*(-135-5))/exp(1))}) ; \draw[blue] (-135+3:1.14) node {$1$};
    \draw[->,blue] ({135+5}:{1.06*exp(cos(2*(135+5))/exp(1))}) to ({135-5}:{1.06*exp(cos(2*(135-5))/exp(1))}) ; \draw[blue] (135-2:1.13) node {$1$};
   
   \draw[->,red] ({22.5+5}:{1.06*exp(cos(2*(22.5+5)+45)/exp(1))}) to ({22.5-5}:{1.06*exp(cos(2*(22.5-5)+45)/exp(1))}) ; \draw[red] (22.5:1.2) node {$\tau^{-1}$};
   \draw[->,red] ({-67.5+5}:{1.06*exp(cos(2*(-67.5+5)+45)/exp(1))}) to ({-67.5-5}:{1.06*exp(cos(2*(-67.5-5)+45)/exp(1))}) ; \draw[red] (-67.5-2:1.13) node {$1$};
   \draw[->,red] ({-157.5+5}:{1.06*exp(cos(2*(-157.5+5)+45)/exp(1))}) to ({-157.5-5}:{1.06*exp(cos(2*(-157.5-5)+45)/exp(1))}) ; \draw[red] (-157.5+3:1.12) node {$1$};
   \draw[->,red] ({112.5+5}:{1.06*exp(cos(2*(112.5+5)+45)/exp(1))}) to ({112.5-5}:{1.06*exp(cos(2*(112.5-5)+45)/exp(1))}) ; \draw[red] (112.5-3:1.12) node {$1$};

   \draw[magenta] (0,1) node {$\bullet$};
   \draw[magenta] (0,1.13) node {$b$};
   \foreach \n in {1,2,3,4}
{\draw ({45-\n*90-180/16}: 0.4) circle (0.04cm);
 \draw[magenta,->] (({45-\n*90-180/16}: 0.4)++(-0.13,0) arc (-180:180:0.13);
 \draw[magenta] ({135-\n*90-180/16}: 0.18) node {$S_{\n}$};
}
\begin{scope}[magenta,decoration={markings,mark=at position 0.75 with {\arrow{>}}}]  
 \draw[postaction={decorate}] (0,1) to [out angle=-70, in angle=80,curve through={(55:0.6)}]($({45-180/16}:0.4)+(80:0.13)$);
 \draw[postaction={decorate}] (0,1) to [out angle=-50, in angle=40,curve through={(0:0.55)}]($({-45-180/16}:0.4)+(40:0.13)$);
  \draw[postaction={decorate}] (0,1) to [out angle=-30, in angle=-40,curve through={(0:0.7)(-90:0.5)}]($({-135-180/16}:0.4)+(-40:0.13)$);
  \draw[postaction={decorate}] (0,1) to [out angle=-10, in angle=-150,curve through={(0:0.8)(-90:0.8)(180:0.55)}]($({135-180/16}:0.4)+(-150:0.13)$);
 \end{scope}
\end{tikzpicture}
\caption{The Stokes local system before Fourier transform. If we choose a base point $b$ at the boundary, it is determined by the monodromies around the punctures at the singular directions (given by the Stokes factors $S_i$) and the formal monodromy $h$ along the boundary circle $\partial$ (represented here by the dotted circle). The corresponding graded local system can be encoded by dividing every exponent circle into four intervals $J_i$ and $J'_i$, and associating a one-dimensional vector space to each interval, with gluing maps for passing from one to the next. We can choose bases on these intervals such that the graded local systems there are glued by three identity maps and a single nontrivial map.\\
Note that -- in contrast to Figure~\ref{fig:SLS} -- this picture shows the real blow-up of $\PP^1$ at $\infty$ from a different perspective: here, the affine line is in the interior (i.e.\ the origin is at the center of this picture).}
\label{fig:GaussSLSbefore}
\end{figure}
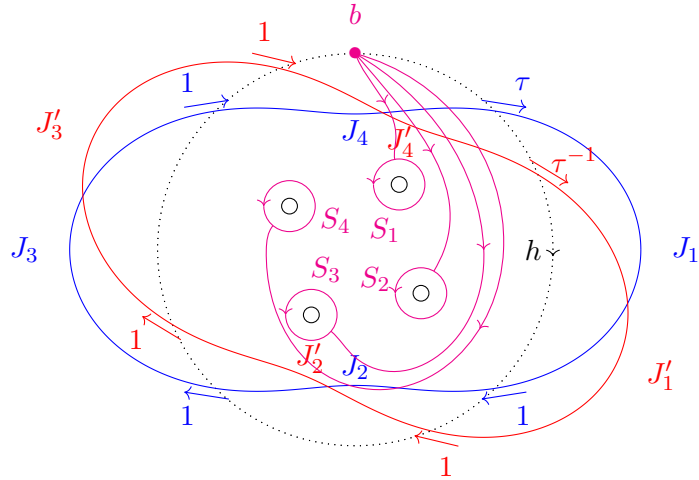

The initial situation and Stokes data, before applying the Fourier transformation, are illustrated in Figure~\ref{fig:GaussSLSbefore}.

Now, we want to determine the corresponding picture after Fourier transform. The formal Fourier transform of $\Theta$ is given (via the Legendre transform) by 
\[
\widehat{\Theta}=\left\langle -\frac{1}{2}w^2\right\rangle+\left\langle \frac{-1+i}{2\sqrt{2}}w^2\right\rangle.
\]
The space of Stokes representations and the wild character varieties with irregular class $\widehat\Theta$ have a similar form as those with irregular class $\Theta$.

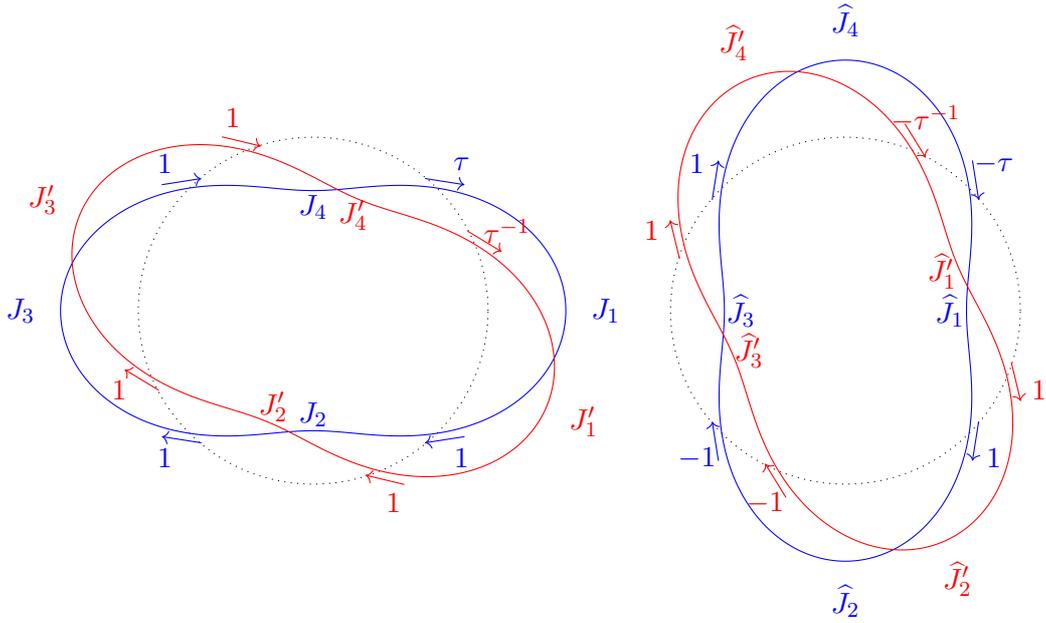
\begin{figure}
\centering
 \begin{tikzpicture}
     \begin{scope}[scale=2.3]
\draw[dotted] (0,0) circle (1);
  \draw[blue,domain=0:(360),scale=1,samples=1000] plot (\x:{exp(cos(2*\x)/exp(1))});
   \draw[red,domain=0:(360),scale=1,samples=1000] plot (\x:{exp(cos(2*\x+45)/exp(1))});
    \foreach \n in {1,3} \draw[blue] ({-90*\n+90}:{exp(1.2*cos(-2*(90*\n+90))/exp(1))}) node {$J_{\n}$};
   \draw[blue] ({-90*2+90}:{exp(1.4*cos(-2*(90*2+90))/exp(1))}) node {$J_{2}$};
   \draw[blue] ({-90*4+90+5}:{exp(1.45*cos(-2*(90*4+90))/exp(1))}) node {$J_{4}$};
    \foreach \n in {1,3} \draw[red] ({-90*\n+90-45/2}:{exp(1.2*cos(2*(90*\n+90-45/2)+45)/exp(1))}) node {$J'_{\n}$};
    \draw[red] ({-90*2+90-45/2}:{exp(1.4*cos(2*(90*2+90-45/2)+45)/exp(1))}) node {$J'_{2}$};
    \draw[red] ({-90*4+90-45/2+5}:{exp(1.5*cos(2*(90*4+90-45/2)+45)/exp(1))}) node {$J'_{4}$};
    \draw[->,blue] ({45+5}:{1.06*exp(cos(2*(45+5))/exp(1))}) to ({45-5}:{1.06*exp(cos(2*(45-5))/exp(1))}) ; \draw[blue] (45+2:1.11) node {$\tau$};
    \draw[->,blue] ({-45+5}:{1.06*exp(cos(2*(-45+5))/exp(1))}) to ({-45-5}:{1.06*exp(cos(2*(-45-5))/exp(1))}) ; \draw[blue] (-45-2:1.14) node {$1$};
    \draw[->,blue] ({-135+5}:{1.06*exp(cos(2*(-135+5))/exp(1))}) to ({-135-5}:{1.06*exp(cos(2*(-135-5))/exp(1))}) ; \draw[blue] (-135+3:1.14) node {$1$};
    \draw[->,blue] ({135+5}:{1.06*exp(cos(2*(135+5))/exp(1))}) to ({135-5}:{1.06*exp(cos(2*(135-5))/exp(1))}) ; \draw[blue] (135-2:1.13) node {$1$};
   \draw[->,red] ({22.5+5}:{1.06*exp(cos(2*(22.5+5)+45)/exp(1))}) to ({22.5-5}:{1.06*exp(cos(2*(22.5-5)+45)/exp(1))}) ; \draw[red] (22.5+1:1.21) node {$\tau^{-1}$};
   \draw[->,red] ({-67.5+5}:{1.06*exp(cos(2*(-67.5+5)+45)/exp(1))}) to ({-67.5-5}:{1.06*exp(cos(2*(-67.5-5)+45)/exp(1))}) ; \draw[red] (-67.5-2:1.13) node {$1$};
   \draw[->,red] ({-157.5+5}:{1.06*exp(cos(2*(-157.5+5)+45)/exp(1))}) to ({-157.5-5}:{1.06*exp(cos(2*(-157.5-5)+45)/exp(1))}) ; \draw[red] (-157.5+3:1.12) node {$1$};
   \draw[->,red] ({112.5+5}:{1.06*exp(cos(2*(112.5+5)+45)/exp(1))}) to ({112.5-5}:{1.06*exp(cos(2*(112.5-5)+45)/exp(1))}) ; \draw[red] (112.5-3:1.12) node {$1$};
\end{scope}
\begin{scope}[xshift=7cm,scale=2.3]
\draw[dotted] (0,0) circle (1);
  \draw[blue,domain=0:(360),scale=1,samples=1000] plot (\x:{exp(cos(2*\x+180)/exp(1))});
   \draw[red,domain=0:(360),scale=1,samples=1000] plot (\x:{exp(cos(2*\x+180-45)/exp(1))});
   \foreach \n in {1,3} \draw[blue] ({-90*\n+90}:{exp(1.4*cos(-2*(90*\n+90)+180)/exp(1))}) node {$\widehat{J}_{\n}$};
   \foreach \n in {2,4} \draw[blue] ({-90*\n+90}:{exp(1.27*cos(-2*(90*\n+90)+180)/exp(1))}) node {$\widehat{J}_{\n}$};
   
    \foreach \n in {1,3} \draw[red] ({-90*\n+90+45/2+8}:{exp(1.4*cos(2*(90*\n+90+45/2)-45+180)/exp(1))}) node {$\widehat{J}'_{\n}$};
    \foreach \n in {2,4} \draw[red] ({-90*\n+90+45/2}:{exp(1.27*cos(2*(90*\n+90+45/2)-45+180)/exp(1))}) node {$\widehat{J}'_{\n}$};
    \draw[->,blue] ({45+5}:{1.06*exp(cos(2*(45+5)+180)/exp(1))}) to ({45-5}:{1.06*exp(cos(2*(45-5)+180)/exp(1))}) ; \draw[blue] (45-4:1.17) node {$-\tau$};
    \draw[->,blue] ({-45+5}:{1.06*exp(cos(2*(-45+5)+180)/exp(1))}) to ({-45-5}:{1.06*exp(cos(2*(-45-5)+180)/exp(1))}) ; \draw[blue] (-45+3:1.12) node {$1$};
    \draw[->,blue] ({-135+5}:{1.06*exp(cos(2*(-135+5)+180)/exp(1))}) to ({-135-5}:{1.06*exp(cos(2*(-135-5)+180)/exp(1))}) ; \draw[blue] (-135-3:1.12) node {$1$};
    \draw[->,blue] ({135+5}:{1.06*exp(cos(2*(135+5)+180)/exp(1))}) to ({135-5}:{1.06*exp(cos(2*(135-5)+180)/exp(1))}) ; \draw[blue] (135+2:1.12) node {$1$};
    \draw[->,red] ({67.5+5}:{1.06*exp(cos(2*(67.5+5)+180-45)/exp(1))}) to ({67.5-5}:{1.06*exp(cos(2*(67.5-5)+180-45)/exp(1))}) ; \draw[red] (67.5-7:1.2) node {$-\tau^{-1}$};
    \draw[->,red] ({-22.5+5}:{1.06*exp(cos(2*(-22.5+5)+180-45)/exp(1))}) to ({-22.5-5}:{1.06*exp(cos(2*(-22.5-5)+180-45)/exp(1))}) ; \draw[red] (-22.5+3:1.12) node {$1$};
    \draw[->,red] ({-112.5+5}:{1.06*exp(cos(2*(-112.5+5)+180-45)/exp(1))}) to ({-112.5-5}:{1.06*exp(cos(2*(-112.5-5)+180-45)/exp(1))}) ; \draw[red] (-112.5-3:1.12) node {$1$};
    \draw[->,red] ({157.5+5}:{1.06*exp(cos(2*(157.5+5)+180-45)/exp(1))}) to ({157.5-5}:{1.06*exp(cos(2*(157.5-5)+180-45)/exp(1))}) ; \draw[red] (157.5+2:1.12) node {$1$};
\end{scope}
 \end{tikzpicture}
	\caption{The transformation of the formal monodromy: For each exponent, the four intervals before (shown on the left) and after Fourier transform (shown on the right) are in bijection via the Legendre transform (as indicated by the numbering). The gluing maps either remain the same (for transitions from a sector with negative exponent to one with a positive exponent in the Fourier transform) or change sign (otherwise). Of course, here this means that the local systems on the circles remain unchanged globally since we introduce an even number of sign changes, but the information about the local gluings will be important in the course of the algorithm.}
	\label{fig:GaussFormalFourier}
\end{figure}

Let us first describe the new formal data. Each of the active circles of $\widehat{\Theta}$ is again subdivided into four distinguished intervals. Via the Legendre transform, these intervals are identified with the ones before Fourier transform, and by the first part of the transformation rule of Proposition~\ref{thm:trafoRule}, the new graded local systems is encoded by having the same vector spaces on corresponding intervals, but introducing some signs in the gluing maps (see Figure~\ref{fig:GaussFormalFourier}).

In particular, the conjugacy class of the new formal monodromy is $\widehat{\mathcal C}=\mathcal C$.

Next, let us compute the deformation data of the Fourier transform using the second part of Proposition~\ref{thm:trafoRule}.
Therefore, we have to do the following: For any Stokes arrow associated to $\widehat{\Theta}$, we draw the corresponding path in the irregular curve of the Fourier transform, we transform it (by applying the Legendre transform to its endpoints) to a path in the irregular curve of the original system, and we determine the map associated to it explicitly since the Stokes local system is explicitly given on the left-hand side.

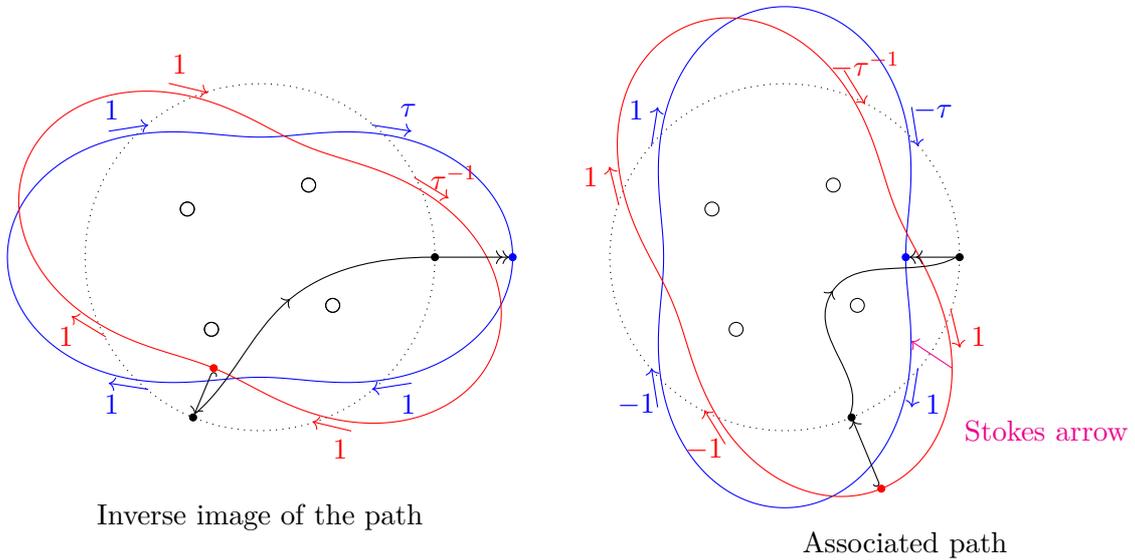
\begin{figure}
	\centering
	\begin{tikzpicture}
\begin{scope}[scale=2.3]
\draw[dotted] (0,0) circle (1);
  \draw[blue,domain=0:(360),scale=1,samples=1000] plot (\x:{exp(cos(2*\x)/exp(1))});
   \draw[red,domain=0:(360),scale=1,samples=1000] plot (\x:{exp(cos(2*\x+45)/exp(1))});
   \foreach \n in {1,2,3,4} 
\draw ({45-\n*90-180/16}: 0.5) circle (0.04cm);
    \draw[->,blue] ({45+5}:{1.06*exp(cos(2*(45+5))/exp(1))}) to ({45-5}:{1.06*exp(cos(2*(45-5))/exp(1))}) ; \draw[blue] (45+2:1.11) node {$\tau$};
    \draw[->,blue] ({-45+5}:{1.06*exp(cos(2*(-45+5))/exp(1))}) to ({-45-5}:{1.06*exp(cos(2*(-45-5))/exp(1))}) ; \draw[blue] (-45-2:1.14) node {$1$};
    \draw[->,blue] ({-135+5}:{1.06*exp(cos(2*(-135+5))/exp(1))}) to ({-135-5}:{1.06*exp(cos(2*(-135-5))/exp(1))}) ; \draw[blue] (-135+3:1.14) node {$1$};
    \draw[->,blue] ({135+5}:{1.06*exp(cos(2*(135+5))/exp(1))}) to ({135-5}:{1.06*exp(cos(2*(135-5))/exp(1))}) ; \draw[blue] (135-2:1.13) node {$1$};
   \draw[->,red] ({22.5+5}:{1.06*exp(cos(2*(22.5+5)+45)/exp(1))}) to ({22.5-5}:{1.06*exp(cos(2*(22.5-5)+45)/exp(1))}) ; \draw[red] (22.5+1:1.21) node {$\tau^{-1}$};
   \draw[->,red] ({-67.5+5}:{1.06*exp(cos(2*(-67.5+5)+45)/exp(1))}) to ({-67.5-5}:{1.06*exp(cos(2*(-67.5-5)+45)/exp(1))}) ; \draw[red] (-67.5-2:1.13) node {$1$};
   \draw[->,red] ({-157.5+5}:{1.06*exp(cos(2*(-157.5+5)+45)/exp(1))}) to ({-157.5-5}:{1.06*exp(cos(2*(-157.5-5)+45)/exp(1))}) ; \draw[red] (-157.5+3:1.12) node {$1$};
   \draw[->,red] ({112.5+5}:{1.06*exp(cos(2*(112.5+5)+45)/exp(1))}) to ({112.5-5}:{1.06*exp(cos(2*(112.5-5)+45)/exp(1))}) ; \draw[red] (112.5-3:1.12) node {$1$};
   
   \draw[red, fill=red] (-112.5:{exp(-1/exp(1))}) circle (0.02cm);
   \draw[fill=black] (-112.5:1) circle (0.02cm);
    \draw[fill=black] (0:1) circle (0.02cm);
    \draw[blue, fill=blue] (0:{exp(1/exp(1))}) circle (0.02cm);
   \draw[{Hooks[round,right]}->] (-112.5:{exp(-1/exp(1))+0.02}) to (-112.5:0.98); 
   \draw[decoration={markings,mark=at position 0.5 with {\arrow{>}}},postaction={decorate}] (-112.5:0.98) to [out angle=40, in angle=180,curve through={(-60:0.3)}](0:0.98);
   \draw[->>] (0:1.02) to (0:{exp(1/exp(1))-0.02});
   \draw (-90:1.5) node {Inverse image of the path};
\end{scope}
\begin{scope}[scale=2.3, xshift=3cm]
\draw[dotted] (0,0) circle (1);
  \draw[blue,domain=0:(360),scale=1,samples=1000] plot (\x:{exp(cos(2*\x+180)/exp(1))});
   \draw[red,domain=0:(360),scale=1,samples=1000] plot (\x:{exp(cos(2*\x+180-45)/exp(1))});  
\foreach \n in {1,2,3,4} 
\draw ({45-\n*90+180/16}: 0.5) circle (0.04cm);
    \draw[->,blue] ({45+5}:{1.06*exp(cos(2*(45+5)+180)/exp(1))}) to ({45-5}:{1.06*exp(cos(2*(45-5)+180)/exp(1))}) ; \draw[blue] (45-4:1.17) node {$-\tau$};
    \draw[->,blue] ({-45+5}:{1.06*exp(cos(2*(-45+5)+180)/exp(1))}) to ({-45-5}:{1.06*exp(cos(2*(-45-5)+180)/exp(1))}) ; \draw[blue] (-45+3:1.12) node {$1$};
    \draw[->,blue] ({-135+5}:{1.06*exp(cos(2*(-135+5)+180)/exp(1))}) to ({-135-5}:{1.06*exp(cos(2*(-135-5)+180)/exp(1))}) ; \draw[blue] (-135-3:1.12) node {$1$};
    \draw[->,blue] ({135+5}:{1.06*exp(cos(2*(135+5)+180)/exp(1))}) to ({135-5}:{1.06*exp(cos(2*(135-5)+180)/exp(1))}) ; \draw[blue] (135+2:1.12) node {$1$};
    \draw[->,red] ({67.5+5}:{1.06*exp(cos(2*(67.5+5)+180-45)/exp(1))}) to ({67.5-5}:{1.06*exp(cos(2*(67.5-5)+180-45)/exp(1))}) ; \draw[red] (67.5-7:1.2) node {$-\tau^{-1}$};
    \draw[->,red] ({-22.5+5}:{1.06*exp(cos(2*(-22.5+5)+180-45)/exp(1))}) to ({-22.5-5}:{1.06*exp(cos(2*(-22.5-5)+180-45)/exp(1))}) ; \draw[red] (-22.5+3:1.12) node {$1$};
    \draw[->,red] ({-112.5+5}:{1.06*exp(cos(2*(-112.5+5)+180-45)/exp(1))}) to ({-112.5-5}:{1.06*exp(cos(2*(-112.5-5)+180-45)/exp(1))}) ; \draw[red] (-112.5-3:1.12) node {$1$};
    \draw[->,red] ({157.5+5}:{1.06*exp(cos(2*(157.5+5)+180-45)/exp(1))}) to ({157.5-5}:{1.06*exp(cos(2*(157.5-5)+180-45)/exp(1))}) ; \draw[red] (157.5+2:1.12) node {$1$};
   
   \draw[red, fill=red] (-67.5:{exp(1/exp(1))}) circle (0.02cm);
   \draw[fill=black] (-67.5:1) circle (0.02cm);
    \draw[fill=black] (0:1) circle (0.02cm);
    \draw[blue, fill=blue] (0:{exp(-1/exp(1))}) circle (0.02cm);
   \draw[{Hooks[round,right]}->] (-67.5:{exp(1/exp(1))-0.02}) to (-67.5:1.02); 
   \draw[decoration={markings,mark=at position 0.5 with {\arrow{>}}},postaction={decorate}] (-67.5:0.98) to [out angle=70, in angle=-150,curve through={(-45:0.35)}](0:0.98);
   \draw[->>] (0:0.98) to (0:{exp(-1/exp(1))+0.02}); 
   \draw[->, magenta] ({-45+180/16}: {exp(cos(2*(-45+180/16)+180-45)/exp(1))}) to ({-45+180/16}: {exp(cos(2*(-45+180/16)+180)/exp(1))});
   \draw[magenta] ({-45+180/16+7}: 1.7) node {Stokes arrow};
   \draw  (-67.5:1.8) node  {Associated path};
\end{scope}
\end{tikzpicture}
\caption{Transporting the distinguished path associated to a Stokes arrow to the initial Stokes diagram.}
\label{fig:GaussStokesArrow}
\end{figure}

We illustrate this with one Stokes arrow (see Figure~\ref{fig:GaussStokesArrow}).
We draw the path associated to the Stokes arrow and pull it back to the picture before Fourier transform.

Now, comparing the picture on the left of Figure~\ref{fig:GaussStokesArrow} with Figure~\ref{fig:GaussSLSbefore}, we note that the parallel transport along this path in the initial Stokes local system (not taking into account the inclusion and projection for the moment) corresponds to
$$\begin{pmatrix}\tau&0\\0&\tau^{-1}\end{pmatrix}S_2^{-1}\begin{pmatrix}1\times \tau&0\\0&1\times \tau^{-1}\end{pmatrix}^{-1} = \begin{pmatrix}1&-s_2 \tau^2\\0&1\end{pmatrix}$$
(The first matrix is the part of the formal monodromy connecting the endpoint of the path with our basepoint $b$, and the third matrix connects the startpoint of our path with $b$.)
Of course, we use here -- and in all the following computations -- the fact that the parallel transport morphism only depends on the homotopy class of a path between two points.

Together with the inclusion of the red strand (the second exponent) at the beginning and the projection to the blue strand (the first exponent) at the end of the path, this means that the deformation datum associated to the Stokes arrow in Figure~\ref{fig:GaussSLSbefore} corresponds to the $(1,2)$ entry of this matrix, which is $-s_2 \tau^2$. If we do the same with the other deformation data (remembering that in cases where the orientation of the path is clockwise, we need to change the sign at the end), we get the full collection of deformation data for the Fourier transform, shown in Figure~\ref{fig:GaussFourDefData}.

\begin{figure}
	\centering
	\begin{tikzpicture}[scale=3.4]
	\draw[dotted] (0,0) circle (1);
	\draw[blue,domain=0:(360),scale=1,samples=1000] plot (\x:{exp(cos(2*\x+180)/exp(1))});
	\draw[red,domain=0:(360),scale=1,samples=1000] plot (\x:{exp(cos(2*\x+180-45)/exp(1))});  
	\foreach \n in {1,2,3,4} 
	\draw ({45-\n*90+180/16}: 0.5) circle (0.04cm);
    \draw[->,blue] ({45+5}:{1.06*exp(cos(2*(45+5)+180)/exp(1))}) to ({45-5}:{1.06*exp(cos(2*(45-5)+180)/exp(1))}) ; \draw[blue] (45-2:1.16) node {$-\tau$};
    \draw[->,blue] ({-45+5}:{1.06*exp(cos(2*(-45+5)+180)/exp(1))}) to ({-45-5}:{1.06*exp(cos(2*(-45-5)+180)/exp(1))}) ; \draw[blue] (-45+3:1.12) node {$1$};
    \draw[->,blue] ({-135+5}:{1.06*exp(cos(2*(-135+5)+180)/exp(1))}) to ({-135-5}:{1.06*exp(cos(2*(-135-5)+180)/exp(1))}) ; \draw[blue] (-135-3:1.12) node {$1$};
    \draw[->,blue] ({135+5}:{1.06*exp(cos(2*(135+5)+180)/exp(1))}) to ({135-5}:{1.06*exp(cos(2*(135-5)+180)/exp(1))}) ; \draw[blue] (135+2:1.12) node {$1$};
    \draw[->,red] ({67.5+5}:{1.06*exp(cos(2*(67.5+5)+180-45)/exp(1))}) to ({67.5-5}:{1.06*exp(cos(2*(67.5-5)+180-45)/exp(1))}) ; \draw[red] (67.5-5:1.19) node {$-\tau^{-1}$};
    \draw[->,red] ({-22.5+5}:{1.06*exp(cos(2*(-22.5+5)+180-45)/exp(1))}) to ({-22.5-5}:{1.06*exp(cos(2*(-22.5-5)+180-45)/exp(1))}) ; \draw[red] (-22.5+3:1.12) node {$1$};
    \draw[->,red] ({-112.5+5}:{1.06*exp(cos(2*(-112.5+5)+180-45)/exp(1))}) to ({-112.5-5}:{1.06*exp(cos(2*(-112.5-5)+180-45)/exp(1))}) ; \draw[red] (-112.5-3:1.12) node {$1$};
    \draw[->,red] ({157.5+5}:{1.06*exp(cos(2*(157.5+5)+180-45)/exp(1))}) to ({157.5-5}:{1.06*exp(cos(2*(157.5-5)+180-45)/exp(1))}) ; \draw[red] (157.5+2:1.12) node {$1$};
	
	\draw[red, fill=red] (-67.5:{exp(1/exp(1))}) circle (0.02cm);
	\draw[fill=black] (-67.5:1) circle (0.02cm);
	\draw[fill=black] (0:1) circle (0.02cm);
	\draw[blue, fill=blue] (0:{exp(-1/exp(1))}) circle (0.02cm);
	\draw[{Hooks[round,right]}->] (-67.5:{exp(1/exp(1))-0.02}) to (-67.5:1.02); 
	\draw[decoration={markings,mark=at position 0.5 with {\arrow{>}}},postaction={decorate}] (-67.5:0.98) to [out angle=70, in angle=-150,curve through={(-45:0.4)}](0:0.98);
	\draw[->>] (0:0.98) to (0:{exp(-1/exp(1))+0.02});

	\draw[blue, fill=blue] (-90:{exp(1/exp(1))}) circle (0.02cm);
	\draw[fill=black] (-90:1) circle (0.02cm);
	\draw[fill=black] (-157.5:1) circle (0.02cm);
	\draw[red, fill=red] (-157.5:{exp(-1/exp(1))}) circle (0.02cm);
	\draw[{Hooks[round,right]}->] (-90:{exp(1/exp(1))-0.02}) to (-90:1.02); 
	\draw[decoration={markings,mark=at position 0.5 with {\arrow{>}}},postaction={decorate}] (-90:0.98) to [out angle=110, in angle=-30,curve through={(-135:0.4)}](-157.5:0.98);
	\draw[->>] (-157.5:0.98) to (-157.5:{exp(-1/exp(1))+0.02});
	
	\draw[red, fill=red] (112.5:{exp(1/exp(1))}) circle (0.02cm);
	\draw[fill=black] (112.5:1) circle (0.02cm);
	\draw[fill=black] (180:1) circle (0.02cm);
	\draw[blue, fill=blue] (180:{exp(-1/exp(1))}) circle (0.02cm);
	\draw[{Hooks[round,right]}->] (112.5:{exp(1/exp(1))-0.02}) to (112.5:1.02); 
	\draw[decoration={markings,mark=at position 0.5 with {\arrow{>}}},postaction={decorate}] (112.5:0.98) to [out angle=-110, in angle=30,curve through={(135:0.4)}](180:0.98);
	\draw[->>] (180:0.98) to (180:{exp(-1/exp(1))+0.02});
	
	\draw[blue, fill=blue] (90:{exp(1/exp(1))}) circle (0.02cm);
	\draw[fill=black] (90:1) circle (0.02cm);
	\draw[fill=black] (22.5:1) circle (0.02cm);
	\draw[red, fill=red] (22.5:{exp(-1/exp(1))}) circle (0.02cm);
	\draw[{Hooks[round,right]}->] (90:{exp(1/exp(1))-0.02}) to (90:1.02); 
	\draw[decoration={markings,mark=at position 0.5 with {\arrow{>}}},postaction={decorate}] (90:0.98) to [out angle=-70, in angle=150,curve through={(45:0.4)}](22.5:0.98);
	\draw[->>] (22.5:0.98) to (22.5:{exp(-1/exp(1))+0.02});   
	
	\draw ({45+180/16}:0.28) node {$-s_1/\tau$};
	\draw ({-30+180/16}:0.3) node {$-s_2 \tau^2$};
	\draw ({-125+180/16}:0.26) node {$-s_3/\tau^2$};
	\draw ({145+180/16}:0.28) node {$-s_4 \tau^2$};
\end{tikzpicture}
	\caption{Deformation data of the Fourier transform.\\
    (Note that the fact that three of them involve $\tau^2$ and one only involves $\tau$ is not a mistake, but comes from the fact that the gluing maps for the formal local system between the endpoints of the Stokes paths are trivial in three cases and involve $\tau^{\pm 1}$ in the remaining case.)}
	\label{fig:GaussFourDefData}
\end{figure}

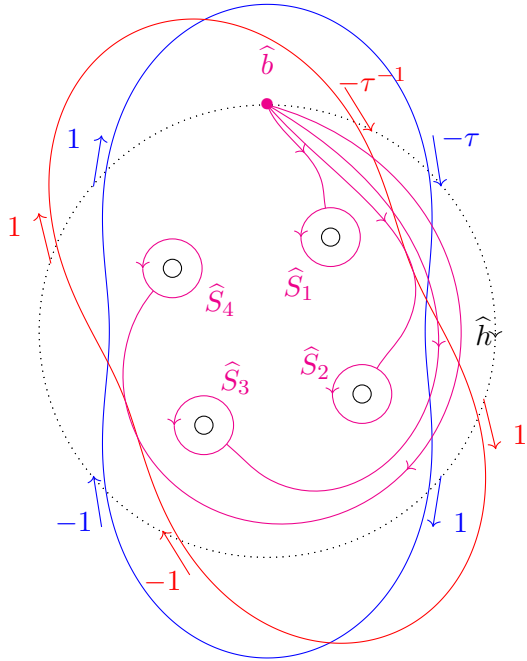
\begin{figure}
	\centering
	\begin{tikzpicture}[scale=3.4]
	\draw[dotted] (0,0) circle (1);
	\begin{scope}[decoration={
			markings,
			mark=at position 0 with {\arrow{<}}}
		]  
		\draw [postaction={decorate},dotted] (0,0) circle (1);
	\end{scope}  
	\draw (1.1,0) node {$\widehat{h}$};
	\draw[blue,domain=0:(360),scale=1,samples=1000] plot (\x:{exp(cos(2*\x+180)/exp(1))});
	\draw[red,domain=0:(360),scale=1,samples=1000] plot (\x:{exp(cos(2*\x+180-45)/exp(1))});  
    \draw[->,blue] ({45+5}:{1.06*exp(cos(2*(45+5)+180)/exp(1))}) to ({45-5}:{1.06*exp(cos(2*(45-5)+180)/exp(1))}) ; \draw[blue] (45-2:1.16) node {$-\tau$};
    \draw[->,blue] ({-45+5}:{1.06*exp(cos(2*(-45+5)+180)/exp(1))}) to ({-45-5}:{1.06*exp(cos(2*(-45-5)+180)/exp(1))}) ; \draw[blue] (-45+3:1.12) node {$1$};
    \draw[->,blue] ({-135+5}:{1.06*exp(cos(2*(-135+5)+180)/exp(1))}) to ({-135-5}:{1.06*exp(cos(2*(-135-5)+180)/exp(1))}) ; \draw[blue] (-135-3:1.12) node {$1$};
    \draw[->,blue] ({135+5}:{1.06*exp(cos(2*(135+5)+180)/exp(1))}) to ({135-5}:{1.06*exp(cos(2*(135-5)+180)/exp(1))}) ; \draw[blue] (135+2:1.12) node {$1$};
    \draw[->,red] ({67.5+5}:{1.06*exp(cos(2*(67.5+5)+180-45)/exp(1))}) to ({67.5-5}:{1.06*exp(cos(2*(67.5-5)+180-45)/exp(1))}) ; \draw[red] (67.5-5:1.19) node {$-\tau^{-1}$};
    \draw[->,red] ({-22.5+5}:{1.06*exp(cos(2*(-22.5+5)+180-45)/exp(1))}) to ({-22.5-5}:{1.06*exp(cos(2*(-22.5-5)+180-45)/exp(1))}) ; \draw[red] (-22.5+3:1.12) node {$1$};
    \draw[->,red] ({-112.5+5}:{1.06*exp(cos(2*(-112.5+5)+180-45)/exp(1))}) to ({-112.5-5}:{1.06*exp(cos(2*(-112.5-5)+180-45)/exp(1))}) ; \draw[red] (-112.5-3:1.12) node {$1$};
    \draw[->,red] ({157.5+5}:{1.06*exp(cos(2*(157.5+5)+180-45)/exp(1))}) to ({157.5-5}:{1.06*exp(cos(2*(157.5-5)+180-45)/exp(1))}) ; \draw[red] (157.5+2:1.12) node {$1$};
	
	\draw[magenta] (0,1) node {$\bullet$};
	\draw[magenta] (0,1.12) node {$\widehat{b}$};  
	\foreach \n in {1,2,3,4} 
	{\draw ({45-\n*90+180/16}: 0.5) circle (0.04cm);
		\draw[magenta] ({135-\n*90+180/16}: 0.25) node {$\widehat{S}_{\n}$};
		\draw[magenta, ->] (({45-\n*90+180/16}: 0.5)++(-0.13,0) arc (-180:180:0.13);
	}
	\begin{scope}[magenta,decoration={markings,mark=at position 0.5 with {\arrow{>}}}]  
		\draw[postaction={decorate}] (0,1) to [out angle=-80, in angle=100,curve through={(70:0.65)}]($({45+180/16}:0.5)+(100:0.13)$);
		\draw[postaction={decorate}] (0,1) to [out angle=-50, in angle=60,curve through={(0:0.57)}]($({-45+180/16}:0.5)+(60:0.13)$);
		\draw[postaction={decorate}] (0,1) to [out angle=-30, in angle=-40,curve through={(0:0.75)(-90:0.65)}]($({-135+180/16}:0.5)+(-40:0.13)$);
		\draw[postaction={decorate}] (0,1) to [out angle=-15, in angle=-130,curve through={(0:0.85)(-90:0.85)(180:0.57)}]($({135+180/16}:0.5)+(-130:0.13)$);
	\end{scope}
\end{tikzpicture}
	\caption{The monodromies of the Stokes local system of the Fourier transform with respect to a choice of basepoint $\widehat{b}$.}
	\label{fig:GaussFourSLS}
\end{figure}

Finally, to obtain an explicit Stokes matrix description of the Stokes local system of the Fourier transform, let us choose a basepoint $\widehat{b}$ at the boundary in the picture for the Fourier transform. The Stokes factors $\widehat{S}_i$ correspond to the monodromies along the paths in Figure~\ref{fig:GaussFourSLS}.

We compare Figure~\ref{fig:GaussFourDefData} and Figure~\ref{fig:GaussFourSLS} to determine the matrices $\widehat{S}_i$ from the deformation data. For example, consider the deformation datum with value $-s_3/\tau^2$. The parallel transport along the corresponding path (without the inclusion and projection) is also given by the matrix
$$\begin{pmatrix}(-1)\times 1 \times(-\tau) &0\\0&(-1)\times 1 \times(-\tau^{-1})\end{pmatrix}\widehat{S}_3\begin{pmatrix} 1 \times(-\tau)&0\\0& 1 \times(-\tau^{-1})\end{pmatrix}^{-1}.$$
The deformation datum itself, obtained by pre-composing with the inclusion of the first component and post-composing with the projection to the second component, hence corresponds to the $(2,1)$ component of this matrix.

Since we know a priori that $\widehat{S}_3$ is of the form $\begin{pmatrix}1&0\\\widehat s_3&1\end{pmatrix}$, it follows that $\widehat s_3=s_3$. Doing this for all deformation data paths, we get:
\begin{align*}
	\widehat{S}_1 = \begin{pmatrix}1&0\\s_1&1\end{pmatrix}&& \widehat{S}_2 = \begin{pmatrix}1&s_2\\0&1\end{pmatrix}&& \widehat{S}_3 = \begin{pmatrix}1&0\\s_3&1\end{pmatrix}&& \widehat{S}_4 = \begin{pmatrix}1&s_4\\0&1\end{pmatrix} && \widehat{h} = \begin{pmatrix}\tau&0\\0&\tau^{-1}\end{pmatrix}.
\end{align*}
We thus obtain, with this choice of parametrisation, that $\widehat{S}_i=S_i$ and $\widehat h=h$, that is, the Stokes factors and formal monodromy after Fourier transform are the same as the initial ones. In particular, we clearly have here a well-defined isomorphism between the spaces of Stokes representations on both sides, which in turn induces (through symplectic reduction by $H$) isomorphisms  at the level of the wild character varieties. The Poisson/symplectic structures are obviously preserved.

To summarise, we have established the following.
\begin{prop}
Keeping the previous notations, let $\mathcal C=\{h\}$ be a conjugacy class for $\Theta$ and $\widehat{\mathcal C}=\mathcal C$. The Fourier transform induces the isomorphism
$\Phi\colon \mathcal M_B(\Theta,\mathcal C) \to \mathcal M_B(\widehat\Theta, \widehat{\mathcal C}) $ given by 
\[
\widehat S_i=S_i,\qquad \widehat h=h
\]
with respect to the above choices of the base point and framing at this base point.
\end{prop}

\begin{rem}
	Let us remark that this result is consistent with previous studies of Fourier transforms of Stokes data for the Gaussian case \cite{Sab16,Ho22}: Although we consider different conditions on the exponents here, we still get a result in the same spirit: With respect to suitable choices, the representations of Stokes data on both sides of the Fourier transform coincide.
\end{rem}

\subsection{The Airy case}

Let us now discuss an example featuring a ramified irregular class.
Consider the irregular class $\Theta=\langle z^3/3\rangle$. Its Stokes structure is trivial (there are no Stokes arrows). In turn, a wild character variety with this irregular class is nonempty only if its formal monodromy is the identity, so there exists a unique rank one Stokes local system $\mathbb V$ with irregular class $\Theta$. Let us determine its Fourier transform. 

The formal Fourier transform of $\Theta$ is $\widehat{\Theta}=\langle \frac{2}{3}w^{3/2}\rangle$, which is closely related to the Airy differential equation.

\begin{figure}
	\centering
\begin{tikzpicture}[scale=2.3]
	\begin{scope}
		\draw[dotted] (0,0) circle (1);
		\draw[blue,domain=0:(360),scale=1,samples=1000] plot (\x:{exp(-cos(3*\x+180)/exp(1))});  
		\foreach \n in {1,2,...,6} \draw[blue] ({120-(\n+1)*60}:{exp(-cos(3*(120-(\n+1)*60)+180)/exp(1))-0.12}) node {$J_\n$};
		\foreach \n in {1,3,5} 
		{\draw[->,blue] ({90-\n*60+4}:{1.06*exp(cos(3*(90-(\n+1)*60+4)+180)/exp(1))}) to ({90-\n*60-4}:{1.06*exp(cos(3*(90-(\n+1)*60-4)+180)/exp(1))});
			\draw[blue] ({90-\n*60+4}:1.15) node {$1$};
		}
        \foreach \n in {2,4,6} 
		{\draw[->,blue] ({90-\n*60+4}:{1.06*exp(cos(3*(90-(\n+1)*60+4)+180)/exp(1))}) to ({90-\n*60-4}:{1.06*exp(cos(3*(90-(\n+1)*60-4)+180)/exp(1))});
			\draw[blue] ({90-\n*60-4}:1.15) node {$1$};
		} 
	\end{scope}
	
	\begin{scope}[xshift=3.5cm]
		\draw[dotted] (0,0) circle (1);
		\draw[blue,domain=0:(2*360),scale=1,samples=1000] plot (\x:{exp(cos(3/2*\x)/exp(1))});  
		\foreach \n in {1,3,5} \draw[blue] ({420-\n*120+60}:{exp(1.55*cos(3/2*(420-\n*120)+90)/exp(1))}) node {$\widehat J_\n$};
        \foreach \n in {2,4,6} \draw[blue] ({420-\n*120+60}:{exp(1.3*cos(3/2*(420-\n*120)+90)/exp(1))}) node {$\widehat J_\n$};
        \draw[blue,->] ({0*120+60+5}:0.75) to[bend left] ({0*120+60-5}:0.85);
			\draw[blue,->] ({0*120+60+5}:0.85) to[bend left] ({0*120+60-5}:0.75);
			\draw[blue] ({0*120+60-8}:0.87) node {{\tiny{$1$}}};
			\draw[blue] ({0*120+60-12+5}:0.7) node {{\tiny{$-1$}}};
        \draw[blue,->] ({1*120+60+5}:0.75) to[bend left] ({1*120+60-5}:0.85);
			\draw[blue,->] ({1*120+60+5}:0.85) to[bend left] ({1*120+60-5}:0.75);
			\draw[blue] ({1*120+60-8}:0.87) node {{\tiny{$1$}}};
			\draw[blue] ({1*120+60-12+4}:0.68) node {{\tiny{$-1$}}};
        \draw[blue,->] ({2*120+60+5}:0.75) to[bend left] ({2*120+60-5}:0.85);
			\draw[blue,->] ({2*120+60+5}:0.85) to[bend left] ({2*120+60-5}:0.75);
			\draw[blue] ({2*120+60-8}:0.87) node {{\tiny{$1$}}};
			\draw[blue] ({2*120+60-12}:0.73) node {{\tiny{$-1$}}};
        
		\draw[blue,->] ({60+5}:1.1) to[bend left] ({60-5}:1.1);
		\draw[blue] ({60+0-4}:1.4) node {{\tiny $\begin{pmatrix} 0 & 1 \\ -1 & 0\end{pmatrix}$}};
		\draw[blue,->] ({120+60+5}:1.1) to[bend left] ({120+60-5}:1.1);
		\draw[blue] ({120+60}:1.43) node {{\tiny $\begin{pmatrix} -1 & 0 \\ 0 & 1\end{pmatrix}$}};
		\draw[blue,->] ({-120+60+5}:1.1) to[bend left] ({-120+60-5}:1.1);
		\draw[blue] ({-120+60+5}:1.4) node {{\tiny $\begin{pmatrix} -1 & 0 \\ 0 & 1\end{pmatrix}$}};
		
		\draw({cos(35)},{sin(35)}) node {$\bullet$};
		\draw({cos(35)+0.1},{sin(35)+0.1}) node {$\widehat{b}$}; 
		\draw[teal] (0:{1.1*exp(-1/exp(1))}) node {\tiny{$2$}};
		\draw[teal] (0:{0.93*exp(1/exp(1))}) node {\tiny{$1$}};
		\draw[teal] (-120:{1.1*exp(-1/exp(1))}) node {\tiny{$1$}};
		\draw[teal] (-120:{0.93*exp(1/exp(1))}) node {\tiny{$2$}};
		\draw[teal] (120:{1.1*exp(-1/exp(1))}) node {\tiny{$2$}};
		\draw[teal] (120:{0.93*exp(1/exp(1))}) node {\tiny{$1$}};
	\end{scope}
\end{tikzpicture}
	\caption{On the left: The formal local system before Fourier transform.\\ On the right: The formal local system after Fourier transform.\\ The numbers in the circles show the identification of sectors given by Legendre transform. The green numbers show the numbering of the strands, given by choosing a numbering near a basepoint $\widehat{b}$ and continuing it compatibly in clockwise direction. The matrices are the partial formal monodromies given by the transformation rule (they come from the formal monodromies on the left, but change signs whenever one transits from a positive to a negative strand).}
	\label{fig:AiryFormalFourier}
\end{figure}
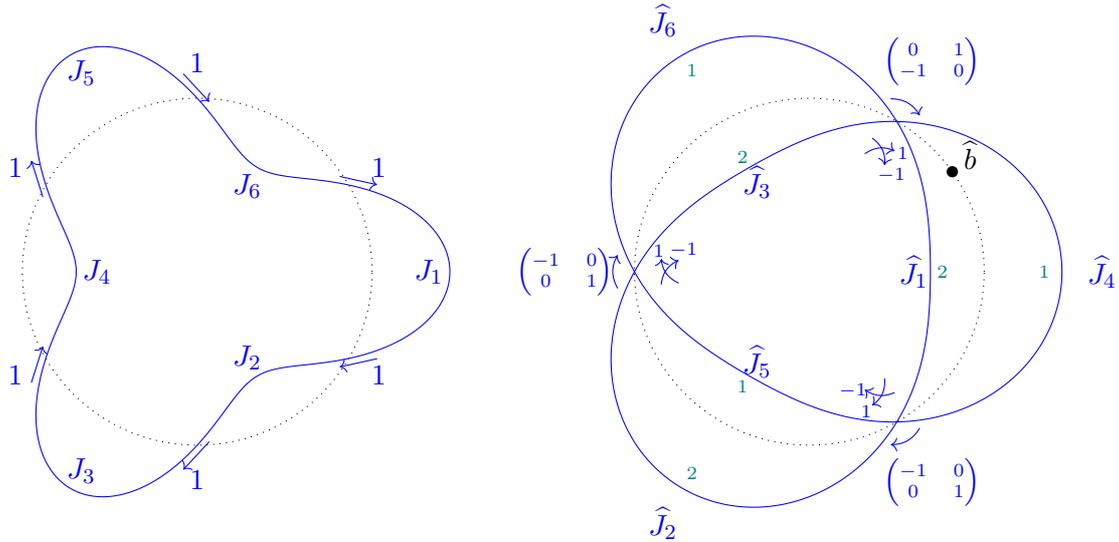

The formal local system associated to $\mathbb V$ is the local system on $\langle z^3/3\rangle$ with trivial monodromy, which we view again as glued from one-dimensional constant sheaves (with stalk $\C$) on 6 sectors. The Legendre transform maps the circle $\langle z^3/3\rangle$ onto $\langle \frac{2}{3}w^{3/2}\rangle$. This identifies the sectors on both sides, and the algorithm gives us the induced transition maps defining the formal monodromy of the new formal local system, as indicated in Figure~\ref{fig:AiryFormalFourier}: If we choose a basepoint and an ordering of the strands, we get three transition matrices (in blue) for the formal local system of the Fourier transform (considered as a rank $2$ local system on the boundary circle), the third one taking into account the twist (permutation of the strands).

There are three Stokes arrows in the Airy case (i.e.\ in the picture after Fourier transform). Each of them corresponds to a path on the Stokes diagram of $\widehat{\Theta}$. Via the Legendre transform, this path can be transported to  a path in the initial Stokes diagram, the parallel transport along which defines a linear map (a number in this explicit case). Let us detail this process in the example drawn in Figure~\ref{fig:AiryStokesPath}.

\begin{figure}
\centering
\begin{tikzpicture}[scale=2]
	\begin{scope}
		\draw[dotted] (0,0) circle (1);
		\draw[blue,domain=0:(360),scale=1,samples=1000] plot (\x:{exp(-cos(3*\x+180)/exp(1))}); 
		
		\draw[blue, fill=blue] (-60:{exp(-1/exp(1))}) circle (0.02cm);
		\draw[fill=black] (-60:1) circle (0.02cm);
		\draw[fill=black] (120:1) circle (0.02cm);
		\draw[blue, fill=blue] (120:{exp(1/exp(1))}) circle (0.02cm);
		
		\draw[{Hooks[round,right]}->] (-60:{exp(-1/exp(1))+0.02}) to (-60:0.98); 
		\draw[decoration={markings,mark=at position 0.5 with {\arrow{>}}},postaction={decorate}] (-60:0.98) to [out angle=110, in angle=-40,curve through={(0:0.3)}](120:0.98);
		\draw[->>] (120:0.98) to (120:{exp(1/exp(1))-0.02}); 
		
		\draw (0,-2) node {Inverse image of path};
	\end{scope}
	
	\begin{scope}[xshift=3cm]
		\draw[dotted] (0,0) circle (1);
		\draw[blue,domain=0:(2*360),scale=1,samples=1000] plot (\x:{exp(cos(3/2*\x)/exp(1))}); 
		
		\draw[blue, fill=blue] (-120:{exp(1/exp(1))}) circle (0.02cm);
		\draw[fill=black] (-120:1) circle (0.02cm);
		\draw[blue, fill=blue] (-120:{exp(-1/exp(1))}) circle (0.02cm);
		
		\foreach \th in {0,120,240} \draw (\th:0.5) circle (0.04cm);
		
		\draw[{Hooks[round,right]}->] (-120:{exp(1/exp(1))-0.02}) to (-120:1.02); 
		\draw[decoration={markings,mark=at position 0.5 with {\arrow{>}}},postaction={decorate}] (-120:0.98) to [out angle=20, in angle=100,curve through={(-120:0.4)}](-120:0.98);
		\draw[->>] (-120:1.02) to (-120:{exp(-1/exp(1))+0.02});  
		\draw (0,-2) node {Path associated to one Stokes arrow};
	\end{scope}
\end{tikzpicture}
	\caption{A Stokes arrow between a positive and a negative sector corresponds to a path connecting the midpoints of the two sectors, going around at least one singular direction (composed with the embedding of the positive strand and the projection to the negative strand). Via the Legendre transform it has a corresponding path on the left-hand side, in the picture before Fourier transform. Notice that this path does not need to correspond to a Stokes arrow.}
	\label{fig:AiryStokesPath}
\end{figure}

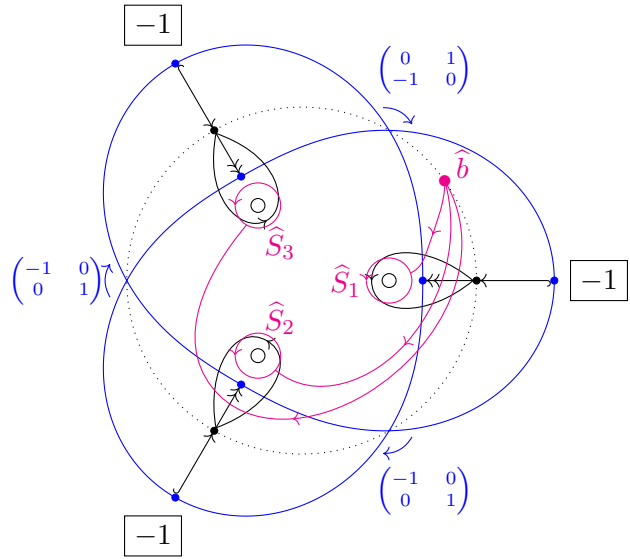
\begin{figure}
\centering
\begin{tikzpicture}[scale=3]
	\draw[dotted] (0,0) circle (1);
	\draw[blue,domain=0:(2*360),scale=1,samples=1000] plot (\x:{exp(cos(3/2*\x)/exp(1))});  
	
	\draw[blue,->] ({60+5}:1.1) to[bend left] ({60-5}:1.1);
		\draw[blue] ({60+0-4}:1.4) node {{\scriptsize $\begin{pmatrix} 0 & 1 \\ -1 & 0\end{pmatrix}$}};
		\draw[blue,->] ({120+60+5}:1.1) to[bend left] ({120+60-5}:1.1);
		\draw[blue] ({120+60}:1.43) node {{\scriptsize $\begin{pmatrix} -1 & 0 \\ 0 & 1\end{pmatrix}$}};
		\draw[blue,->] ({-120+60+5}:1.1) to[bend left] ({-120+60-5}:1.1);
		\draw[blue] ({-120+60+5}:1.4) node {{\scriptsize $\begin{pmatrix} -1 & 0 \\ 0 & 1\end{pmatrix}$}};
	
	\draw[blue, fill=blue] (0:{exp(1/exp(1))}) circle (0.02cm);
	\draw[fill=black] (0:1) circle (0.02cm);
	\draw[blue, fill=blue] (0:{exp(-1/exp(1))}) circle (0.02cm);
	
	\draw[{Hooks[round,right]}->] (0:{exp(1/exp(1))-0.02}) to (0:1.02); 
	\draw[decoration={markings,mark=at position 0.7 with {\arrow{>}}},postaction={decorate}] (0:0.98) to [out angle=140, in angle=220,curve through={(0:0.325)}](0:0.98);
	\draw[->>] (0:1.02) to (0:{exp(-1/exp(1))+0.02});  
	
	\draw[blue, fill=blue] (120:{exp(1/exp(1))}) circle (0.02cm);
	\draw[fill=black] (120:1) circle (0.02cm);
	\draw[blue, fill=blue] (120:{exp(-1/exp(1))}) circle (0.02cm);
	
	\draw[{Hooks[round,right]}->] (120:{exp(1/exp(1))-0.02}) to (120:1.02); 
	\draw[decoration={markings,mark=at position 0.7 with {\arrow{>}}},postaction={decorate}] (120:0.98) to [out angle=-100, in angle=-20,curve through={(137:0.6)(120:0.34)(103:0.6)}](120:0.98);
	\draw[->>] (120:1.02) to (120:{exp(-1/exp(1))+0.02});  
	
	\draw[blue, fill=blue] (240:{exp(1/exp(1))}) circle (0.02cm);
	\draw[fill=black] (240:1) circle (0.02cm);
	\draw[blue, fill=blue] (240:{exp(-1/exp(1))}) circle (0.02cm);
	
	\draw[{Hooks[round,right]}->] (240:{exp(1/exp(1))-0.02}) to (240:1.02); 
	\draw[decoration={markings,mark=at position 0.3 with {\arrow{>}}},postaction={decorate}] (240:0.98) to [out angle=20, in angle=100,curve through={(240:0.335)}](240:0.98);
	\draw[->>] (240:1.02) to (240:{exp(-1/exp(1))+0.02}); 
	
	\foreach \th in {0,120,240} \draw (\th:0.5) circle (0.04cm);
	\foreach \n in {1,2,3} \draw ({\n*120}: 1.7) node[draw] {$-1$};
	
	\draw[magenta] ({cos(35)},{sin(35)}) node {$\bullet$};
	\draw[magenta] ({cos(35)+0.1},{sin(35)+0.1}) node {$\widehat{b}$};  
	\foreach \n in {1,2,3} 
	{\draw[magenta] ({120-\n*120}: 0.25) node {$\widehat{S}_{\n}$};
		\draw[magenta, ->] (({120-\n*120}: 0.5)++(-0.13,0) arc (-180:180:0.13);
	}
	\begin{scope}[magenta,decoration={markings,mark=at position 0.55 with {\arrow{>}}}]  
		\draw[postaction={decorate}] ({cos(35)},{sin(35)}) to [out angle=-90, in angle=20,curve through={(15:0.75)}]($({0}:0.5)+(20:0.13)$);
		\draw[postaction={decorate}] ({cos(35)},{sin(35)}) to [out angle=-75, in angle=-40,curve through={(-30:0.7)}]($({-120}:0.5)+(-40:0.13)$);
		\draw[postaction={decorate}] ({cos(35)},{sin(35)}) to [out angle=-60, in angle=-130,curve through={(-60:0.75)(-130:0.8)}]($(-240:0.5)+(-120:0.13)$);
	\end{scope}
\end{tikzpicture}
	\caption{The monodromy data defining the Stokes local system with respect to a basepoint (in magenta) and the deformation data (in black) for the system after Fourier transform.}
	\label{fig:AiryStokesOfFourier}
\end{figure}

Notice that the path associated to the Stokes arrow is closed here, but its inverse image in the picture before Fourier transform is not. The path on the left corresponds to the identity map. The deformation data map associated to this Stokes arrow is therefore given by $-1$ (the transformation rule tells us to change the sign). The computation of the deformation data corresponding to other two Stokes arrows is similar. All three deformation data are indicated in Figure~\ref{fig:AiryStokesOfFourier}.

To explicitly compute the Stokes matrices for Airy, we now need to relate them to these deformation data. To this end, we choose a basepoint $\widehat{b}$, and we get the cycles also shown in Figure~\ref{fig:AiryStokesOfFourier}.

Let us, for instance, detail how to compute $S_2$. We know that it is of the form
$$\widehat S_2=\begin{pmatrix}1&a\\0&1\end{pmatrix}.$$
\newpage Now, the parallel transport along the path we considered above corresponds to the matrix
$$\begin{pmatrix}-1&0\\0&1\end{pmatrix}\widehat S_2\begin{pmatrix}-1&0\\0&1\end{pmatrix}^{-1}$$
(go back to the basepoint via a partial formal monodromy, follow the loop associated with $\widehat S_2$ and go away from the basepoint again),
and the number associated to its deformation datum corresponds to the $(1,2)$-entry of this matrix, so we get the equation
$$-a=-1,$$
and hence
$$\widehat S_2=\begin{pmatrix}1&1\\0&1\end{pmatrix}.$$

Analogously, we get
$$\widehat S_1=\begin{pmatrix}1&0\\-1&1\end{pmatrix}, \qquad \widehat S_3=\begin{pmatrix}1&0\\-1&1\end{pmatrix}.$$
The new formal monodromy is the product of the partial formal monodromies
$$\widehat h = \begin{pmatrix}0&1\\-1&0\end{pmatrix}\begin{pmatrix}1&0\\0&-1\end{pmatrix}\begin{pmatrix}-1&0\\0&1\end{pmatrix}=\begin{pmatrix}0&-1\\1&0\end{pmatrix}.$$

We easily check that
$$\widehat h \widehat S_3\widehat S_2\widehat S_1=\id,$$
so the algorithm provides -- as expected -- a valid set of Stokes data $(\widehat h, \widehat S_1,\widehat S_2,\widehat S_3)\in \mathcal M_B(\widehat\Theta, -1)$ for $\widehat{\Theta}$. 

In this case, the Betti moduli spaces $\mathcal M_B(\Theta, 1)$ and $\mathcal M_B(\widehat\Theta, -1)$ both consist of just a point, so the algorithm induces an isomorphism as expected (and there is no symplectic form to be considered).

\subsection{The case of $\cir{z^{5/3}}$ and $\cir{w^{5/2}}$}

We now discuss an example where one obtains a more interesting isomorphism. Let us consider the irregular class $\Theta=\langle{z^{5/3}}\rangle$ with only one active circle, of slope $5/3$. Its image under Fourier transform is  an irregular class with one active circle (up to a positive real coefficient) $\widehat\Theta=\langle w^{5/2}\rangle$. Let us once again describe the isomorphism between the two corresponding wild character varieties induced by the Fourier transform. This time we will give less details about the passage between the Stokes matrices and the deformation data descriptions, since everything is analogous to the computations seen in the previous two examples.

The Stokes diagram corresponding to $\Theta=\cir{z^{5/3}}$ is drawn in Figure~\ref{fig:stokes_diagram_5_3} below. There are ten singular directions and ten Stokes arrows. 
We choose a reference direction $b\in \partial$ as indicated in the figure. Up to choosing an appropriate framing, one may assume that the formal monodromy has the form
\[
h=\begin{pmatrix}
	0 & 0 & \tau\\
	1 & 0 & 0\\
	0 & 1 & 0
\end{pmatrix}
\]
with $\tau\in \mathbb C^*$. Any conjugacy class $\mathcal C$ for the formal monodromy corresponds to one value of $\tau$. With the choice of numbering of the strands at $b$ indicated in Figure~\ref{fig:stokes_diagram_5_3}, the Stokes matrices are of the form

\[
\begin{array}{cccc}
	S_1=\begin{pmatrix}
		1 & 0 & 0\\
		0 & 1 & 0\\
		s_1 & 0 & 1
	\end{pmatrix},&
	S_2=\begin{pmatrix}
		1 & 0 & 0\\
		0 & 1 & 0\\
		0 & s_2 & 1
	\end{pmatrix},&
	S_3=\begin{pmatrix}
		1 & s_3 & 0\\
		0 & 1 & 0\\
		0 & 0 & 1
	\end{pmatrix},&
	S_4=\begin{pmatrix}
		1 & 0 & t_4\\
		0 & 1 & 0\\
		0 & 0 & 1
	\end{pmatrix},\\
	S_5=\begin{pmatrix}
		1 & 0 & 0\\
		0 & 1 & s_5\\
		0 & 0 & 1
	\end{pmatrix},&
	S_6=\begin{pmatrix}
		1 & 0 & 0\\
		s_6 & 1 & 0\\
		0 & 0 & 1
	\end{pmatrix},&
	S_7=\begin{pmatrix}
		1 & 0 & 0\\
		0 & 1 & 0\\
		s_7 & 0 & 1
	\end{pmatrix},&
	S_8=\begin{pmatrix}
		1 & 0 & 0\\
		0 & 1 & 0\\
		0 & s_8 & 1
	\end{pmatrix},\\
	S_9=\begin{pmatrix}
		1 & s_9 & 0\\
		0 & 1 & 0\\
		0 & 0 & 1
	\end{pmatrix},&
	S_{10}=\begin{pmatrix}
		1 & 0 & s_{10}\\
		0 & 1 & 0\\
		0 & 0 & 1
	\end{pmatrix}.& &
\end{array}
\]
Let us denote by $\Sto_1,\dots, \Sto_{10}$ the corresponding Stokes groups. It follows that the Betti moduli space (taking into account that here $\breve{H}=\mathbb C^*$ acts trivially on $\mathcal E_\Theta$) is explicitly given by 

\[ 
\mathcal M_B(\cir{z^{5/3}},\tau) =\left\lbrace (S_1,\dots, S_{10})\in \prod_{i=1}^{10}\widehat{\Sto}_i \mid h S_{10} \dots S_1=\id\right\rbrace.
\]

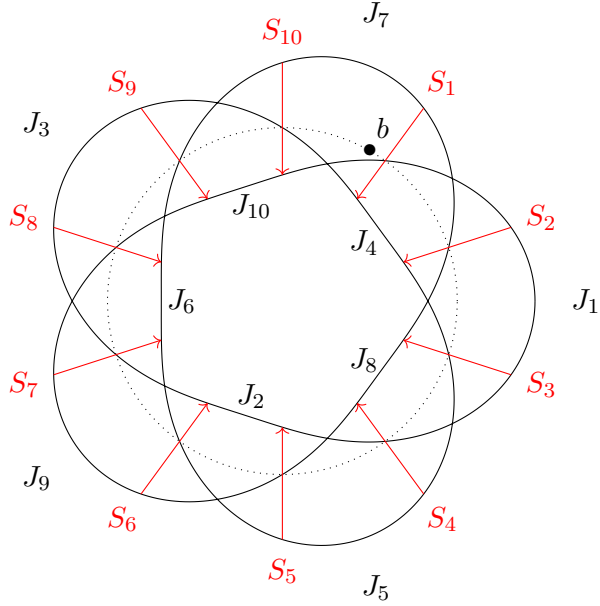
\begin{figure}
	\begin{center}
		\begin{tikzpicture}[scale=2.3]
			\draw[dotted] (0,0) circle (1);
			\draw (60:1) node {$\bullet$};
			\draw (60:1.15) node {$b$};
			\draw[domain=0:(3*360),scale=1,samples=1000] plot (\x:{exp(cos(-5/3*\x)/exp(1))});
			\foreach \th in {1,3,5,7,9}
			\draw ({(-\th+6)*108+180}:{exp(1.25*cos(-(5/3*\th*108+180))/exp(1))}) node {${J}_{\th}$};
            \foreach \th in {2,4,6,8,10}
			\draw ({(-\th+6)*108+180}:{exp(1.5*cos(-(5/3*\th*108+180))/exp(1))}) node {${J}_{\th}$};
			\foreach \th in {0,1,2,3,4,5,6,7,8,9}
			\draw[red,->] ({\th*36+18}:{exp(cos(-(5/3*18))/exp(1))}) to ({\th*36+18}:{exp(cos(-(5/3*90))/exp(1))});
			\foreach \th in {1,2,3,4,5,6,7,8,9,10}
			\draw[red] ({(-\th+2)*36+18}:{exp(1.2*1/exp(1))}) node {$S_{\th}$};

            \draw[teal] (65:1.49) node {\tiny{$1$}};
            \draw[teal] (65:0.95) node {\tiny{$2$}};
            \draw[teal] (65:0.72) node {\tiny{$3$}};
		\end{tikzpicture}
	\end{center}
	\caption{The Stokes diagram of $\cir{z^{5/3}}$, with the distinguished intervals and the Stokes arrows.}
	\label{fig:stokes_diagram_5_3}
\end{figure}

The Stokes diagram for $\cir{w^{5/2}}$ is drawn in Figure~\ref{fig:stokes_diagram_5_2} below. There are five Stokes arrows on five singular directions. We choose a reference direction $\widehat{b}$ as indicated in Figure~\ref{fig:stokes_diagram_5_2}. Again, if we choose an appropriate framing (which is indeed given by the Legendre transform), we can assume that the formal monodromy has the form
\[
\widehat h=\begin{pmatrix}
	0 & \widehat{\tau}\\
	-1 & 0
\end{pmatrix}
\]
with a conjugcacy class for the formal monodromy corresponding to one value of $\widehat \tau$. The Stokes matrices have (with the numbering of strands as shown in Figure~\ref{fig:stokes_diagram_5_2}) the form
\[
\widehat S_1=\begin{pmatrix}
	1 & 0 \\
	t_1 & 1
\end{pmatrix},
\quad
\widehat S_2=\begin{pmatrix}
	1 & t_2\\
	0 & 1
\end{pmatrix},
\quad
\widehat S_3=\begin{pmatrix}
	1 & 0\\
	t_3 & 1
\end{pmatrix},
\quad
\widehat S_4=\begin{pmatrix}
	1 & t_4\\
	0 & 1
\end{pmatrix},
\quad
\widehat S_5=\begin{pmatrix}
	1 & 0\\
	t_5 & 1
\end{pmatrix}.
\]
Let us denote by $\widehat{\Sto}_i$, $i=1,\dots, 5$ the corresponding Stokes groups. The Betti moduli space associated to $\cir{w^{5/2}}$ with formal monodromy determined by $\widehat \tau$ is thus given by

\[ 
\mathcal M_B(\cir{w^{5/2}},\widehat \tau) =\left\lbrace (\widehat S_1,\dots, \widehat S_5)\in \prod_{i=1}^5\Sto_i \mid h\widehat S_5 \widehat S_4 \widehat S_3 \widehat S_2 \widehat S_1=\id\right\rbrace.
\]

\begin{figure}
	\begin{center}
		\begin{tikzpicture}[scale=2.3]
			\draw[dotted] (0,0) circle (1);
			\draw (-120:1) node {$\bullet$};
			\draw (-125:0.87) node {$\widehat{b}$};
			\draw[domain=0:(2*360),scale=1,samples=1000] plot (\x:{exp(cos(-5/2*\x)/exp(1))});
			\foreach \th in {1,2,3,4,5}
			\draw[->,blue] (\th*72:{exp(1/exp(1))}) to ((\th*72:{exp(-1/exp(1))});
			\foreach \th in {1,2,3,4,5}
			\draw[blue] ({(-\th-1)*72}:{0.9*exp(2*1/exp(1))}) node {$\widehat S_{\th}$};
			\foreach \th in {1,3,5,7,9}
			\draw ({(-\th+6)*72}:{exp(1.5*cos(-5/2*\th*72)/exp(1))}) node {$\widehat J_{\th}$};
            \foreach \th in {2,4,6,8,10}
			\draw ({(-\th+6)*72}:{exp(1.25*cos(-5/2*\th*72)/exp(1))}) node {$\widehat J_{\th}$};

            \draw[teal] (-120:1.28) node {\tiny{$1$}};
            \draw[teal] (-120:0.77) node {\tiny{$2$}};
		\end{tikzpicture}
	\end{center}
	\caption{The Stokes diagram of $\cir{w^{5/2}}$, with the distinguished intervals and the Stokes arrows.}
	\label{fig:stokes_diagram_5_2}
\end{figure}
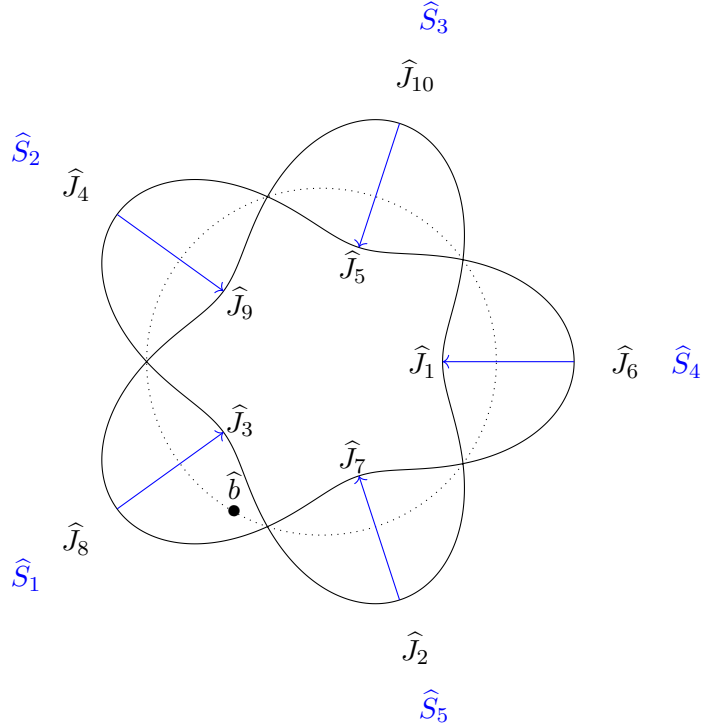

A first observation is the following:

\begin{lemma}
	$\mathcal M_B(\cir{z^{5/2}},\tau)$ is nonempty if and only if $\tau=1$. $\mathcal{M}_B(\cir{w^{5/3}},\widehat \tau)$ is nonempty if and only if $\widehat\tau=1$.
\end{lemma}

\begin{proof}
	For both moduli spaces, considering the determinant of the product of the Stokes matrices immediately implies that the condition $\tau=1$ or $\widehat \tau=1$ is necessary. It is then possible to check directly that in this case the moduli space is nonempty.
\end{proof}

\begin{rem}
	Notice that this condition for the moduli spaces to be nonempty is consistent with the sign change of the formal monodromy in the formal Fourier transform.
\end{rem}

Both circles $\cir{z^{5/3}}$ and $\cir{w^{5/2}}$ are divided into 10 distinguished intervals $J_1,\dots J_{10}$ and $\widehat J_1,\dots, \widehat J_{10}$, respectively, as indicated in Figures~\ref{fig:stokes_diagram_5_3} and \ref{fig:stokes_diagram_5_2}, with the Legendre transform sending $J_i$ to $\widehat{J}_i$. Applying the algorithm, we obtain:

\begin{prop}
	The Fourier transform induces between $\mathcal M_B(\cir{z^{5/2}},1)$ and $\mathcal{M}_B(\cir{w^{5/3}},1)$ the isomorphism $\Phi$ given by 
	
	\[
	\left\lbrace \begin{array}{ll}
		t_1&=-s_6\\
		t_2&=-s_9\\
		t_3&=s_7\\
		t_4&=-s_5\\
		t_5&=-s_8
	\end{array}
	\right.
	\]
	using the choices of basepoints and parametrisation of the Stokes matrices described above. Furthermore, this isomorphism is symplectic.
\end{prop}

\begin{proof}
	The table below summarises the main ingredients of the computation for each coefficient $t_i$: It indicates the Stokes arrow corresponding to the coefficient, the corresponding entry of a product of the initial Stokes matrices obtained by transporting back the Stokes arrows on the initial Stokes diagram, as well as the extra sign coming from the changes of signs in the formal monodromy.
	
	\begin{center}
		\begin{tabular}{|c|c|c|c|}
			\hline
			coefficient & Stokes arrow & Entry of matrix product & extra sign\\
			\hline
			
			$t_1$ & $8\to 3$ & $s_6$ & $-$\\
			$t_2$ & $4 \to 9$ & $-s_9$ & $+$\\
			$t_3$ & $10 \to 5$ & $-s_7$ & $-$\\
			$t_4$ & $6 \to 1$ & $-s_5$ & $+$\\
			$t_5$ & $2 \to 7$ & $\tau s_8$ & $-$\\
			\hline
		\end{tabular}
	\end{center}
	
	The procedure is completely analogous to the one detailed in the previous examples.
	Let us just give some details on the computation of the first coefficient: For $t_1$, the corresponding Stokes arrow is $8\to 3$. To go from $J_8$ to $J_3$ on the initial picture, one needs to consider the product of Stokes factors $S_8\cdots S_4$. Then, taking into account the numbering of the strands (from high to low at $b$), one considers its entry at position $(2,1)$, which is $s_6$. The ``extra sign'' $\varepsilon_{8\to 3}=+$ comes from the signs visible in Theorem~\ref{thm:trafoRule}: The corresponding deformation datum for the Fourier transform will be $-s_6$ in this case, and when reconstructing the Stokes matrix $\widehat{S}_1$ from this deformation datum, one takes into account the signs induced by the Legendre transform on the level of formal monodromy, which in total gives another negative sign. In the case of $t_5$, an additional factor of $\tau$ occurs in this last step, which cancels with the one in the matrix entry. The remaining coefficients are obtained in a similar way. 
	
	We can check explicitly that this indeed provides a well-defined isomorphism. The fact that this isomorphism preserves the symplectic structure can also be checked by direct computation.
\end{proof}

\end{document}